\documentclass{article}
\usepackage[letterpaper,margin=1.0in]{geometry}
\usepackage{hyperref}
\usepackage{enumitem}
\usepackage[utf8]{inputenc} 
\usepackage[T1]{fontenc}    
\usepackage{hyperref}       
\usepackage{url}            
\usepackage{booktabs}       
\usepackage{amsfonts}       
\usepackage{nicefrac}       
\usepackage{microtype}      
\usepackage{graphicx}       
\usepackage[numbers,sort]{natbib}
\usepackage{multirow}
\usepackage{bbding}
\graphicspath{{media/}}     
\usepackage{color}
\usepackage{algorithmic,algorithm}
\usepackage{nicefrac}
\usepackage{cancel}
\usepackage{graphicx} 
\usepackage{booktabs}
\usepackage{makecell}
\usepackage{multirow}
\usepackage{amssymb,amsmath,amsthm}
\usepackage{enumitem}
\usepackage[table]{xcolor}

\usepackage{authblk}
  
\usepackage{amsmath,amsfonts,bm}

\def\eqref#1{equation~\ref{#1}}

\def\1{\bm{1}}

\def\eps{{\epsilon}}

\def\vzero{{\bm{0}}}
\def\vone{{\bm{1}}}

\def\va{{\bm{a}}}
\def\vb{{\bm{b}}}

\def\ve{{\bm{e}}}

\def\vg{{\bm{g}}}

\def\vu{{\bm{u}}}
\def\vv{{\bm{v}}}
\def\vw{{\bm{w}}}
\def\vx{{\bm{x}}}
\def\vy{{\bm{y}}}
\def\vz{{\bm{z}}}

\def\mA{{\bm{A}}}
\def\mB{{\bm{B}}}

\def\mD{{\bm{D}}}

\def\mG{{\bm{G}}}
\def\mH{{\bm{H}}}
\def\mI{{\bm{I}}}
\def\mJ{{\bm{J}}}

\def\mL{{\bm{L}}}

\def\mO{{\bm{O}}}
\def\mP{{\bm{P}}}

\def\mR{{\bm{R}}}
\def\mS{{\bm{S}}}

\def\mU{{\bm{U}}}

\def\mW{{\bm{W}}}

\def\mY{{\bm{Y}}}

\DeclareMathAlphabet{\mathsfit}{\encodingdefault}{\sfdefault}{m}{sl}
\SetMathAlphabet{\mathsfit}{bold}{\encodingdefault}{\sfdefault}{bx}{n}

\def\gA{{\mathcal{A}}}

\def\gC{{\mathcal{C}}}

\def\gF{{\mathcal{F}}}
\def\gG{{\mathcal{G}}}

\def\gO{{\mathcal{O}}}
\def\gP{{\mathcal{P}}}

\def\gX{{\mathcal{X}}}
\def\gY{{\mathcal{Y}}}
\def\gZ{{\mathcal{Z}}}

\def\sB{{\mathbb{B}}}

\def\sN{{\mathbb{N}}}
\def\sO{{\mathbb{O}}}

\def\sR{{\mathbb{R}}}

\newcommand{\E}{\mathbb{E}}

\usepackage{xcolor}
\usepackage{nicefrac}
\usepackage{cancel}
\newtheorem{thm}{Theorem}[section]
\newtheorem{dfn}{Definition}[section]
\newtheorem{exmp}{Example}[section]

\newtheorem{pro}{Property}[section]
\newtheorem{lem}{Lemma}[section]

\usepackage{thmtools}
\usepackage{thm-restate}

\title{On the Condition Number Dependency \\ in 
Bilevel Optimization\textsuperscript{*}}

\author{
  Lesi Chen~$^{\dagger~1~2}$ \qquad
  Kaiyi Ji~$^{\dagger~3}$ \qquad
  Jingzhao Zhang~$^{1~2}$ \\
{\small \texttt{chenlc23@mails.tsinghua.edu.cn, kaiyiji@buffalo.edu, jingzhaoz@mail.tsinghua.edu.cn}} \\ \vspace{2mm}
  $^1$IIIS, Tsinghua University, Beijing, China \\
  $^2$Shanghai Qi Zhi Institute, Shanghai, China \\
  $^3$Department of CSE, University at Buffalo, New York, USA 
}

\begin{document}
\maketitle
\begingroup
\begin{NoHyper}
\renewcommand\thefootnote{*}
\footnotetext{This paper is a merge to two concurrent works of  \citet{chen2025condition} and \citet{ji2025lower} with substantial improvements.}
\end{NoHyper}
\begin{NoHyper}
\renewcommand\thefootnote{$\dagger$}
\footnotetext{The first two authors contribute equally to this work.}
\end{NoHyper}
\endgroup



\begin{abstract}
Bilevel optimization minimizes an objective function, defined by an upper-level problem whose feasible region is the solution of a lower-level problem. We study the oracle complexity of finding an $\epsilon$-stationary point with first-order methods when the upper-level problem is nonconvex, and the lower-level problem is strongly convex. 
Recent works achieve a $\tilde{\mathcal{O}}(\bar \kappa_\vy^{7/2} \epsilon^{-2})$ upper bound that is near-optimal in $\epsilon$.
In this work, we establish a new $\Omega(\kappa_\vy^{5/2} \epsilon^{-2})$ lower bound, where $\kappa_\vy \le \bar \kappa_\vy$ is the lower-level condition number. Our lower bound establishes the first provable gap {in terms of condition number dependency} between bilevel problems and minimax problems in this setup, and \textit{is tight up to logarithmic factors when the lower-level function is quadratic.}
Our lower bounds can be extended to various settings. 
(1) For second-order and arbitrarily smooth problems, we show lower bounds of $\Omega(\kappa_\vy^{9/4} \epsilon^{-7/4})$ and $\Omega(\kappa_\vy^{13/6} \epsilon^{-5/3})$, respectively. 
(2) For convex--strongly-convex problems, we improve the previously best lower bound (Ji and Liang, JMLR 2022) from $\Omega(\kappa_\vy /\sqrt{\epsilon})$ to $\Omega(\kappa_\vy^{3/2} / \sqrt{\epsilon})$. 
(3) For stochastic nonconvex--strongly-convex problems, we also show the lower bounds of $\Omega(\kappa_\vy^4 \epsilon^{-4})$  and $\Omega(\kappa_\vy^{9/2} \epsilon^{-4})$ for stochastic Hessian-vector-product and stochastic first-order methods, respectively.  
\end{abstract}

\section{Introduction}
We study the first-order oracle complexity for solving the bilevel optimization
\begin{align} \label{prob:main-BLO}
    \min_{\vx \in \sR^{d_x}} F(\vx) = f(\vx,\vy^*(\vx)) , \quad \vy^*(\vx) = \arg \min_{\vy \in \sR^{d_y}} g(\vx,\vy).
\end{align}
The goal is to minimize the hyper-objective $F(\vx)$, which is implicitly defined by two explicit objectives $f(\vx,\vy)$ and $g(\vx,\vy)$ and the $\arg \min$ operator in $g(\vx,\vy)$ with respect to~$\vy$. 
This formulation originates from the two-player general-sum Stackelberg game \citep{stackelberg1934marktform}, and reflects the sequential nature of the decision process of two players $\vx$ and $\vy$. Bilevel optimization of this form has received extensive attention in the machine learning community due to its wide applications, such as meta-learning \citep{rajeswaran2019meta}, hyper-parameter tuning \citep{bao2021stability,franceschi2018bilevel,mackay2019self},  generative adversarial networks \citep{goodfellow2020generative}, reinforcement learning \citep{konda1999actor,zeng2024two,hong2023two}, and robotics \citep{wang2024efficient}.


To ensure that Problem (\ref{prob:main-BLO}) is well-defined, we follow \citep{ghadimi2018approximation,ji2021bilevel,hong2023two} to study   \textbf{nonconvex--strongly-convex} (NC--SC) problems, where the lower-level function $g(\vx,\vy)$ is strongly convex in $\vy$ while the upper-level function $f(\vx,\vy)$ can be nonconvex.
In NC--SC bilevel problems, a first-order method applied to the 
 hyper-objective $F(\vx)$ 
 {requires second-order information of $g$}
because 
\begin{align} \label{eq:hyper-grad}
    \nabla F(\vx) = \nabla_\vx f(\vx,\vy^*(\vx)) - \nabla_{\vx \vy}^2 g(\vx,\vy^*(\vx)) [ \nabla_{\vy \vy}^2 g(\vx,\vy^*(\vx))]^{-1} \nabla_\vy f(\vx,\vy^*(\vx)).
\end{align}
In large-scale problems, directly inverting the Hessian matrix $\nabla_{\vy \vy}^2 g(\vx,\vy)$ is very expensive. There are mainly two lines of work to address this issue. 
The first line studies algorithms based on Hessian-vector product (HVP) oracles. Among them, the best-known 
theoretical guarantees \citep{ji2021bilevel,arbel2022amortized} show that they
can find an $\epsilon$-stationary point of $F(\vx)$ with $\gO(\bar \kappa_\vy^{4} \epsilon^{-2})$ and $\gO(\bar \kappa_\vy^9 \epsilon^{-4})$  first-order oracle complexity for smooth deterministic and stochastic problems, respectively, 
where $\bar \kappa_\vy$ is the global condition number for the problem (see Definition~\ref{dfn:cond-number}). In addition,
\citet{yang2023accelerating} showed an improved upper bound of $\tilde \gO(\bar \kappa_\vy^{13/4} \epsilon^{-7/4})$ for second-order smooth deterministic problems.


Since HVP oracles are more expensive than gradient oracles,
another line of work
\citep{liu2020generic,liu2021value,liu2021towards,liu2022bome, kwon2023fully,shen2023penalty,lu2024first,pan2024scalebio}
directly uses gradients to optimize $F(x)$ by solving the penalty function:
\begin{align} \label{eq:Lag}
   \min_{\vx \in \sR^{d_x}, \vy \in \sR^{d_y}}  f(\vx,\vy) + \lambda (g(\vx,\vy) - g(\vx,\vy^*(\vx))),
\end{align}
where $\lambda = \Omega(\epsilon^{-1})$ is a large penalty coefficient.
Among them, 
\citet{kwon2023fully} 
developed the first fully first-order method that can find an $\epsilon$-stationary point of $F(\vx)$ with $\tilde \gO({\rm poly}(\bar \kappa_\vy) \epsilon^{-3})$ first-order oracles for smooth deterministic problems and $\tilde \gO({\rm poly}(\bar \kappa_\vy) \eps^{-7} )$ for stochastic problems.
In a follow-up work, \citet{chen2023near} improved the complexity to $\tilde \gO({\bar \kappa_\vy^{4}} \epsilon^{-2})$ and $\tilde \gO({\bar \kappa_\vy^{11}} \epsilon^{-6})$ for deterministic and stochastic problems, respectively. \citet{chen2023near} also extended their method to achieve the $\tilde \gO({\bar \kappa_\vy^{15/4}} \epsilon^{-7/4})$ complexity for second-order smooth deterministic problems . Moreover, as we additionally show in Appendix~\ref{apx:ub-refined}, by using a naive application of Nesterov
acceleration in the inner loop, the deterministic method by \citet{chen2023near} can be improved to $\tilde \gO(\bar \kappa_\vy^{7/2} \epsilon^{-2}) $ and $\tilde \gO(\bar \kappa_\vy^{13/4} \epsilon^{-7/4})$ for first- and second-order smooth problems, respectively. 

Although both HVP- and gradient-based methods can achieve (near)-optimal $\epsilon^{-2}$ and $\epsilon^{-7/4}$ rates for first- and second-order smooth NC--SC bilevel problems, their upper bounds have high~$\bar \kappa_\vy$~dependencies compared with 
NC--SC minimax problems, which are special bilevel problems with $f = -g$. For instance, 
\citet{lin2020near,wang2024efficient} proposed algorithms with a $\tilde \gO(\sqrt{\kappa_\vy} \epsilon^{-2})$ and $\tilde \gO(\sqrt{\kappa_\vy} \eps^{-7/4})$ complexity for first- and second-order smooth NC--SC 
deterministic minimax problems, respectively, where $\kappa_\vy \le \bar \kappa_\vy$ is the condition number of the lower-level problem (see Definition \ref{dfn:cond-number}). Moreover, nearly matching lower bounds of  $\Omega(\sqrt{\kappa_\vy} \epsilon^{-2})$ and $\Omega(\sqrt{\kappa_\vy} \eps^{-7/4})$ are shown by \citet{li2021complexity,zhang2021complexity} and \citet{wang2024efficient}, respectively. However, the state-of-the-art methods 
for deterministic smooth NC--SC bilevel problems require an $\tilde \gO(\bar \kappa_\vy^{7/2} \epsilon^{-2})$ complexity.
Given this, we wonder  \textit{whether NC--SC bilevel problems are provably more challenging than NC--SC minimax problems measured by the condition number dependency.}

To answer this question, we establish new lower bounds for NC--SC bilevel problems, where almost all previously applicable lower bounds are attained in minimax ({not bilevel}) problems where~$f = -g$. 
In this work, by leveraging the asymmetry of the objectives $f$ and $g$, we modify the nested chain construction originally designed for minimax problems \citep{li2021complexity} to create a longer effective chain in $\vx$, leading to larger $\kappa_\vy$ dependencies in various setups:
\vspace{-1mm}
\begin{enumerate}
    \item For deterministic NC--SC bilevel problems, we show a lower bound of $\Omega(\kappa_\vy^{5/2} \epsilon^{-2})$, improving the lower bound of $\Omega(\sqrt{\kappa_\vy} \epsilon^{-2})$ from NC--SC minimax problems \citep{li2021complexity,zhang2021complexity} by a factor of $\kappa_\vy^2$.
\vspace{-1mm}
\item Motivated by the very inspiring work by \citet{zhou2026sharp}, we also provide lower bounds under additional $p$th-order smoothness. Specifically, we establish lower bounds of $\Omega(\kappa_\vy^{9/4}\epsilon^{-7/4})$ for $p=2$ and  $\Omega(\kappa_\vy^{13/6}\epsilon^{-5/3})$ for $p \ge 3$, respectively. 
\vspace{-1mm}
\item For convex--strongly-convex (C--SC) and strongly-convex--strongly-convex (SC--SC) bilevel problems, we also show lower bounds of 
$\Omega(\kappa_\vy^{3/2} / \sqrt{\epsilon})$ and $\Omega(\kappa_\vy^{3/2} \sqrt{\kappa_\vx})$, respectively, improving 
prior results \citep{ji2023lower}
of $\tilde \Omega(\min \{ \kappa_\vy, \epsilon^{-3/2} \} / \sqrt{\epsilon} )$ and~$\tilde \Omega(\kappa_\vy \sqrt{\kappa_\vx})$ by a factor of $\sqrt{\kappa_\vy}$.
\vspace{-1mm}
\item For stochastic NC--SC bilevel problems, we prove a lower bound of 
$\Omega(\kappa_\vy^4 \epsilon^{-4})$ and $\Omega(\kappa_\vy^{9/2} \epsilon^{-4})$ for stochastic HVP and first-order methods, respectively, improving the lower bound of $\Omega(\kappa_\vy^{1/3} \epsilon^{-4})$ 
from minimax problems \citep{li2021complexity}.
\end{enumerate}

To summarize, we compare our lower bounds with existing upper bounds under the above three setups in Tables \ref{tab:det}-\ref{tab:stoc-NCSC}. Note that we prove our lower bounds on lower-level quadratic problems, which form a simple yet representative setup.
In the nonconvex--quadratic (NC--Q) setting, we show that existing upper bounds can also be refined to 
$\tilde \gO(\kappa_\vy^{5/2} \epsilon^{-2})$ and $\tilde \gO(\bar \kappa_\vy^{9/4} \epsilon^{-7/4})$ for first- and second-order smooth deterministic problems, respectively, showing that our lower bounds are tight up to logarithmic factors. Similarly, our lower bounds $\Omega(\kappa_\vy^{3/2} / \sqrt{\epsilon})$ and $\Omega(\kappa_\vy^{3/2} \sqrt{\kappa_\vx})$ for convex--quadratic (C--Q) and strongly-convex--quadratic (SC--Q) problems are also nearly achievable by existing methods \citep{ji2023lower}. Consequently, our results yield a near-tight characterization of the complexity of deterministic bilevel problems with quadratic lower-level functions. However, the optimality in more general setups remains an open problem, forming a sharp contrast to minimax optimization, where lower bounds for SC--Q and NC--Q problems \citep{li2021complexity,zhang2021complexity,zhang2022lower} are also tight for general SC--SC and NC--SC problems.

\begin{table*}[htbp] 
    \centering
     \caption{We present the oracle complexities
      of \textbf{deterministic} methods for finding an $\epsilon$-stationary point for 
     smooth nonconvex--strongly-convex (NC--SC) and nonconvex--quadratic (NC--Q) bilevel optimization under first-order or HVP oracles.
     }
    \label{tab:det}
      \begin{tabular}{c c c}
    \hline 
     Setup  &  Reference & Complexity \\ 
    \hline \addlinespace 
    Upper Bound (NC--SC) & \citet{ji2021bilevel,chen2023near}  & $\tilde \gO(\bar \kappa_\vy^{7/2}\epsilon^{-2})$  \\ \addlinespace
    Upper Bound (NC--Q) & \citet{ji2021bilevel,chen2023near} & $\widetilde{\mathcal{O}}(\kappa_\vy^{5/2}\epsilon^{-2})$  \\ \addlinespace
    \hline \addlinespace
Lower Bound (NC--Q) 
     & This Paper (Theorem \ref{thm:NC-SC}) & { $\Omega(\kappa_\vy^{5/2} \epsilon^{-2})$}  \\ \addlinespace
    \hline 
    \end{tabular}
\end{table*}

\begin{table*}[htbp] 
    \centering
     \caption{We present the oracle complexities
      of \textbf{deterministic} methods for finding an $\epsilon$-stationary point for 
     \textbf{second-order smooth} NC--SC and NC--Q bilevel optimization under first-order or HVP oracles.
     }
    \label{tab:det-second-order}
      \begin{tabular}{c c c}
    \hline 
     Setup  &  Reference & Complexity \\ 
    \hline \addlinespace 
    Upper Bound (NC--SC) & \citet{yang2023accelerating,chen2023near}  & $\tilde \gO(\bar \kappa_\vy^{13/4}\epsilon^{-7/4})$  \\ \addlinespace
    Upper Bound (NC--Q) & \citet{yang2023accelerating,chen2023near} & $\widetilde{\mathcal{O}}(\kappa_\vy^{9/4}\epsilon^{-7/4})$  \\ \addlinespace
    \hline \addlinespace
Lower Bound (NC--Q) 
     & This Paper (Theorem \ref{thm:highly-smooth-NC-Q}) & { $\Omega(\kappa_\vy^{9/4} \epsilon^{-7/4})$}  \\ \addlinespace
    \hline 
    \end{tabular}
\end{table*}

\begin{table*}[htbp] 
    \centering
     \caption{We present the oracle complexities
      of \textbf{deterministic} methods for finding an $\epsilon$-stationary point for 
     smooth convex--strongly-convex (C--SC) and convex--quadratic (C--Q) bilevel optimization under first-order or HVP oracles.
     }
    \label{tab:det-c-sc}
      \begin{tabular}{c c c}
    \hline 
     Setup  &  Reference & Complexity \\ 
    \hline \addlinespace 
    Upper Bound (C--SC) & \citet{ji2023lower}  & $\tilde \gO(\bar \kappa_\vy^{2}\epsilon^{-1/2})$  \\ \addlinespace
    Upper Bound (C--Q) & \citet{ji2023lower} & $\widetilde{\mathcal{O}}(\kappa_\vy^{3/2}\epsilon^{-1/2})$  \\ \addlinespace
    \hline \addlinespace
   Lower Bound (C--Q)
    & \citet{ji2023lower} & $\Omega(\kappa_\vy^{5/4} \eps^{-1/2})$ \\ \addlinespace
    Lower Bound (C--Q)  & This Paper (Theorem \ref{thm:C-SC}) & { $\Omega(\kappa_\vy^{3/2} \epsilon^{-1/2})$} \\ \addlinespace
    \hline 
    \end{tabular}
\end{table*}
\begin{table*}[htbp] 
    \centering
     \caption{We present the oracle complexities
      of \textbf{stochastic} methods for finding an $\epsilon$-stationary point for 
     smooth nonconvex--strongly-convex (NC--SC) bilevel optimization under stochastic HVP and first-order oracles.}
     \vspace{0.2cm}
    \label{tab:stoc-NCSC}
    \begin{tabular}{c c c c }
    \hline 
    Oracle & Setup &  Reference & Complexity  \\ 
    \hline \addlinespace 
    \multirow{3}{*}{HVP} & Upper Bound & \citet{ji2021bilevel} & 
    $\tilde \gO(\bar \kappa_\vy^9 \epsilon^{-4})$ \\ \addlinespace 
    & Lower Bound  & This Paper (Theorem \ref{thm:NC-SC-stoc})  & $\Omega(\kappa_\vy^4 \epsilon^{-4})$ \\ \addlinespace
    \hline  \addlinespace
    \multirow{5}{*}{First-Order} & Upper Bound  & 
    \citet{kwon2023fully}
    & $\tilde{\mathcal{O}}({\rm poly}(\bar \kappa_\vy) \epsilon^{-7})$ \\ \addlinespace
    & Upper Bound  & 
    \citet{chen2023near}
    & $\tilde{\mathcal{O}}(\bar \kappa_\vy^{11}\epsilon^{-6})$ \\ \addlinespace
    \cline{2-4} \addlinespace
    & Lower Bound  & This Paper (Theorem \ref{thm:stoc-NC-Q}) & { $\Omega(\kappa_\vy^{9/2} \epsilon^{-4})$}  \\  \addlinespace
    \hline  
    \end{tabular}
\end{table*}

\paragraph{Notations.} For a vector $\vx \in \sR^d$, 
 we let $x_i $ be its $i$th coordinate; define its support by ${\rm supp}(\vx) = \{ i \mid \vx_i \ne 0 \}$ and its progress by ${\rm prog}(\vx) = \max\{ i \ge 0 \mid \vx_i \ne 0 \}$.
 We use $\|\vx\| = \sqrt{\sum_{i=1}^d x_i^2}$ and 
$\|\vx\|_\infty = \max_{1 \le i \le d} |x_i|$ to denote the $\ell_2$ and $\ell_\infty$ norms, respectively. Moreover, the notation $\Vert \,\cdot \Vert$ also denotes the operator norm of a matrix. Given an event $A$, we let $\vone_A$ denote its indicator function. We let $\ve_i$ be the vector with its $i$th element being $1$ and all the other elements being $0$. For a convex compact set $\gX \subseteq \sR^d$, we denote $\gP_\gX$ as the projection operator onto $\gX$.
 Given functions $f,g : \gX \rightarrow [0,+\infty)$, we use $f = \gO(g)$ to denote there exists a numerical constant $C< +\infty$
such that $f(x) \le C g(x)$ for all $x \in \gX$,  $f = \Omega(g)$ to denote $g = \gO(f)$,  and $f = \tilde \gO(g)$ to denote $f = \gO(g \max\{1, \log g \})$. We also denote $f \asymp g$ when $f = \gO(g)$ and $f = \Omega(g)$ both hold, and denote $\Omega_p$ to hide the constants that only depend on $p$. 
\subsection{Additional Related Work}
In addition to the works discussed above, we further review the closely related work on lower bounds, including lower bounds for minimization, minimax, and bilevel problems.

\paragraph{Lower bounds for minimization problems.} In the textbooks, \citet{nemirovskij1983problem,nesterov2018lectures} established tight lower bounds for convex optimization using chain-like quadratic functions. Recently, \citet{carmon2020lower,carmon2021lower} further abstracted the above idea into the concept of zero-chain as a generic approach for proving lower bounds. By carefully constructing zero-chains, \citet{carmon2020lower} proved the tight $\Omega(\epsilon^{-p/(p+1)})$  lower bound for 
$p$th-order methods to find $\epsilon$-stationary points of $p$th-order smooth nonconvex functions. In the same series of work, \citet{carmon2021lower} proved $\Omega(\epsilon^{-12/7})$ and  $\Omega(\epsilon^{-8/5})$  lower bounds for 
first-order methods to optimize second-order smooth and arbitrarily smooth nonconvex functions, respectively.
Very recently, in an inspiring work, \citet{zhou2026sharp} introduced the concept of block-chain and improved the lower bounds to 
$\Omega(\epsilon^{-7/4})$ and  $\Omega(\epsilon^{-5/3})$, which are achieved by the upper bound in \citet{carmon2017convex} up to logarithmic factors.
For nonconvex stochastic first-order optimization,
\citet{arjevani2023lower} proved the tight $\Omega(\epsilon^{-4})$ and $\Omega(\epsilon^{-3})$ lower bounds under the general expectation and the mean-squared smooth setting, respectively; 
\citet{fang2018spider,li2021page} proved the tight $\Omega(N + \sqrt{N} \epsilon^{-2})$ lower bound for nonconvex finite-sum first-order optimization; \citet{arjevani2020second} proved 
the tight $\Omega(\epsilon^{-7/2})$ and $\Omega(\epsilon^{-3})$ lower bounds for first-order stochastic nonconvex optimization on second-order smooth objectives, and
high-order stochastic nonconvex optimization on arbitrarily smooth objectives, respectively.

\paragraph{Lower bounds for minimax problems.} Inspired by the lower bounds in minimization problems, researchers also proved tight lower bounds for (first-order) minimax optimization under many setups. 
\citet{nemirovski2004prox} gave tight $\Omega(\epsilon^{-1})$ lower bound for convex--concave minimax problems. \citet{ouyang2021lower,zhang2022lower} proved the tight 
 $\Omega(\sqrt{\kappa_\vy/ \epsilon})$ and
$\tilde \Omega(\sqrt{\kappa_\vx \kappa_\vy})$
 lower bounds for C--SC and SC--SC problems, respectively.
\citet{han2024lower} extended those results to the finite-sum setting.
\citet{li2021complexity,zhang2021complexity} proved the tight $\Omega(\sqrt{\kappa_\vy} \epsilon^{-2})$ lower bound for nonconvex--strongly-convex (NC--SC) problems.
Their results have also been extended to
$\Omega(\sqrt{\kappa_\vy} \epsilon^{-12/7})$, $\Omega(N + \sqrt{N \kappa_\vy} \epsilon^{-2})$, and  
$\Omega(\kappa_\vy^{1/3} \epsilon^{-4})$ for 
second-order smooth problems \citep{wang2024efficient},
finite-sum problems \citep{zhang2021complexity}, and 
general stochastic problems \citep{li2021complexity}, respectively.

\paragraph{Lower bounds for bilevel problems.} Although the lower bounds for both minimization and minimax problems have been extensively studied, we find there are only a few works on the lower complexity bounds of bilevel optimization, except for the following ones. In their novel work,
\citet{ji2023lower} proved lower bounds of $\tilde \Omega(\min \{ \kappa_\vy, \epsilon^{-3/2} \} / \sqrt{\epsilon}  )$ and 
$\tilde \Omega (\kappa_\vy \sqrt{\kappa_\vx})$ for C--SC and SC--SC bilevel problems, respectively, but they left the lower bounds for nonconvex bilevel problems unexplored.
\citet{kwon2024complexity} established tight upper and lower bounds for $r$-$\vy^*(\vx)$-aware optimization, 
where the oracle only gives  $r$-estimates of $\vy^*(\vx)$ along with $r$-locally reliable gradients. This setting recovers stochastic bilevel optimization if $r = +\infty$, but requires $r = 
\gO(\epsilon)$, which is different from the setting of bilevel optimization with standard oracles.
\citet{chen2024finding} considered bilevel problems without lower-level strong convexity, and constructed a hard instance such that any first-order algorithm will get stuck at $\vx_0$ within finite iterations, which can be viewed as the limit of $\kappa_\vy \uparrow +\infty$ for the results in this paper but with an additional constraint $\Vert \vy^*(\vx) \Vert \le D$.

\section{Problem Setup}

In this section, we introduce the bilevel function class considered in this paper, the algorithm class to which our lower bound applies, and the complexity measurement of the algorithms.

\subsection{Function Class}

Our main focus is the general high-order smooth nonconvex--strongly-convex (NC--SC) bilevel problems as follows, where the case $p=1$ is introduced by \citet{ghadimi2018approximation}, and the case $p=2$ is recently introduced by \citet{huang2025efficiently}.

\begin{dfn}[$p$th-order smooth NC--SC  problem] \label{dfn:NC-SC-func}
Given $p \in \sN_+$, $L_0, \cdots,L_{p+1} >0$, $\Delta>0$, and $\mu_\vy \in (0,L_1]$, 
we use $\gF^{\text{nc--sc}}(L_0,\cdots,L_{p+1},\mu_\vy, \Delta)$ to denote the set of all bilevel problems $(f,g)$, where the upper-level function $f:\sR^{d_x} \times \sR^{d_y} \rightarrow \sR$ and lower-level function $g: \sR^{d_x} \times \sR^{d_y} \rightarrow \sR$ satisfy the following assumptions:
\begin{enumerate}
    \item $f$ has $L_q$-Lipschitz continuous $q$th-order derivative for all $q = 0,\cdots,p$; and $g$ has $L_q$-Lipschitz continuous $q$th-order derivative for all $q = 1,\cdots,p+1$, that is, for every $\vx_1, \vx_2 \in \sR^{d_x}$ and $\vy_1,\vy_2 \in \sR^{d_y}$, we have
    \begin{align*}
        \Vert \nabla^q f(\vx_1,\vy_1) - \nabla^q f(\vx_2,\vy_2) \Vert &\le L_q \Vert (\vx_1,\vy_1) - (\vx_2,\vy_2) \Vert, \quad q=0,\cdots,p~; \\
        \Vert \nabla^q g(\vx_1,\vy_1) - \nabla^q g(\vx_2,\vy_2) \Vert &\le L_q \Vert (\vx_1,\vy_1) - (\vx_2,\vy_2) \Vert, \quad q=1,\cdots,p+1.
    \end{align*}
    \item For any fixed $\vx \in \sR^{d_x}$, $g$ is $\mu_\vy$-strongly convex in $\vy$, that is, for any $\vy_1,\vy_2 \in \sR^{d_y}$,
    \begin{align*}
        g(\vx,\vy_2) \ge g(\vx,\vy_1) + \langle \nabla_\vy g(\vx,\vy_1), \vy_2 - \vy_1 \rangle + \frac{\mu_\vy}{2} \Vert \vy_1 - \vy_2 \Vert^2.
    \end{align*}
    \item The hyper-objective $F(\vx) = f(\vx,\vy^*(\vx))$ in Eq. (\ref{prob:main-BLO}) is lower bounded and 
    \begin{align*}
        F(\vzero) - \inf_{\vx \in \sR^{d_x}} F(\vx) \le \Delta.
    \end{align*}
\end{enumerate}
\end{dfn}

For these problems, we define the lower-level and global condition numbers $\kappa_\vy$ and $\bar \kappa_\vy$
as follows. By definition, we have $\kappa_\vy \le \bar \kappa_\vy$. When $L_1 \asymp \bar L$ we know that $\kappa_\vy \asymp \bar \kappa_\vy$ . Otherwise, these two condition numbers can have very
different magnitudes.


\begin{dfn} \label{dfn:cond-number}
    For the $p$th-order smooth NC--SC function class in Definition \ref{dfn:NC-SC-func}, we define $\kappa_\vy: =L_1/\mu_\vy$ as the condition number of $g$  and $\bar \kappa_\vy: = \bar L / \mu_\vy$ as the global condition number of the bilevel problem, where $\bar L = \max_{q = 0,\cdots,p+1} L_q$.
\end{dfn}


In our lower bound analyses, we construct bilevel problems with quadratic lower-level functions, which are formally defined as the following nonconvex--quadratic (NC--Q) problems as follows. Note that for this subclass, we do not make assumptions on the boundedness of $L_0$ since neither the lower nor the upper bounds depend on it.

\begin{dfn}[$p$th-order smooth NC--Q problem] \label{dfn:qg-sub-class}
For $p \in \sN_+$,~$L_1,\cdots,L_p>0$, $\mu_\vy \in (0,L_1]$, and $\Delta>0$, we denote
$\gF^{\text{nc--q}}(L_1,\cdots,L_{p},\mu_\vy, \Delta) \subseteq \gF^{\text{nc--sc}}(+\infty, L_1,\cdots,L_p,0,\mu_\vy, \Delta)$
as the subset of NC--SC problems, where the lower-level function takes the form of
\begin{align} \label{eq:quadratic-g}
    g(\vx,\vy) = \frac{1}{2} \vy^\top \mH \vy+ \vx^\top \mJ \vy +  \vb^\top \vy,
\end{align}
where $\mu_\vy \mI_{d_y} \preceq \mH \preceq L_1 \mI_{d_y}$ and $-L_1 \mI_{d_y} \preceq \mJ \preceq L_1 \mI_{d_y}$.
\end{dfn}

From a complexity theory perspective, it is natural to consider quadratic functions as optimizing convex quadratics is usually as hard as general convex functions \citep[Theorem 2.1.13]{nesterov2018lectures}. Therefore, 
understanding the complexity of quadratic problems is essential to tackle more complex setups.
From a practical perspective, this problem class already covers many useful applications in reinforcement learning, meta-learning, robust optimization, and graph neural network training, as detailed in Appendix \ref{apx:example-Q}.

\subsection{Algorithm Class}

Following the seminal framework proposed by \citet{nemirovskij1983problem}, we model the algorithmic complexity by the number of oracle accesses, where first-order oracles reveal the gradient information at the query points, defined as follows.

\begin{dfn}[First-order oracle] 
Let $f$ and $g$ be two continuously differentiable functions. We define the first-order oracle 
for the bilevel problem $(f,g)$
as the mapping 
\[
\sO^{\rm fo}(\vx,\vy) = (f(\vx,\vy), g(\vx,\vy), \nabla f(\vx,\vy),  \nabla g(\vx,\vy)).
\]
\end{dfn}





We then define the (deterministic) first-order algorithms as follows, which 
generalizes the algorithm classes defined 
in \citep{agarwal2018lower,arjevani2020second,carmon2020lower,arjevani2023lower}
from a single objective to two objectives.
Our definition  
encompasses almost all existing first-order algorithms for bilevel optimization \citep{liu2022bome,shen2023penalty,kwon2023fully,lu2024first}.

\begin{dfn}[First-order algorithm] \label{dfn:rand-alg}
A first-order algorithm $\texttt{A}$ consists of a sequence of measurable mappings $\{ \texttt{A}^{(t)}\}_{t \in \sN}$ such that $ \texttt{A}^{(t)}$ takes in the first $t-1$ oracle responses to produce the $t$-th query $\vz^{(t)} = (\vx^{(t)},\vy^{(t)})$, such that $\vz^{(0)} = (\vzero,\vzero)$ and 
\begin{align}  \label{eq:rand-seq}
    \vz^{(t)} = \texttt{A}^{(t)} \left(\sO^{\rm fo}(\vz^{(0)}), \cdots, \sO^{\rm fo}(\vz^{(t-1)}) \right), \quad \forall t \in \sN_+.
\end{align}                                                                             
We denote $\gA^{{\rm fo}}$ as the set of all deterministic algorithms that follow protocol (\ref{eq:rand-seq}).
\end{dfn}

In addition to first-order algorithms, all of our lower bounds also apply to the 
following algorithm class in Definition \ref{dfn:rand-alg-hvp}, which includes many standard HVP-based algorithms \citep{ghadimi2018approximation,ji2021bilevel,hong2023two}. 
These algorithms introduce an additional variable $\vv$ that tracks $-(\nabla_{\vy\vy}^2 g(\vx,\vy))^{-1} \nabla_\vy g(\vx,\vy)$ via gradient updates on the quadratic function $\vv^\top \nabla_{\vy\vy}^2 g(\vx,\vy) \vv - 2\nabla_\vy f(\vx,\vy)^\top \vv$ and then use the approximate hyper-gradient $\nabla_\vx f(\vx,\vy) + \nabla_{\vx\vy}^2 g(\vx,\vy) \vv$ to update $\vx$. Formally, we define the deterministic HVP oracle and HVP-based methods as follows.




\begin{dfn}[HVP oracle]  \label{dfn:hvp-oracle}
Let $f$ be differentiable, and $g$ be twice differentiable. We define the HVP oracle for the bilevel problem $(f,g)$ as the mapping
\begin{align*}
 \sO^{\rm hvp}(\vx,\vy,\vv) = \left({\sO^{\rm fo}(\vx,\vy)}, \nabla_{\vx \vy}^2 g(\vx,\vy) \vv, \nabla_{\vy\vy}^2 g(\vx,\vy) \vv \right).
\end{align*}
\end{dfn}

\begin{dfn}[HVP-based algorithm] \label{dfn:rand-alg-hvp}
An HVP-based algorithm $\texttt{A}$ consists of a sequence of measurable mappings $\{ \texttt{A}^{(t)}\}_{t \in \sN}$ such that $ \texttt{A}^{(t)}$ takes in the first $t-1$ oracle responses to produce the $t$-th query $\vw^t = (\vx^{(t)},\vy^{(t)}, \vv^{(t)})$, such that $\vw^{(0)} = (\vzero,\vzero, \vzero)$ and  
\begin{align}  \label{eq:rand-seq-hvp}
    \vw^{(t)} = \texttt{A}^{(t)} \left(\sO^{\rm hvp}(\vw^{(0)}), \cdots, \sO^{\rm hvp}(\vw^{(t-1)}) \right), \quad \forall t \in \sN_+.
\end{align}                                                                             
We denote $\gA^{{\rm hvp}}$ as the set of all deterministic algorithms that follow protocol (\ref{eq:rand-seq-hvp}).

\end{dfn}


Definition \ref{dfn:rand-alg-hvp} subsumes the hyper-gradient-based algorithm class previously defined in \citep[Definition 2]{ji2023lower}.
{Since the HVP oracles for quadratic functions are equivalent to the differences of gradients,} proving lower bounds for first-order algorithms on lower-level quadratic problems implies that the same lower bound applies to HVP-based algorithms. Formally, we have the following statement.



\begin{lem} \label{lem:zr-grad-hvp}
If $g$ is a quadratic function of the form (\ref{eq:quadratic-g}), then
$  \gA^{{\rm hvp}} = \gA^{\rm fo}$.
\end{lem}

\begin{proof} 
The inclusion $\gA^{\rm fo} \subseteq \gA^{\rm hvp}$ is trivial since the HVP oracle includes the first-order oracle. In the following, let us show that the reverse inclusion $ \gA^{\rm hvp} \subseteq \gA^{\rm fo}$ also holds and conclude that $ \gA^{\rm hvp}= \gA^{\rm fo}$. The intuition is that, if $g$ is a quadratic function of the form~(\ref{eq:quadratic-g}), then we can simulate the HVP oracles using the difference of first-order oracles. Note that  $\nabla_\vy g(\vzero,\vzero) = \vb$, $\nabla_{\vy\vy}^2 g(\vx,\vy) =  \mH$, and $\nabla_{\vx \vy}^2 g(\vx,\vy) = \mJ$. Therefore, 
it holds that
\begin{align} \label{eq:rela-gq}
\begin{split}
\begin{cases}
    \nabla_{\vy \vy}^2 g(\vx,\vy) \vv  = \mH \vv = {\nabla_\vy g(\vzero,\vv) - \nabla_\vy g(\vzero,\vzero)}, \\
    \nabla_{\vx \vy}^2 g(\vx,\vy) \vv = \mJ \vv = \nabla_\vx g(\vx,\vv). 
\end{cases}    
\end{split}
\end{align}
{By Definition \ref{dfn:hvp-oracle}, the first two components of the HVP oracle are directly provided by the first-order oracles, and the last two components (the Hessian-vector products) are linear combinations of first-order oracle responses at specific queries $(\vzero, \vv)$ and $(\vx, \vv)$. 

Now, for any sequence $\{ (\vx^{(s)},\vy^{(s)},\vv^{(s)}) \}_{s=0}^t$ generated by an HVP-based algorithm that satisfies Eq. (\ref{eq:rand-seq-hvp}), we define the expanded sequence 
$\{ (\vx^{(s)},\vy^{(s)}), (\vx^{(s)},\vv^{(s)}), (\vzero,\vv^{(s)}) \}_{s=0}^t$ and prove by induction that it satisfies Eq. (\ref{eq:rand-seq}). Note that the induction base ($t=0$) is trivial since we have $\vx^{(0)} = \vy^{(0)} = \vv^{(0)} = \vzero$. 

In addition, assuming that $\{ (\vx^{(s)},\vy^{(s)}), (\vx^{(s)},\vv^{(s)}), (\vzero,\vv^{(s)}) \}_{s=0}^{t}$ satisfies Eq.~(\ref{eq:rand-seq}), then we can know that the new queries \((\vx^{(t+1)},\vy^{(t+1)}), (\vx^{(t+1)},\vv^{(t+1)}), (\vzero,\vv^{(t+1)})\) are a deterministic mapping of previous HVP oracles, which by Eq. (\ref{eq:rela-gq}) is also a deterministic mapping of previous first-order oracles.
Hence, the induction follows, and we can finally obtain that $\gA^{\rm hvp} \subseteq \gA^{\rm fo} $.}
\end{proof}



\subsection{Complexity Measure}
The ultimate goal of bilevel optimization is to optimize the hyper-objective $F(\vx)$ implicitly defined in Eq. (\ref{prob:main-BLO}). Given that $F(\vx)$ is nonconvex in general, finding its approximate minima requires exponential complexity in the worst case \citep{nemirovskij1983problem}. Hence, we follow \citep{carmon2020lower,carmon2021lower,arjevani2023lower} to measure the complexity of an algorithm by the number of oracle calls to find an $\epsilon$-stationary point.

\begin{dfn} \label{dfn:fo-sta}
We call 
a point $\vx \in \sR^{d_x}$ an $\epsilon$-stationary point of $F: \sR^{d_x} \rightarrow \sR$ if 
\begin{align*}
    \Vert \nabla F(\vx) \Vert \le \epsilon.
\end{align*}
\end{dfn}

\begin{dfn}
Given a bilevel function class $\gF$  and the associated first-order oracle $\sO$, we define the worst-case complexity of an algorithm class $\gA$ to find an $\epsilon$-stationary point as
\begin{align*}
    {\rm Compl}_{\epsilon}(\gA, \gF, \sO) = \sup_{(f,g) \in \gF} \inf_{\texttt{A} \in \gA} \inf \{t \in \sN \mid \Vert \nabla F(\vx^{(t)}) \Vert \le \epsilon \},
\end{align*}
where $\{\vz^{(t)}\}_{t \in \sN}$ with $\vz^{(t)} = (\vx^{(t)},\vy^{(t)})$ is the sequence generated by \texttt{A} running on $(f,g)$.
\end{dfn}

\section{Lower Bound Construction} \label{sec:lb}

In this section, we introduce the framework to prove lower bounds. In Section \ref{subsec:zo-chain}, we first recall the technique of zero-chains \citep{nesterov2018lectures,arjevani2019oracle,carmon2020lower,carmon2021lower}, which are commonly used in proving lower bounds in single-level convex or nonconvex optimization. Then, in Section \ref{subsec:nested-chain}, we introduce the nested chains to characterize the complexities of bilevel problems. Finally, in Section \ref{subsec:hs-nested-chain}, we model highly-smooth nonconvex optimization as an implicit bilevel problem and apply the idea of nested chains to prove tight lower bounds  \citep{carmon2021lower,zhou2026sharp}, which are used to show lower bounds for highly smooth NC--SC bilevel optimization in the later sections.

\subsection{Zero-Chain and Zero-Respecting Algorithm}
\label{subsec:zo-chain}

We construct hard instances with chain-like structures to prove optimization lower bounds. {Such functions are known as zero-chains, which are 
well-known and applied extensively in the literature on lower
bounds for optimization \citep{nemirovskij1983problem,nesterov2018lectures,woodworth2016tight,carmon2020lower,carmon2021lower,arjevani2019oracle,arjevani2023lower}. We recall the formal definition \citep[Definition 3]{carmon2020lower} as follows.}
\begin{dfn}[Zero-chain]
A function $f: \sR^T\rightarrow \sR$ is a (first-order) zero-chain if 
\begin{align*}
    {\rm supp}(\vx) \subseteq \{1,\cdots,t-1 \} ~~ \Longrightarrow ~~ {\rm supp}(\nabla f(\vx)) \subseteq \{1,\cdots,t \}, \quad \forall \vx \in \sR^T.
\end{align*}
\end{dfn}

As an accompanying definition of the zero-chain, we denote $\gA^{\rm zr}$ as the class of zero-respecting methods \citep{carmon2020lower,carmon2021lower} that can only place mass on coordinates that
have already been revealed by previous gradient queries. Therefore, on a
zero-chain, such a method must discover the coordinates sequentially. 
\begin{dfn} \label{dfn:zero-respecting}
A zero-respecting method \(\mathsf A\in\mathcal A^{\rm zr}\), initialized at \(\vx^{(0)}=0\) without loss of generality, generates iterates \(\{\vx^{(t)}_{\mathsf A,f}\}_{t\ge0}\) that 
\[
{\rm supp}(\vx^{(t)}_{\mathsf A, f}) \subseteq \bigcup_{s < t} {\rm supp}(\nabla f(\vx^{(s)}_{\mathsf A,f})).
\]
\end{dfn}
Although zero-respecting methods $\gA^{\rm zr}$ form a more restricted algorithm class than all deterministic algorithms $\gA^{\rm det}$ in Definition \ref{dfn:rand-alg}, using a standard resisting-oracle argument of \citet{nemirovskij1983problem,nesterov2018lectures,carmon2020lower}, we can show that lower bounds against $\gA^{\rm zr}$ transfer to deterministic first-order methods $\gA^{\rm fo}$ for orthogonally invariant function classes. 
\begin{lem}[{\citet[Proposition 1]{carmon2020lower}}] \label{lem:carmon-alg-class-reduction}
Let $\gF$ be an orthogonal invariant function class, that is, if for every $f \in \gF$, $ \sR^d \rightarrow \sR $, and every orthogonal matrix $\mU \in \sR^{d' \times d}$ such that $\mU^\top \mU = \mI_d$, the function $f_\mU:\sR^{d'} \rightarrow \sR$ defined by $f_\mU(\vx) = f(\mU^\top \vx)$ belongs to $\gF$. For the deterministic first-order oracle $\sO$, we have
\begin{align*}
{\rm Compl}_{\epsilon}(\gA^{\rm fo}, \gF, \sO) \ge 
{\rm Compl}_{\epsilon}(\gA^{\rm zr}, \gF, \sO).
\end{align*}
\end{lem}
Therefore, it is without loss of generality to prove lower bounds on zero-respecting algorithms, for which we only need to construct zero-chains such that any point $\vx \in \sR^d$ has gradients larger than $\eps$ if the last coordinates remain zero. The optimal nonconvex zero-chains have been constructed by \citet{carmon2020lower,carmon2021lower}, and now let us recall one of their functions as follows.

\begin{dfn}[NC zero-chain \citep{carmon2021lower}] \label{dfn:nc-zo-chain}
Let $\widehat f_{r}^{\rm nc}: \sR^{T}\rightarrow \sR$ be
\begin{align*}
    \widehat f_{r}^{\rm nc}(\vx) = \frac{1}{2} (x_1 - 1)^2 + \frac{1}{2} \sum_{i=2}^T (x_{i}-  x_{i-1})^2
    + \sum_{i=1}^{T-1} \Upsilon_r(x_i),
\end{align*}
where the nonconvex regularizer $\Upsilon_r(x) = 120 \int_1^x \frac{t^2 (t-1)}{1+(t/r)^2} {\rm d}t$.
\end{dfn}

We summarize the key properties of this function in the following lemma.


\begin{lem}[{\citet[Lemma 2 and 3]{carmon2021lower}}] \label{lem:Car-Upsion}
 The functions $\Upsilon_r$ and $\widehat f_{r}^{\rm nc}$ 
 satisfy the following:
 \begin{enumerate}
 \item We have $\Upsilon_r'(0) = \Upsilon_r'(1) = 0$.
 \item For all $x \in \sR$, we have $\Upsilon_r(x) > \Upsilon_r(1) = 0$, and for all $r$, $\Upsilon_r(0) \le 10$.
    \item For every $r \ge 1$ and every $p \ge 1$, the $p$th-order derivatives of $\Upsilon_r$ are $r^{3-p}\ell_p$-Lipschitz continuous, where $\ell_p \le \exp(3p \log(p)/2 + c_0p)$ for a numerical constant $c_0>0$.
\item  Let $r \ge 1$, then for any $\vx \in \sR^T$ such that $x_{T-1} = x_{T}= 0$, we have 
\begin{align*}
    \Vert \nabla \widehat f^{\rm nc}_{r}(\vx) \Vert > 1 / 4.
\end{align*}
 \end{enumerate}
\end{lem}
Invoking this lemma and then properly scaling the hard instance in Definition \ref{dfn:nc-zo-chain}, \citet{carmon2021lower} proved lower bounds of $\Omega(\epsilon^{-2})$, $\Omega(\eps^{-12/7})$, and $\Omega(\epsilon^{-8/5})$ for first-order, second-order, and arbitrarily smooth nonconvex optimization, respectively. In the following, we show how to construct a hard instance for bilevel optimization based on their construction.

\subsection{Nested Chain for Bilevel Problems}
\label{subsec:nested-chain}

To apply the above framework to prove lower bounds for bilevel problems, we first give a natural generalization to the zero-chains and zero-respecting algorithms for a function pair $(f,g)$. In the following, we let $\vz = (\vx,\vy)$ be the aggregated variable.
\begin{dfn}[Zero-chain pair] We say
$f,g: \sR^T \rightarrow \sR$ is a pair of zero-chains if 
\begin{align*}
    {\rm supp}(\vz) \subseteq \{1,\cdots,t-1 \} ~~ \Longrightarrow ~~ {\rm supp}(\nabla f(\vz)) \cup {\rm supp}(\nabla g(\vz)) \subseteq \{1,\cdots,t \}, \quad \forall \vz  \in \sR^T.
\end{align*}
\end{dfn}

\begin{dfn}[Zero-respecting algorithms]
Operating on a function pair $(f,g)$,
we say \(\mathsf A\in\mathcal A^{\rm zr}\) is a zero-respecting method, if it initialized at \(\vz_{{\mathsf A,(f,g)}}^{(0)}=0\) without loss of generality, and generates iterates \(\{\vz^{(t)}_{\mathsf A,(f,g)}\}_{t\ge0}\) that 
\begin{align} \label{eq:dfn-zero-respecting}
{\rm supp}(\vz^{(t)}_{\mathsf A, (f,g)}) \subseteq \bigcup_{s < t} {\rm supp}\left(\nabla f(\vz^{(s)}_{\mathsf A,(f,g)}) \right) \cup {\rm supp} \left( \nabla g(\vz^{(s)}_{\mathsf A,(f,g)})\right).
\end{align}

\end{dfn}
NC--SC bilevel optimization typically requires a double-loop algorithm, where the inner loop solves the strongly convex subproblem to obtain $\vy^*(\vx)$ for a given $\vx$, and the outer loop performs (accelerated) gradient descent to find an approximate stationary point of the nonconvex hyper-objective $f(\vx,\vy^*(\vx))$. Since the optimal zero-chains for both the NC outer and SC inner loops are known \citep{carmon2021lower}, we hope that nesting these two chains to form a larger one can lead to the optimal complexity. 

\begin{figure}[t]
 \vspace{-1cm} 
\centering\includegraphics[width=0.8 \linewidth]{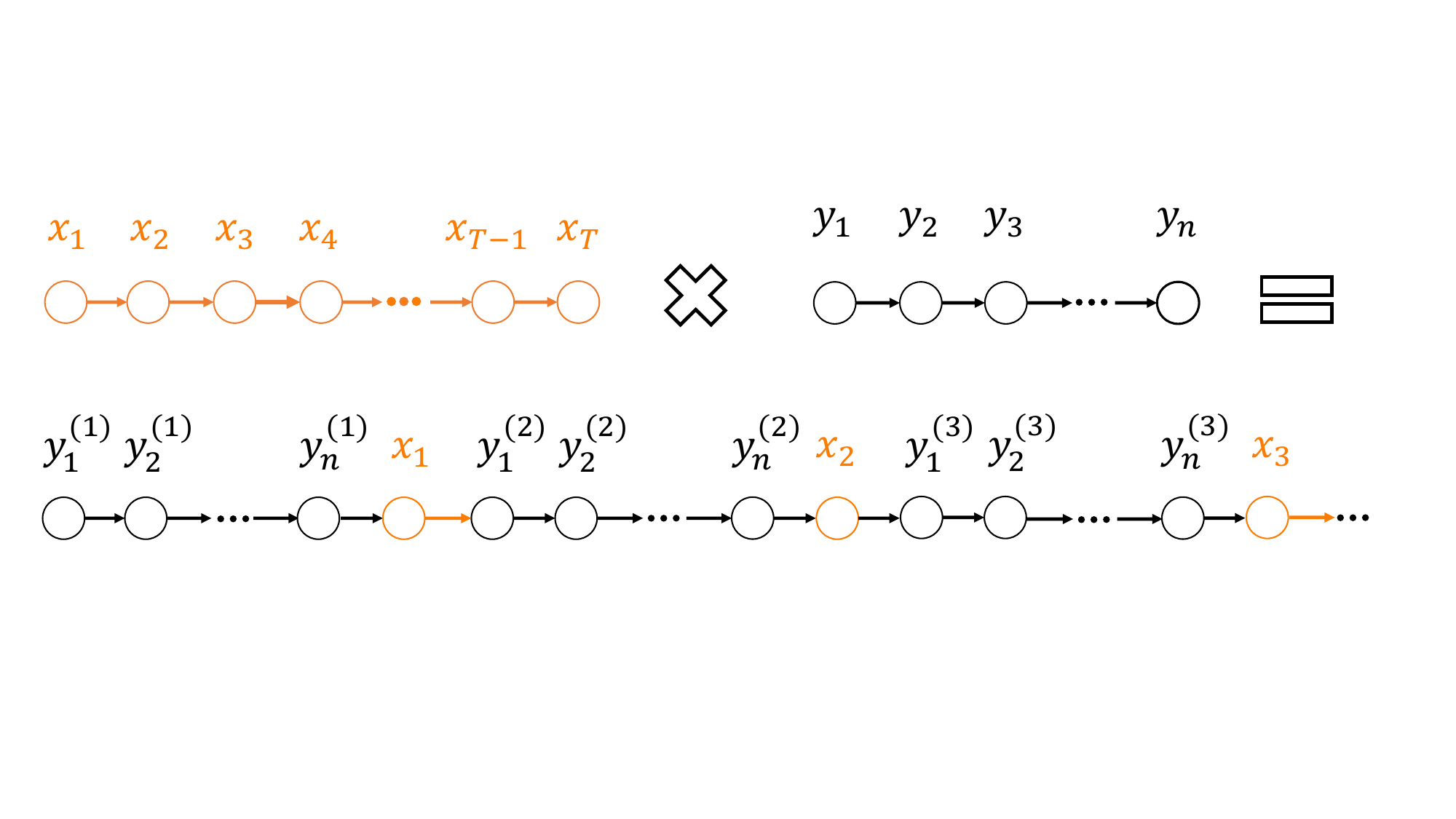}
 \vspace{-2cm}
\caption{Nesting two chains to form a longer zero-chain after reordering.
}
\end{figure} \label{fig:nested-chain}


As illustrated by Figure \ref{fig:nested-chain}, we construct a nested chain for NC--SC bilevel optimization as follows. 
Let $\vx=[x_1,...,x_{T}]\in\mathbb{R}^{T}$ be the upper-level variable, and the block vector ${\widetilde \vy} = [\vy^{(1)},\vy^{(2)},....,\vy^{(T)}]$ with each $\vy^{(i)}\in\mathbb{R}^n$ be the lower-level variable. Then we construct a sequential chain of $T$ independent lower-level sub-problems $ \widehat g(\vx,\widetilde \vy) = \sum_{i=1}^T g_n(x_i, \vy^{(i)})$ that are linked sequentially through the upper-level objective function $f({\widetilde \vy})$. Formally, we define the unscaled upper- and lower-level functions as follows:
\begin{align} \label{eq:our-class-det}
\begin{split}
    \widehat f(\widetilde \vy) =& \frac{1}{2} ( \alpha - 2y_1^{(1)}  )^2 + \frac{1}{2} \sum_{i=1}^{T-1} \left( y_n^{(i)} - 2y_1^{(i+1)} \right)^2 + \sum_{i=1}^{T} \Upsilon_1\left(y_n^{(i)}  \right), \\
    \widehat g(\vx,\widetilde \vy) = & \sum_{i=1}^T \left\{ g_n(x_i,\vy^{(i)}) \triangleq \frac{1}{2} (\vy^{(i)})^\top \mA_n \vy^{(i)} -  x_i y^{(i)}_n \right\},
\end{split}
\end{align}
where $\mA_n \in \sR^{n \times n}$ is a tridiagonal matrix that acts as a discrete Laplacian that enforces sequential
activation \citep{nesterov2018lectures,arjevani2019oracle,carmon2021lower}. In our setting, this matrix $\mA_n$ should be carefully designed to satisfy the following properties.
\begin{pro}[Exact reduction] \label{pro:exact-reduction}
Let $\vy^*(\vx) = \arg \min_{\vy \in \sR^n} g_n(x,\vy)$ be the optimal lower-level solution for $g_n(x,\,\cdot\,)$. We have
\[
\left( \langle \ve_n, \vy^*(x_1) \rangle - 2 \langle \ve_1, \vy^*(x_2) \rangle \right)^2 = \alpha^2 (x_1 - x_2)^2.
\]
Consequently, for $\widetilde \vy^*(\vx) = \arg \min_{\widetilde \vy \in (\sR^n)^T} \widehat g(\vx, \widetilde \vy)$, the hyper-objective $\widehat f(\widetilde \vy^*(\vx)) = \widehat f^{\rm nc}_{1} (\alpha \vx)$ exactly coincides with the nonconvex zero-chain in Definition \ref{dfn:nc-zo-chain} up to a scaling factor of~$\alpha$.
\end{pro}

\begin{pro}[Regular conditions for lower-level functions]
\label{pro:lower-level}
The lower-level function $\widehat g(\vx,\widetilde \vy)$ is $\gO(1)$-smooth and $\Omega(1/\kappa_\vy)$-strongly convex in $\widetilde \vy$.
\end{pro}

A consequence of the exact reduction property is the following lemma, which is analogous to Lemma~\ref{lem:carmon-alg-class-reduction} that reduces the lower bound for deterministic algorithms to the more restricted class of zero-respecting algorithms.

\begin{lem} \label{lem:our-alg-class-reduction}
Let $\gF^{\text{nc--sc}}(L_0,\cdots,L_{p+1},\mu_\vy, \Delta)$ be the NC--SC bilevel function class as per Definition \ref{dfn:NC-SC-func}. When Property \ref{pro:exact-reduction} is satisfied, for the deterministic first-order oracle $\sO$, we have
\begin{align*}
{\rm Compl}_{\epsilon}(\gA^{\rm det}, \gF, \sO) \ge 
{\rm Compl}_{\epsilon}(\gA^{\rm zr}, \gF, \sO).
\end{align*}
\end{lem}

\begin{proof}
The function class is orthogonal-invariant such that the conditions of Lemma \ref{lem:carmon-alg-class-reduction} can be applied. Therefore, the proof is completed by noting that the optimality measurement $\Vert \nabla_\vx f(\widetilde \vy^*(\vx)) \Vert$ is also orthogonal-invariant since $f(\widetilde \vy^*(\vx)) = \widehat f^{\rm nc}_{1} (\alpha \vx)$ under Property \ref{pro:exact-reduction}.
\end{proof}

Therefore, to establish lower bounds using the Construction (\ref{eq:our-class-det}), we only need to design the matrix $\mA_n$ to satisfy the above two properties, which is technical and left to the subsequent sections. As we will show in Section \ref{subsec:lb-det-nc}, Construction (\ref{eq:our-class-det}) can lead to the tight lower bound of $\Omega(\kappa_\vy^{5/2} \eps^{-2})$ for first-order smooth NC--Q bilevel problems. However, the same construction in the higher-order smooth settings results in suboptimal lower bounds of $\Omega( \kappa_\vy^{31/14} \eps^{-12/7} )$ and $\Omega(\kappa_\vy^{21/10} \eps^{-8/5})$ for second-order and arbitrarily smooth NC--Q problems, respectively. To show tight lower bounds for highly smooth problems, we require the following improved nonconvex hard instance based on a new nested chain.

{ \subsection{Nested Chain for Highly-Smooth Nonconvex Problems}
\label{subsec:hs-nested-chain}}

\begin{table}[tp]
   \caption{Comparisons of the first-order lower bounds by different works for $p$th-order smooth nonconvex optimization and the techniques they used.}
    \centering
    \begin{tabular}{c c cccc}
    \hline
  Bound  & Reference & $p=1$ & $p=2$ & $p \ge 3$  & Technique \\
    \hline \addlinespace
\multirow{3}{*}{Upper} & \citet{carmon2017convex} & $\gO(\epsilon^{-2})$ & $\tilde \gO(\epsilon^{-7/4})$ & $\tilde \gO(\epsilon^{-5/3})$ \\ \addlinespace
 & \citet{li2023restarted} & $\gO(\epsilon^{-2})$ & $\gO(\epsilon^{-7/4})$ & -  & \\ \addlinespace
 \hline \addlinespace
\multirow{4}{*}{Lower} & \citet{carmon2021lower} &
 {$\Omega(\epsilon^{-2})$} & $\Omega(\epsilon^{-12/7})$ & $\Omega(\epsilon^{-8/5})$ & Zero-Chain  \\ \addlinespace
& \citet{zhou2026sharp} & {$\Omega(\epsilon^{-2})$} & {$\Omega(\epsilon^{-7/4})$} 
& {$\Omega(\epsilon^{-5/3})$} & Block-Chain \\ \addlinespace
 & Ours & {$\Omega(\epsilon^{-2})$} & {$\Omega(\epsilon^{-7/4})$} 
& {$\Omega(\epsilon^{-5/3})$} & Nested Chain  \\
\hline 
    \end{tabular}
    \label{tab:lower-bound-NC}
\end{table}

To extend our lower bound for highly smooth problems, we need to first construct highly smooth nested chain for single-level nonconvex optimization, yielding the same results as \citet{zhou2026sharp}. The lemmas and structures within this construction will be used later.

Our key observation is that minimizing a nonconvex function $f(\vx)$ is naturally a bilevel optimization where the upper- and lower-level functions coincide
\[
\min_{\vx \in \sR^d} f(\vx) = \min_{\bar \vx \in \sR^{T+1}, \underline{\vx} \in (\sR^m)^T} f(\bar \vx, \underline{\vx}).
\]
This happens when the full dimension \(d = (m+1)T + 1\) and 
$\vx$ partitions into two groups $\bar \vx \in \sR^{T+1}$, $\underline{\vx} = (\underline{\vx}^{(1)}, \underline{\vx}^{(2)}, \cdots, \underline{\vx}^{(T)}) \in (\sR^{m})^T$, where \(\bar \vx\) is the outer variable and each \(\underline{\vx}^{(i)}\in\sR^m\) forms an inner chain connecting two consecutive outer variables \(\bar x_i\) and \(\bar x_{i+1}\). For highly smooth nonconvex optimization, we construct the (unscaled) nested chain as follows:
\begin{align} \label{eq:our-coupling-chain}
\begin{split}
&\widehat f^\text{hs-nc}_{\nu,r}(\vx) = \frac{\nu}{2} \left( 1 - \bar x_1 \right)^2  + \frac{1}{2} \sum_{i=1}^T \widehat Q_m(\bar x_{i}, \underline{\vx}^{(i)} , \bar x_{i+1} )  + \nu \sum_{i=1}^T \Upsilon_r \left( 
\bar x_i
\right), \\
&\widehat Q_m(x_-, \vx, x_+) =\frac{1}{m} (x_1 - \theta x_-)^2 + \sum_{j=1}^{m-1} (x_{j+1} - x_j)^2 + \frac{1}{m} (x_m - \theta x_+)^2,
\end{split}
\end{align}
which embeds a quadratic chain \citep{nesterov2018lectures} in the nonconvex zero-chain in Definition \ref{dfn:nc-zo-chain} to form a nested chain. The outer variables \(\bar x_i\) carry the same nonconvex stationarity obstruction as the classical hard instance of \citet{carmon2021lower}, whose chain length corresponds to the optimal complexity of conceptual high-order methods for nonconvex optimization \citep{birgin2017worst,carmon2020lower}. The inner variables \(\underline{\vx}^{(i)}\) force first-order methods to reveal the transition from \(\bar x_i\) to \(\bar x_{i+1}\) through an \(m\)-step quadratic chain, which corresponds to the optimal first-order complexity of solving the high-order regularized models \citep{agarwal2017finding,carmon2017convex,carmon2018accelerated}.

The following lemma shows that Construction (\ref{eq:our-coupling-chain}) also satisfies the exact reduction Property \ref{pro:exact-reduction}, because the embedded quadratic bridge \(\widehat Q_m(x_-,\,\cdot\,,x_+)\) reduces exactly to the scalar quadratic edge \((x_--x_+)^2\) up to scaling after partial minimization.

{
\begin{lem}[Exact Reduction to $\widehat f^{\rm nc}_r$] \label{lem:argmin-close-form}
Let $\nu =  \theta^2 /(3m-1)$. For fixed $x_-, x_+$, we have
\begin{align*}
    \min_{\vx \in \sR^m} \widehat Q_m(x_-,\vx, x_+) = \nu \left( {x_- - x_+} \right)^2.
\end{align*}
Hence, the hyper-objective $F(\bar \vx) = \min_{\underline{\vx} \in (\sR^m)^T} \widehat f^\text{hs-nc}_{\nu,r}(\vx) $  reduces to the nonconvex zero-chain $\widehat f^{\rm nc}_{r}$ in Definition \ref{dfn:nc-zo-chain} multiplied by $\nu$, that is, $F(\bar \vx) = \nu  \widehat f^{\rm nc}_{r}(\bar \vx)$.
\end{lem}}

\begin{proof}
See Appendix \ref{apx:proof-close-argmin} for the complete proof.
\end{proof}

Then we can directly apply Lemma \ref{lem:Car-Upsion} to show lower bounds for finding stationary points of the hyper-objective $F(\bar \vx)$. However, our goal is to show lower bounds on the gradient norm of the original function $\widehat f^\text{hs-nc}_{\nu,r}(\vx)$ instead of $F(\bar \vx)$. Therefore, we require an additional lemma to translate the two stationarities, which replaces the strong convexity constraint for bilevel optimization (Property \ref{pro:lower-level}) with the following one.

{
\begin{lem}[$f$-Stationarity implies $F$-stationarity] \label{lem:rela-stationarity}
If $\theta \sqrt{m} \le 1$, then the construction in Eq. (\ref{eq:our-coupling-chain}) satisfies \(\Vert \nabla \widehat f^\text{hs-nc}_{\nu,r}(\vx) \Vert \ge \Vert \nabla F(\bar \vx) \Vert /2\).
\end{lem}}

\begin{proof}
See Appendix \ref{apx:proof-rela-sta} for the complete proof.
\end{proof}


Finally, our nested chain construction in Eq. (\ref{eq:our-coupling-chain}) can give the following tight lower bounds for nonconvex optimization.

{
\begin{thm}[$p$th-order smooth NC lower bound] \label{thm:lb-highligh-smooth-NC}
For any integer \(p\ge1\) and sufficiently small precision \(\epsilon>0\), the 
first-order oracle complexity of finding an \(\epsilon\)-stationary point over the class of $p$th-order smooth nonconvex functions is at least
\begin{align} \label{eq:NC-lower-bound}
\begin{cases}
\Omega\!\left(\Delta L_1\epsilon^{-2}\right),
& p=1,\\[1mm]
\Omega\!\left(\Delta L_1^{1/2}L_2^{1/4}\epsilon^{-7/4}\right),
& p=2,\\[1mm]
\Omega\!\left(\Delta L_1^{1/2}L_3^{1/6}\epsilon^{-5/3}\right),
& p=3,\\[1mm]
\Omega_p\!\left(
\Delta L_1^{1/2}
\displaystyle\min_{1\le q\le p}
L_q^{1/(2q)}
\epsilon^{-5/3}
\right),
& p\ge4.
\end{cases}    
\end{align}
\end{thm}}

\begin{proof}
See Appendix \ref{apx:lb-highly-smooth-nc} for the complete proof.
\end{proof}

Our nested chain construction recovers the recently proved results by \citet{zhou2026sharp}, which are achieved by the algorithms in \citet{carmon2017convex} up to logarithmic factors. A slight difference is that \citet{zhou2026sharp} relies on the concept of block-chain, which we compare below.


\paragraph{Comparison to the Block-Chain.} Very recently,
\citet{zhou2026sharp} established first sharp lower bounds of 
$\Omega(\epsilon^{-7/4})$ and $\Omega(\epsilon^{-5/3})$ via the block-chain of the form:
\begin{align} \label{eq:zhou-block-chain}
\widehat f(\vx) =  \sum_{i=1}^T \frac{1}{2} \vx^{(i)} \mH_m \vx^{(i)} - \lambda \vu^{(i-1)} x^{(i)}_1 + \tau \left(\vu^{(i)}\right)^2+ \nu \Upsilon_r\left(\vu^{(i)} \right),
\end{align}
where $\vx = (\vx^{(1)}, \cdots, \vx^{(T)}) \in (\sR^m)^T$ is a block vector, $\mH_m \in \sR^{m \times m}$ is a tridiagonal matrix,  $\vu^{(i)} = \langle \va,\vx^{(i)} \rangle$ is a suffix readout for some $\va \in \sR^m$ such that ${\rm supp}(\va) = \{\lceil m/2 \rceil, \cdots,m \}$. It introduces a intricate mechanism to decouple the smoothness calibration of the nonconvex regularizer from the oracle-revelation process through a linear pullback induced by the vector $\va$.
In contrast, we design inner chains such that the hyper-objective is identical to the existing hard instance in Definition \ref{dfn:nc-zo-chain}, allowing us to directly reuse many established analyses in \citet{carmon2017convex} and provide an alternative proof of \citet{zhou2026sharp}. Finally, we compare our results with the lower bounds by \citet{carmon2017convex,zhou2026sharp} and the different techniques we use in Table~\ref{tab:lower-bound-NC}.

Compared to the block-chain, our construction may be more convenient to used in bilevel problems, because the exact reduction Property \ref{pro:exact-reduction} also immediately applies to the quadratic links $\sum_{j=1}^{m-1} ( x_{j+1} - x_j )^2$ in Construction (\ref{eq:our-coupling-chain}).

\section{Lower Bounds for Deterministic Bilevel Optimization} \label{subsec:lb-det}

In this section, we prove lower bounds for deterministic bilevel optimization by explicitly constructing the tridiagonal matrix $\mA_n$ in Eq. (\ref{eq:our-coupling-chain}).

\subsection{Lower Bounds for Nonconvex Bilevel Problems} \label{subsec:lb-det-nc}

A first natural idea is to set $\mA_n =  \mB_n + \mI_n / n^2$, where $\mB_n \in \sR^{n \times n}$ is the tridiagonal discrete Laplacian matrix used in classical lower-bound constructions for convex optimization \citep{nesterov2018lectures,woodworth2016tight}  or nonconvex--strongly-concave  (NC--SC) minimax optimization \citep{li2021complexity,zhang2021complexity,wang2024efficient}:
{\begin{align} \label{eq:finite-A}
   \mB_n
    =
    {
    \begin{bmatrix}
        1 & -1 \\
        -1 & 2 & -1 \\
        & -1 & \ddots & \ddots  \\
        & & \ddots & 2 & -1 \\
        & &  & -1 & 1\\
    \end{bmatrix}}.
\end{align}
By \citet[Lemma 7]{li2021complexity}, this matrix satisfies the required Property \ref{pro:exact-reduction} with $\alpha = \sqrt{\kappa_\vy}$. However, the scaling of $\alpha = \sqrt{\kappa_\vy}$ is not tight and can only result in suboptimal lower bounds of $\Omega(\kappa_\vy^{3/2} \eps^{-2})$ for NC--SC bilevel problems. This is because the influence of the first coordinate is roughly uniform across all indices \citep[Lemma 7]{li2021complexity}, which leads to suboptimal endpoint values of the inverse matrix $\mA_n^{-1}$.

To address the above issue, we propose an improved matrix to push the mass of its inverse to the endpoints, which can improve both $(\mA_n^{-1})_{1,1}$ and $(\mA_n^{-1})_{1,n}$ from $\Omega(\sqrt{\kappa_\vy})$ to $ \Omega(\kappa_\vy)$. Specifically, we embed the prior Laplacian matrix into the middle coordinates and connect it to the two endpoints through carefully
chosen negative weights $-q$ and positive endpoint diagonal $\omega$. We let
{\begin{align} \label{eq:finite-A-endpoint}
   \mA_n
    =
    {
    \begin{bmatrix}
    {\omega} &{-q}&0&\cdots&0\\
    {-q} &2&-1&\cdots&0\\
    0&-1&2&\ddots&\vdots\\
    \vdots&&\ddots&2&{-q}\\
    0&\cdots&0&{-q}&{\omega}
    \end{bmatrix}
   }.
\end{align}
In addition, the parameters $\omega$ and $q$ are carefully selected so that, after eliminating the middle block, the Schur complement on the endpoint coordinates becomes exactly 
\[
\mS  =
\frac{1}{(n-1)^2}
  {\begin{bmatrix}
    2&-1\\
    -1&2
    \end{bmatrix}}.
\]
Consequently, the endpoint block of $\mA_n^{-1}$ is $\mS^{-1}$, which immediately leads to $\Omega(n^2)$ endpoint values of $\mA_n^{-1}$. Combined with the property $\mA_n \succeq \Omega(1/n^2) \mI_n$ shown below, we can successfully conclude that the new matrix $\mA_n$ satisfies Property \ref{pro:exact-reduction} with $\alpha = \Omega(\kappa_\vy)$.
Formally, the properties of this new inverse-endpoint-amplified matrix are summarized as follows.
{ 
\begin{lem} \label{lem:matrix-A}
For $n \ge 5$, let $q = 1 / \sqrt{n-1}$ and $\omega = n / (n-1)^2$. Then the matrix $\mA_n$ in Eq.~(\ref{eq:finite-A-endpoint}) satisfies
    \begin{enumerate}[label=(\alph*)]
        \item $(\mA_n^{-1})_{1,1} = (\mA_n^{-1})_{n,n}=  {2 (n-1)^2}/{3}$  and  $(\mA_n^{-1})_{1 ,n} = (\mA_n^{-1})_{n,1} = {(n-1)^2}/{3}$. 
        \item $\mI_{n} / (2(n-1)^2) \preceq \mA_n \preceq 5 \mI_{n}$.
        \end{enumerate}
\end{lem}}
\begin{proof} 
See Appendix \ref{apx:proof-property-Bk} for the complete proof.
\end{proof}
Then, plugging this matrix into the Construction (\ref{eq:our-class-det}) 
indicates that the (unscaled) hyper-objective is $\widehat f(\widetilde \vy^*(\vx)) = \widehat f^{\rm nc}_r(\alpha \vx)$ under the setting  \(\alpha = 2(n-1)^2/3\), which leads to the following lower bound under the scaling $n = \Theta(\sqrt{\kappa_\vy})$ and $\alpha = \Theta(\kappa_\vy)$.

\begin{thm}[NC--Q lower bound] \label{thm:NC-SC}
There exists a numerical constant $a_0 \in (0,1)$ such that, 
for any $L_1,\Delta >0$ and $\mu_\vy \in (0, a_0 L_1]$, there exists a smooth nonconvex--quadratic bilevel problem $(f, g) \in \gF^{\text{nc--q}}(L_1,\mu_\vy, \Delta)$ such that the complexity of any deterministic first-order algorithm $\texttt{A} \in \gA^{\rm fo}$ or HVP-based algorithm $\texttt{A} \in \gA^{\rm hvp}$  to find an $\epsilon$-stationary point of the  hyper-objective  $F(\vx)$ is at least
\[
\Omega \left( \kappa_\vy^{5/2} 
             L_1 \Delta \epsilon^{-2} \right),
\]  
where $\kappa_\vy = L_1/\mu_\vy$ is the condition number of the lower-level problem.
\end{thm}

\begin{proof}
For the unscaled upper- and lower-level functions $\widehat f$ and $\widehat g$ in Construction (\ref{eq:our-class-det}), we scale them as $f(\widetilde \vy) = \delta \beta^2 \widehat f(\widetilde \vy / \beta)$ and $g(\vx,\widetilde \vy) = \delta \beta^2 \widehat g(\vx / \beta, \widetilde \vy / \beta)$, which are $L_1$-smooth if $\delta = \Theta(L_1)$ according to Lemma \ref{lem:Car-Upsion} and \ref{lem:matrix-A}. Since the hyper-objective 
$f(\widetilde \vy^*(\vx)) = \delta \beta^2 \widehat f^{\rm nc}_1(\alpha \vx / \beta)$ for \(\alpha = 2(n-1)^2/3\), where $\widehat{f}^{\text{nc}}_r(\vx)$ is the $T$-dimensional nonconvex zero-chain function in Definition \ref{dfn:nc-zo-chain}. Moreover, Lemma \ref{lem:matrix-A}(b) indicates that the lower-level function $g(\vx,\widetilde \vy)$ is $\mu_\vy$-strongly convex under the setting $n = \Theta(\sqrt{\kappa_\vy})$, leading to the scaling $\alpha = \Theta(\kappa_\vy)$. Moreover, for the hyper-objective function, Lemma \ref{lem:Car-Upsion} suggests that 
\[
\vx_{T-1} = \vx_T = 0 \quad \Longrightarrow \quad
\Vert \nabla_\vx f(\widetilde \vy^*(\vx)) \Vert \ge \delta \alpha \beta/4.
\]
Recalling Lemma \ref{lem:our-alg-class-reduction}, the lower bound for deterministic algorithms can be reduced to zero-respecting algorithms, for which the previous $n(T-1)$ queries satisfy $\Vert \nabla_\vx f(\widetilde \vy^*(\vx)) \Vert >\eps$ under the parameter setting
\[
\beta = \frac{4 \eps}{\alpha \delta} = \Theta \left(\frac{\eps}{\kappa_\vy L_1} \right).
\]
It remains to select parameters $T$ as large as possible, such that the function pair $(f,g) \in \mathcal{F}^\text{nc--q}(L_1,\mu_\vy,\Delta)$, and then conclude the lower bound of $n (T-1) = \Omega( \sqrt{\kappa_\vy} T)$. From Lemma~\ref{lem:Car-Upsion}, we know $\Upsilon_1(0) \le 10$ and $\Upsilon_1(1) = 0$. For the starting point $\vx^{{(0)}} = \vzero$, we have 
$$f(\widetilde \vy^*(\vzero)) = \delta \beta^2
\widehat{f}^{\text{nc}}_1(\vzero) = \delta \beta^2 \left(\frac{1}{2} + \sum_{i=1}^{T} \Upsilon_1(0) \right) 
\le \delta \beta^2 \left(\frac{1}{2} + 10T\right).
$$
Since $\widehat{f}^{\text{nc}}_1(\vx) \ge 0$ everywhere, the suboptimal gap is bounded by $\Delta$ under the setting 
\[
T = \left \lfloor \frac{\Delta}{10 \delta \beta^2} - \frac{1}{2} \right \rfloor = \Theta \left( \frac{\Delta L_1 \kappa_\vy^2}{\eps^{2}} \right).
\] 
Putting all the pieces together leads to the final lower bound of \(\Omega (\Delta L_1 \kappa_\vy^{5/2} \eps^{-2}) \).
\end{proof}


Compared with the known lower bound of $\Omega(\sqrt{\kappa_\vy} \eps^{-2})$ for NC--Q minimax problems 
\citep{li2021complexity,zhang2021complexity}, our results show that NC--Q bilevel optimization is provably more challenging than NC--Q minimax optimization, with the increased $\kappa_\vy$ dependency on $\kappa_\vy^{5/2}$ reflecting a fundamental difficulty stemming from the bilevel problem structure. Compared with the upper bounds of $\gO(\kappa_\vy^{5/2} \eps^{-2}) $ and $\gO(\bar \kappa_\vy^{7/2} \eps^{-2}) $ for NC--Q and general NC--SC
problems achieved by applying a lower-level Nesterov acceleration in the method of \citet{chen2023near} as shown in Appendix \ref{apx:ub-refined}, we can conclude that our lower bound is tight up to logarithmic factors for NC--Q problems, but still has a gap of $\kappa_\vy$ between our lower bounds for general NC--SC problems.

\subsection{Lower Bounds for Highly-Smooth Nonconvex Bilevel Problems}

Then, we consider the generalization of the above lower bound to all orders of smoothness. We nest the same lower-level chain in Eq. (\ref{eq:our-class-det}) with the optimal highly-smooth zero-chain $\widehat f^\text{hs-nc}_{\nu,r}$ introduced in Eq. (\ref{eq:our-coupling-chain}). Formally,
we introduce one lower-level block $\bar \vy^{(i)}, \underline{\vy}^{(i,j)} \in \sR^n $ for 
every scalar upper-level coordinate $\bar \vx_i, \underline{\vx}^i_j \in \sR$, and define 
\begin{align}  \label{eq:our-bilevel-fg-highly-smooth}
\begin{split}
&\widehat f(\vy) = \frac{\nu}{2} \left( \left( \alpha - 2\bar y_1^{(1)} \right)^2 \right) + \nu \sum_{i=1}^T \Upsilon_r \left( 
\bar y_n^{(i)}
\right) \\
&+ \frac{1}{2} \sum_{i=1}^T \left(\frac{1}{m} \left(\theta \bar y_n^{(i)} - 2\underline{y}^{(i,1)}_1 \right)^2 + \sum_{j=1}^{m-1} \left( \underline{y}^{(i,j+1)}_n - 2\underline{y}^{(i,j)}_1 \right)^2  + \frac{1}{m} \left( \underline{y}^{(i,m)}_n - 2\theta \bar y^{(i+1)}_1 \right)^2 \right), \\
& \widehat g(\vx,\vy) = \sum_{i=1}^T \left( g_n (\bar x_i, \bar \vy^{(i)}) + \sum_{j=1}^m g_n(\underline{x}_j^{(i)}, \underline{\vy}^{(i,j)}) \right)
\end{split}
\end{align}
Similarly, using Lemma \ref{lem:matrix-A} (a), we can also deduce that the hyper-objective $f(\widetilde \vy^*(\vx))$ reduces to the highly-smooth nonconvex nested-chain $\widehat f^\text{hs-nc}_{\nu,r}(\alpha \vx)$ in Eq. (\ref{eq:our-coupling-chain}) after applying a rescaling by $\alpha$, which, combined with Theorem \ref{thm:lb-highligh-smooth-NC}, leads to the following lower bounds.

\begin{thm}[$p$th-order NC--Q lower bound] \label{thm:highly-smooth-NC-Q}
There exists a numerical constant $a_0 \in (0,1)$ such that, 
for any $p \in \sN_+$, $L_1,\cdots,L_p>0$, $\Delta >0$ and $\mu_\vy \in (0, a_0 L_1]$, there exists a $p$th-order smooth nonconvex--quadratic bilevel problem $(f, g) \in \gF^{\text{nc--q}}(L_1,\cdots,L_p, \mu_\vy, \Delta)$ such that 
the complexity of any deterministic first-order algorithm $\texttt{A} \in \gA^{\rm fo}$ or HVP-based algorithm $\texttt{A} \in \gA^{\rm hvp}$
to find an $\epsilon$-stationary point of the  hyper-objective  $F(\vx)$ is at least
\[
\begin{cases}
\Omega\!\left(
    \kappa_\vy^{5/2} \,\Delta L_1\epsilon^{-2}
\right),
& p=1,\\[1mm]
\Omega\!\left(
    \kappa_\vy^{9/4}\,\Delta L_1^{1/2}L_2^{1/4}\epsilon^{-7/4}
\right),
& p=2,\\[1mm]
\Omega\!\left(
    \kappa_\vy^{13/6} \,\Delta L_1^{1/2}L_3^{1/6}\epsilon^{-5/3}
\right),
& p=3,\\[1mm]
\Omega_p\!\left(
   \kappa_\vy^{13/6} \,\Delta L_1^{1/2}
    \displaystyle\min_{1\le q\le p} L_q^{1/(2q)}
    \epsilon^{-5/3}
\right),
& p\ge4,
\end{cases}
\]
where $\eps \downarrow 0$ and $\kappa_\vy = L_1/\mu_\vy$ is the condition number of the lower-level problem.
\end{thm}

\begin{proof} 
Lemma \ref{lem:matrix-A}(b) indicates that the lower-level function $g(\vx,\widetilde \vy)$ is $\Omega(1/\kappa_\vy)$-strongly convex under the setting $n = \Theta(\sqrt{\kappa_\vy})$, leading to the scaling $\alpha = 2 (n-1)^2 / 3 = \Theta(\kappa_\vy)$. Then, for $\widehat f(\widetilde \vy^*(\vx)) = \widehat f^\text{hs-nc}_{\nu,r}(\alpha \vx )$, in order to guarantee $ \Vert \nabla_\vx \widehat f(\widetilde \vy^*(\vx)) \Vert \le \epsilon $, it is equivalent to ensure that $\Vert \nabla \widehat f^\text{hs-nc}_{\nu,r}(\alpha \vx)  \Vert \le \eps / \alpha$. Therefore, in the bilevel construction, we use the same parameter choices (up to numerical constants) as
in Theorem~\ref{thm:lb-highligh-smooth-NC}, but with the target precision
\(\epsilon/ \alpha \) instead of \(\epsilon\). To be more precise, we verify that this parameter setting ensures that our construction lies in the function class.

The initial gap is preserved under the rescaling from  \(\vx\) to \(\alpha \vx\). Indeed,
\[
    \widehat f(\vzero)-\inf_{\vx \in (\sR^m)^T} \widehat f(\widetilde \vy^*(\vx))
    =
    \widehat f^\text{hs-nc}_{\nu,r}(\vzero)-\inf_{\vx \in (\sR^m)^T} \widehat f^\text{hs-nc}_{\nu,r}(\vx).
\]
The smoothness constraints are also satisfied. For the unscaled upper- and lower-level functions $\widehat f$ and $\widehat g$ in Eq.~(\ref{eq:our-bilevel-fg-highly-smooth}), the additional quadratic edges contribute at most \(\gO(1)\) to the first-order
smoothness compared to the highly smooth nonconvex nested chain in Eq.~(\ref{eq:our-coupling-chain}). The remaining higher-order derivatives come only from the
regularizer \(\Upsilon_r\), which is exactly as in Theorem
\ref{thm:lb-highligh-smooth-NC}. 

Therefore, up to numerical constants, the largest possible value of $m (T-1)$ can be set as Eq. (\ref{eq:NC-lower-bound}) with $\epsilon$ replaced by $\eps / \alpha$.
Putting all the pieces together leads to the $n m(T-1)$ lower bound of $\Omega( n (\alpha/\eps)^{2} ) $, $\Omega( n (\alpha/\eps)^{7/4} ) $, and $\Omega( n (\alpha/\eps)^{5/3})$ for first-order, second-order, and arbitrarily smooth problems. Recalling that $n = \Theta(\sqrt{\kappa_\vy})$ and $\alpha = \Theta(\kappa_\vy)$, we obtain the lower bound of $\Omega(\kappa_\vy^{5/2} \epsilon^{-2})$, $\Omega(\kappa_\vy^{9/4} \eps^{-7/4} )$, and $\Omega( \kappa_\vy^{13/6} \eps^{-5/3})$, respectively.
\end{proof}

In Appendix \ref{apx:ub-refined}, we show that the method by \citet{chen2023near} combined with a lower-level Nesterov acceleration achieves the upper bound of $\tilde \gO(\kappa_\vy^{9/4} \eps^{-7/4})$ and $\tilde \gO(\bar \kappa_\vy^{13/4} \eps^{-7/4})$ for second-order smooth NC--Q and NC--SC problems, respectively. By comparison, we know that our lower bound for $p=2$ is also nearly tight for NC--Q problems. We conjecture that our lower bound of $\Omega(\kappa_\vy^{13/6} \epsilon^{-5/3})$ is also nearly tight for the case $p\ge3$, and it is possible to achieve the bound up to logarithmic factors by extending \citet{carmon2017convex}.

\subsection{Lower Bounds for Convex Bilevel Problems} \label{subsec:lb-det-sc}

\citet{ji2023lower} proved lower bounds of $\tilde \Omega(\min \{ \kappa_\vy, \epsilon^{-3/2} \} / \sqrt{\epsilon}  )$ and
$\tilde \Omega (\kappa_\vy \sqrt{\kappa_\vx})$   
for the convex--quadratic (C--Q) and
strongly-convex--quadratic (SC--Q) problems, which demonstrates the first larger $\kappa_\vy$ dependency than the lower bounds of  $\Omega( \sqrt{\kappa_\vy/ \epsilon} )$ and $\tilde \Omega(\sqrt{\kappa_\vy \kappa_\vx})$ for the corresponding minimax problems \citep{ouyang2021lower,zhang2022lower}.
However, we show that the lower bounds in \citet{ji2023lower} are still not tight and can be improved to $\Omega(\kappa_\vy^{3/2} /\sqrt{\epsilon} )  $ and $\Omega( \kappa_\vy^{3/2} \sqrt{\kappa_\vx} )$ using our construction. To begin with, we formally define smooth (S)C--SC follows. 




\begin{dfn}[(S)C--SC problems] \label{dfn:SC-SC-func}
Given
$L_0,L_1, L_2>0$, $\mu_\vx \in [0,L_1]$, $\mu_\vy \in (0,L_1]$, and $D>0$, we use $\gF^{\text{sc--sc}}(L_0,L_1,L_2,\mu_\vx,\mu_\vy, D)$ to denote the subset of all bilevel problems $(f,g)$ 
that satisfy items 1 and 2 in Definition \ref{dfn:NC-SC-func}, \textit{i.e.}, 
$\gF^{\text{sc--sc}}(L_0,L_1,L_2,\mu_\vx,\mu_\vy, D)$ $\subset $ $\gF^{\text{nc--sc}} (L_0,L_1,L_2, \mu_\vy, +\infty )$, and the hyper-objective
$F(\vx)$ also satisfies
\begin{enumerate}
    \item $F(\vx)$ is $\mu_\vx$-strongly convex in $\vx$. 
    \item $\min\{ \Vert \vx \Vert : \vx \in  \arg \min_{\vx \in \sR^{d_x}} F(\vx)  \} \le D$.
\end{enumerate}
If $\mu_\vx = 0$, we call it a C--SC problem: $\gF^{\text{c--sc}}(L_0,L_1,L_2,\mu_\vy, D) = \gF^{\text{sc--sc}}(L_0,L_1,L_2,0,\mu_\vy, D)$.
\end{dfn}

We note that the (S)C--SC problems are not defined by the convexity of the upper-level function $f$, but instead by the convexity of the hyper-objective function $F$, making the definition potentially difficult to check in practice, as discussed in~\citet{hong2023two}.
However, in theory, it provides a simple yet important setting to understand the complexity of bilevel optimization \citep{ji2023lower}. In addition, (S)C--Q problems are the subset of (S)C--SC problems when the lower-level function is quadratic.


\begin{dfn}[(S)C--Quadratic problems]  \label{dfn:SC-SC-func-quadra}
Given $L_1>0$, $\mu_\vx \in [0,L_1)$,$\mu_\vy \in (0,L_1]$, and $D>0$, we denote
$\gF^{\text{sc--q}}(L_1,\mu_\vx,\mu_\vy, D) \subseteq \gF^{\text{sc--sc}}(+\infty, L_1, 0 ,\mu_\vx,\mu_\vy, D)$ as the subset of SC--SC problems such that both the upper- and lower-level functions are quadratic. If $\mu_\vx=0$, we call them C--Q problems: $\gF^\text{c--q}(L_1,\mu_\vy,D) = \gF^\text{sc--q}(L_1,0,\mu_\vy,D)$.
\end{dfn}

In the following, we show improved lower bounds for C--SC and SC--SC bilevel problems by reusing the construction in Eq. (\ref{eq:our-class-det}) but simply removing the nonconvex regularizers $\Upsilon_r$. We first note that if the regularizers $\Upsilon_r$ in the nonconvex zero-chain $\widehat f_{r}^{\rm nc}$ in Definition~\ref{dfn:nc-zo-chain} are removed, the function reduces to a finite version of Nesterov's convex zero-chain in Definition \ref{dfn:finite-Nes-func} and has the following property stated in Lemma \ref{lem:Car-large-grad-convex}.

\begin{dfn}[Convex zero-chain \citep{carmon2021lower}] \label{dfn:finite-Nes-func}
Let $\widehat f^{\rm c}: \sR^T \rightarrow \sR$ be the convex function  with the nonconvex regularizer $\Upsilon_r$ in Definition \ref{dfn:nc-zo-chain} removed, \textit{i.e.},
\begin{align*}
    \widehat f^{\rm c}(\vx)= \frac{1}{2} (x_1 - 1)^2 + \frac{1}{2} \sum_{i=2}^T (x_{i}-  x_{i-1})^2.
\end{align*}
\end{dfn}

\begin{lem}[{\citet[Lemma 1]{carmon2021lower}}] \label{lem:Car-large-grad-convex}
For any $\vx \in \sR^T$ such that $x_{T}= 0$, 
\begin{align*}
    \Vert \nabla \widehat f^{\rm c}(\vx) \Vert > 1/ T^{3/2}.
\end{align*}
\end{lem}

Hence, if we remove the nonconvex regularizer in Eq. (\ref{eq:our-class-det}) for NC--Q problems and keep all the other components unchanged, the hyper-objective $f(\widetilde \vy^*(\vx)) = \widehat f^{\rm c}(\alpha \vx)$ under the same setting \(\alpha = 2(n-1)^2/3\). This leads to the following lower bound for C--Q problems.


\begin{thm}[C--Q bilevel lower bound] \label{thm:C-SC}
There exists a numerical constant $a_0 \in (0,1)$ such that, 
for any $ L_1>0$, $D >0$ and $\mu_\vy \in (0, a_0  L_1]$, there exists a quadratic   bilevel problem $(f, g) \in \gF^{\text{c--q}}(L_1, \mu_\vy, D)$ such that, the complexity of
any deterministic first-order algorithm $\texttt{A} \in \gA^{\rm fo}$ or HVP-based algorithm $\texttt{A} \in \gA^{\rm hvp}$ to find an $\epsilon$-stationary point of the  hyper-objective  $F(\vx)$ is at least
    \begin{align*}
        \Omega\left( \kappa_\vy^{3/2} \sqrt{  L_1  D/ \epsilon} \right),
    \end{align*}
where $\kappa_\vy = L_1/\mu_\vy$ is the condition number of the lower-level problem.
\end{thm}

\begin{proof}
For the unscaled upper- and lower-level functions $\widehat f$ and $\widehat g$ in Eq. (\ref{eq:our-class-det}), we scale them as $f(\widetilde \vy) = \delta \beta^2 \widehat f(\widetilde \vy / \beta)$ and $g(\vx,\widetilde \vy) = \delta \beta^2 \widehat g(\vx / \beta, \widetilde \vy / \beta)$, which are $L_1$-smooth if $\delta = \Theta(L_1)$ according to Lemma \ref{lem:Car-Upsion} and \ref{lem:matrix-A}. Let
$\widehat{f}^{\text{c}}(\vx)$ be the $T$-dimensional convex zero-chain function in Definition \ref{dfn:finite-Nes-func}.
According to Lemma \ref{lem:matrix-A}(a), our constructed hard instance ensures that the hyper-objective 
$f(\widetilde \vy^*(\vx)) = \delta \beta^2 \widehat f^{\rm c}(\alpha \vx / \beta)$ for \(\alpha = 2(n-1)^2/3\). Moreover, Lemma~\ref{lem:matrix-A}(b) indicates that the lower-level function $g(\vx,\widetilde \vy)$ is $\mu_\vy$-strongly convex by setting $n = \Theta(\sqrt{\kappa_\vy})$ and thus $\alpha = \Theta(\kappa_\vy)$. Then, for $f(\widetilde \vy^*(\vx)) = \delta \beta^2 \widehat f^{\rm c}(\alpha \vx / \beta)$, Lemma \ref{lem:Car-large-grad-convex} suggests that 
\[
\vx_T = 0 \quad \Longrightarrow \quad
\Vert \nabla_\vx f(\widetilde \vy^*(\vx)) \Vert \ge \delta \alpha \beta/ T^{3/2}.
\]
Recalling Lemma \ref{lem:our-alg-class-reduction} that the lower bound for deterministic algorithms can be reduced to zero-respecting algorithms, for which the previous $n(T-1)$ queries satisfy $\Vert \nabla_\vx f(\widetilde \vy^*(\vx)) \Vert >\eps$ under the parameter setting
\[
\beta = \frac{T^{3/2} \eps}{\alpha \delta} = \Theta \left(\frac{T^{3/2} \eps}{\kappa_\vy L_1} \right).
\]
It remains to select parameters $T$ as large as possible, such that the function pair $(f,g) \in \mathcal{F}^\text{c--q}(L_1,\mu_\vy,D)$, and then conclude the lower bound of $n T = \Omega( \sqrt{\kappa_\vy} \, T)$. By the form of the function $\widehat f^{\rm c}$ in Definition \ref{dfn:finite-Nes-func}, it is easy to see that the minimizer of $f(\widetilde \vy^*(\vx)) = \delta \beta^2 \widehat f^{\rm c}(\alpha \vx/ \beta)$ is $\vx^* = [\beta/\alpha,\beta/\alpha,\cdots,\beta/\alpha ] \in \sR^T$. Thus, we have
\[
\Vert \vx^* \Vert = \frac{\beta \sqrt{T}}{\alpha} = \frac{T^2 \eps}{\alpha^2 \delta} = \Theta\left( 
\frac{T^2 \eps}{\kappa_\vy^2 L_1}
\right).
\]
Consequently, we can ensure the condition $\Vert \vx^* \Vert \le D$ under the parameter setting
\[
T = \left \lfloor \alpha \sqrt{ \delta D/\eps}  \right \rfloor = \Theta\left(\kappa_\vy \sqrt{L_1 D/ \eps} \right).
\]
Putting all the pieces together leads to the lower bound of $nT = \Omega(\kappa_\vy^{3/2} / \sqrt{\eps})$ as claimed. 
\end{proof}


Moreover, due to a black-box reduction from strongly convex problems to convex problems \citep{allen2016optimal} or \citep[Remark 5]{hinder2020near}, we can also convert the above lower bound for C--Q problems to SC--Q problems as follows.

\begin{thm}[SC-Q lower bound] \label{thm:SC-SC}
There exists a numerical constant $a_0 \in (0,1)$ such that, 
for any $ L_1>0$, $D >0$, $\mu_\vx \in (0, a_0 L_1]$, $\mu_\vy \in (0, a_0 L_1]$, there exists a quadratic bilevel problem $( f, g) \in \gF^{\text{sc--q}}(L_1, \mu_\vx, \mu_\vy, D)$ such that, the complexity of any deterministic first-order algorithm $\texttt{A} \in \gA^{\rm fo}$ or HVP-based algorithm $\texttt{A} \in \gA^{\rm hvp}$ 
to find an $\epsilon$-stationary point of the hyper-objective  $F(\vx)$ is at least
    \begin{align*}
        \Omega \left( \kappa_\vy^{3/2}
 \sqrt{\kappa_\vx} \right),
    \end{align*}
where $\kappa_\vx = L_1/\mu_\vx$ and $\kappa_\vy = L_1/\mu_\vy$.
\end{thm}

\begin{proof}
Suppose there was a deterministic first-order/HVP method \texttt{A} for SC--Q bilevel problems $\mathcal{F}^\text{sc--q}(L_1, \mu_\vx, \mu_\vy, D)$ that only requires $o( \kappa_\vy^{3/2} \sqrt{\kappa_\vx} )$ first-order calls to find an $\eps$-hyper-stationary point such that $\Vert \nabla_\vx f(\widetilde \vy^*(\vx)) \Vert \le \eps$. Let $(f,g) \in \mathcal{F}^\text{c--q}(L_1, \mu_\vy, D)$ be the hard instance for C--Q bilevel problems such that any deterministic first-order/HVP method requires at least $\Omega( \kappa_\vy^{3/2} \sqrt{L_1 D / \eps} )$ first-order oracles to find an $\eps$-hyper-stationary point. We apply \texttt{A} on the instance $(\tilde f, g)$, where the upper-level function is regularized by
\[
\tilde f(\vx,\vy) = f(\vy) + \frac{\eps}{4 D} \Vert \vx - \vx^* \Vert^2.
\]
Then, the regularized hard instance $(\tilde f, g) \in \mathcal{F}^\text{sc--q}(L_1+ \eps/(2D), \eps/(2D), \mu_\vy, D)$, which means $\kappa_\vx = 1 + 2 L_1 D / \eps$. Therefore, applying algorithm \texttt{A} on $(\tilde f, g)$ can find an $\eps$-hyper-stationary point such that 
$\Vert \nabla_\vx \tilde f(\widetilde \vy^*(\vx)) \Vert \le \eps/2$ within the first-order oracle complexity of $ o( \kappa_\vy^{3/2} \sqrt{L_1 D/ \eps} ) $. However, it contradicts the lower bound for C--Q problems since it indicates that $\vx$ is an $\eps$-hyper-stationary point such that $\Vert \nabla_\vx f(\widetilde \vy^*(\vx)) \Vert \le \eps$:
\[
\Vert \nabla_\vx f(\widetilde \vy^*(\vx)) \Vert \le \Vert \nabla_\vx \tilde f(\widetilde \vy^*(\vx)) \Vert + \frac{\eps}{2D} \Vert \vx - \vx^* \Vert \le \frac{\eps}{2} + \frac{\eps}{2} \le \eps.
\]
Therefore, the lower bound for SC--Q problems is at least $\Omega(\kappa_\vy^{3/2} \sqrt{\kappa_\vx} )$.
\end{proof}

The above theorems establish lower bounds of $\Omega(\kappa_\vy^{3/2} / \sqrt{\eps})$ and $\Omega(\kappa_\vy^{3/2} \sqrt{\kappa_\vx})$ for C--Q and SC--Q problems, respectively. Our results strictly improve \citet{ji2023lower} by a factor of $\sqrt{\kappa_\vy}$ and achieve a nearly tight characterization for C--Q and SC--Q problems since they match the upper bounds in \citet{ji2023lower} up to logarithmic factors. However, they still have a gap compared with the $\tilde \gO(\bar \kappa_\vy^2 / \sqrt{\eps})$ and $\tilde \gO(\bar \kappa_\vy^2 \sqrt{\kappa_\vx})$ for more general C--SC and SC--SC problems. As all functions in our construction are quadratic, the same lower bounds also apply to objectives with higher-order smoothness. This forms a sharp contrast with NC--SC and NC--Q problems, where higher-order smoothness
can lead to faster convergence as shown by Theorem \ref{thm:highly-smooth-NC-Q}.

\section{Lower Bounds for Stochastic  Bilevel Optimization} \label{subsec:lb-stoc}

For large-scale nonconvex bilevel optimization \citep{ghadimi2018approximation,ji2021bilevel,hong2023two}, a more practical approach is to consider stochastic methods that use a stochastic gradient estimator instead of the deterministic full gradient. In this section, we consider lower bounds for the stochastic first-order (SFO) and stochastic Hessian-vector-product (SHVP) oracles. 

\subsection{Stochastic Oracles and Probability Zero-Chain} \label{subsec:arje-lem}

Formally, SFO is an unbiased estimator of the exact gradient. For bilevel optimization, the algorithm has access to SFO for both the upper- and lower-level objectives $f$ and $g$. 

\begin{dfn}[SFO]
    Let $f$ and $g$ be two continuously differentiable functions.  We define the stochastic first-order oracle (SFO) with variance $\sigma^2>0$
for the bilevel problem $(f,g)$
as the mapping 
\[
\sO^{\rm sfo}(\vx,\vy) = \left(\widehat \nabla f(\vx,\vy), \widehat \nabla {g}(\vx,\vy) \right),
\] 
where $\widehat \nabla f(\vx,\vy)$ and $ \widehat \nabla {g}(\vx,\vy)$ are random vectors such that, {for any possibly random query $(\vx,\vy)$, conditioned on the
$\sigma$-algebra generated by the randomness in $(\vx,\vy)$}, we have
\begin{align*}
    \E \widehat \nabla f(\vx,\vy) =& \nabla f(\vx,\vy), \quad \E \Vert \widehat \nabla f(\vx,\vy) - \nabla  f(\vx,\vy) \Vert^2 \le \sigma^2; \\
    \E \widehat \nabla g(\vx,\vy) =& \nabla g(\vx,\vy), \quad \E \Vert \widehat \nabla g(\vx,\vy) - \nabla  g(\vx,\vy) \Vert^2 \le \sigma^2.
\end{align*}
\end{dfn}

Similarly, we can also define the stochastic Hessian-vector-product (SHVP) oracle as follows.
Although HVP and first-order oracles are equivalent under the deterministic setting (Lemma \ref{lem:zr-grad-hvp}), SHVP oracles are strictly stronger than SFO and can improve the complexity from $\tilde \gO(\epsilon^{-6})$ \citep{chen2023near,kwon2024complexity} to $\tilde \gO(\epsilon^{-4})$ \citep{ji2021bilevel}.

\begin{dfn}[SHVP oracle]  \label{dfn:shvp-oracle}
Let $f$ be differentiable, and $g$ be twice differentiable. We define the stochastic Hessian-vector-product (SHVP) oracle for the bilevel problem $(f,g)$ as the mapping
\begin{align*}
 \sO^{\rm shvp}(\vx,\vy,\vv) = \left({\sO^{\rm sfo}(\vx,\vy)}, \widehat \nabla_{\vx \vy}^2 g(\vx,\vy) \vv, \widehat \nabla_{\vy \vy}^2 g(\vx,\vy) \vv \right),
\end{align*}
where $\widehat \nabla_{\vx \vy}^2 g(\vx,\vy)$ and $ \widehat \nabla_{\vy \vy}^2 g(\vx,\vy)$ are random matrices such that, {for any possibly random query $(\vx,\vy)$, conditioned on the
$\sigma$-algebra generated by the randomness in $(\vx,\vy)$}, we have
\begin{align*}
    \E \widehat \nabla_{\vx \vy}^2 g(\vx,\vy) =& \nabla_{\vx \vy}^2 g(\vx,\vy), \quad \E \Vert \widehat \nabla_{\vx \vy}^2 g(\vx,\vy) - \nabla_{\vx \vy}^2 g(\vx,\vy) \Vert^2 \le \sigma^2; \\
    \E \widehat \nabla_{\vy \vy}^2 g(\vx,\vy) =& \nabla_{\vy \vy}^2 g(\vx,\vy), \quad \E \Vert \widehat \nabla_{\vy \vy}^2 g(\vx,\vy) - \nabla_{\vy \vy}^2 g(\vx,\vy) \Vert^2 \le \sigma^2.
\end{align*}
\end{dfn}
To establish lower bounds in the stochastic setting, we leverage the following notion of probability-$p$ zero-chains \citep{arjevani2023lower,arjevani2020second}. Note that we only set $g$ as a deterministic function in our lower bounds of SHVP methods. Therefore, the following definition only considers the stochasticity of the SFO.

\begin{dfn}[Probability-$p$ zero-chain]
A function $f:\mathcal{X}\to\mathbb{R}$ equipped with an SFO
$\sO:\vx\mapsto \widehat \nabla f(\vx;\xi)$ 
is called a \emph{probability-$p$ zero-chain} if, for any iterate satisfying 
$\operatorname{supp}(\vx)\subset \{1,\ldots,i-1\}$, 
the stochastic gradient obeys
\[
\begin{cases}
\mathbb{P}\!\left(
\operatorname{supp}\big(\widehat \nabla f(\vx;\xi)\big)
\not\subset \{1,\ldots,i-1\}
\right)\le p,\\[4pt]
\mathbb{P}\!\left(
\operatorname{supp}\big(\widehat \nabla f(\vx;\xi)\big)
\subset \{1,\ldots,i\}
\right)=1.
\end{cases}
\]
\end{dfn}
Intuitively, a probability-$p$ zero-chain reveals at most one new coordinate at a time, and the next coordinate becomes accessible with probability at most $p$ at each oracle query. Consequently, the discovery of informative coordinates is significantly slowed down when $p$ is small, leading to stronger stochastic lower bounds than deterministic ones.

\begin{lem}[{\citet[Lemma 1]{arjevani2023lower}}] \label{eq:zr-sfo-tail-lb}
Let $f: \gX \rightarrow \sR$, where $\gX \subseteq \sR^T$ satisfies $\operatorname{supp}(\gP_\gX(\vx)) = \operatorname{supp}(\vx)$ for all $\vx \in \sR^T$, be a probability-$p$ zero-chain with an SFO. For any zero-respecting SFO algorithm that generates the sequence \( \{\vx^{(t)} \}_{t \ge 0}\) such that
\begin{align} \label{eq:zero-respect-sfo}
    \vx^{(t)} \in
\Bigl\{
\gP_{\gX}(\vu):
{\rm supp}(\vu)
\subseteq
\bigcup_{s<t}
{\rm supp}(\vx^{(s)})\cup 
{\rm supp}(\widehat \nabla f(\vx^{(s)}))
\Bigr\}, \quad \vx^{(0)} = \vzero,
\end{align}
with probability at least $1-\gamma$, we have $\vx^{(t)}_T = 0$
 for all $t \le (T -\log(1/\gamma)) / (2p)$.
\end{lem}

\citet{arjevani2023lower} showed that there is also a general mechanism for converting a deterministic zero-chain into a probability-$p$ zero-chain. We recall their result as follows.

\begin{lem}[{\citet[Lemma 3]{arjevani2023lower}}]\label{le:pchain}
Let $f:\mathcal{X}\to\mathbb{R}$ be a zero-chain on $\mathcal{X}\subseteq \mathbb{R}^T$.
For any $\vx\in\mathcal{X}$, let 
$i^*(\vx):=\inf\{\,i\in[T]:x_i=0\,\}$ to be the next coordinate that can be activated.
For any $p\in(0,1]$, define \([\widehat \nabla f(\vx;\xi)]_i\), the $i$th coordinate of the SFO as
\[
[\widehat \nabla f(\vx;\xi)]_i :=
\begin{cases}
\dfrac{\xi}{p} \cdot \dfrac{\partial f(\vx)}{\partial x_i}  , & \text{if } i=i^*(\vx),\\[6pt]
\dfrac{\partial f(\vx)}{\partial x_i}, & \text{otherwise},
\end{cases}
\]
where $\xi\sim\mathrm{Bernoulli}(p)$.
Suppose that there exists $B<\infty$ such that 
$\|\nabla f(\vx)\|_\infty \le B,
\; \forall\,\vx\in\mathcal{X}$.
Then the oracle 
$\sO:\vx\mapsto \widehat \nabla f(\vx;\xi)$ 
is an unbiased SFO with variance bounded by 
$\sigma^2 \le B^2/p$. 
Moreover, $f$ equipped with $\sO$ forms a probability-$p$ zero-chain.
\end{lem}

Based on these lemmas, we construct probability zero-chains for bilevel optimization to establish lower bounds for SHVP and SFO methods in the following subsections.

\subsection{Lower Bounds for SHVP Methods}

We first present a simple warm-up construction obtained by directly embedding the nonconvex zero-chain of \citet{arjevani2020second,arjevani2023lower} into a bilevel problem, which can lead to an $\Omega(\kappa_\vy^4 \eps^{-4})$ lower bound for the following SHVP-based algorithms. 

\begin{dfn}[SHVP-based algorithm] \label{dfn:rand-alg-shvp}
An SHVP-based algorithm $\texttt{A}$ consists of a sequence of measurable mappings $\{ \texttt{A}^{(t)}\}_{t \in \sN}$ such that $ \texttt{A}^{(t)}$ takes in the first $t-1$ oracle responses to produce the $t$-th query $\vw^{(t)} = (\vx^{(t)},\vy^{(t)}, \vv^{(t)})$, such that $\vw^{(0)} = (\vzero,\vzero, \vzero)$ and  
\begin{align}  \label{eq:rand-seq-shvp}
    \vw^{(t)} = \texttt{A}^{(t)} \left(\sO^{\rm shvp}(\vw^{(0)}), \cdots, \sO^{\rm shvp}(\vw^{(t-1)}) \right), \quad \forall t \in \sN_+.
\end{align}                                                         
We denote $\gA^{{\rm shvp}}$ as the set of all stochastic algorithms that follow protocol (\ref{eq:rand-seq-shvp}).
\end{dfn}

\citet{arjevani2020second,arjevani2023lower} introduced the following stochastic nonconvex zero-chain to show the optimal $\Omega(\epsilon^{-4})$ lower bound for single-level nonconvex optimization. We recall their construction as follows.

\begin{dfn}[Stochastic NC zero-chain \citep{arjevani2023lower}] \label{dfn:Arj-NC-chain-stoc}
Let 
{\begin{align*}
    \widehat f^{\text{snc}}(\vx) &=  \left( -\Psi(1) \Phi(\vx_1) + \sum_{i=2}^T [\Psi(-\vx_{i-1}) \Phi(-\vx_i) - \Psi(\vx_{i-1} ) \Phi(\vx_i)] \right).
\end{align*}}
where the component functions $\Psi, \Phi: \sR \rightarrow \sR$ are 
\begin{align*}
    \Psi(x) =  
    \begin{cases}
        \exp \left( 1 -\dfrac{1}{(2x-1)^2} \right), & x > 1/2;\\
        0, & x\le 1/2
    \end{cases}
     \quad {\rm and} \quad \Phi(x) = \sqrt{\rm e} \int_{-\infty}^x \exp(- t^2/2 ) {\rm d}t.
\end{align*}
\end{dfn}
The component functions $\Psi, \Phi: \sR \rightarrow \sR$ are first introduced by \citet{carmon2020lower} and we recall their properties as follows.
\begin{lem}[{\citet[Lemma 1]{carmon2020lower}}]
\label{le:phipsi}
The functions $\Psi$ and $\Phi$ satisfy:
\begin{enumerate}
    \item For all $x\le 1/2$, we have 
        $\Psi^{'}(x)=0$.
    \item For all $x\ge 1$ and $|x'|<1$, we have 
        $\Psi(x)\Phi'(x')>1$.
    \item For all $x\in \sR$, we have
        $|\Psi'(x)|
        \le
        \exp\!\left(2.5\log (4) \right)
     $ and $ 
        |\Phi'(x)|
        \le
        \exp\!\left(1.5\log (1.5)\right)$.
    \item The functions $\Psi$, $\Phi$ and their derivatives $\Psi'$,  and $\Phi'$ satisfy: 
        $0<\Psi<{\rm e},\;
        0<\Psi'<\sqrt{{54}/{\rm e}}$, and $
        0<\Phi<\sqrt{2\pi {\rm e}},
        \; 0<\Phi'<\sqrt {\rm e}$. 
\end{enumerate}
\end{lem}

Based on the above properties of component functions, \citet{arjevani2023lower} showed the following property of the nonconvex zero-chain $\widehat f^{\rm snc}$.
\begin{lem}[{\citet[Lemma 6 and 7]{arjevani2023lower}}] \label{lem:Arj-random}
The nonconvex zero-chain function
$\widehat f^{\text{snc}}: \sR^T \rightarrow \sR$ satisfy:
\begin{enumerate}
    \item We have $\widehat f^{\text{snc}}(\vzero)-  \inf_{\vx \in \sR^T} \widehat f^{\text{snc}}(\vx) \le 12 T $;
    \item $\widehat f^{\text{snc}}(\vx)$ has $152$-Lipschitz continuous gradients;
    \item  For all $\vx \in \sR^T$, we have $\Vert \nabla \widehat f^{\rm snc}(\vx) \Vert_\infty \le 23$.
    \item  For all $\vx \in \sR^T$ such that $\vert x_T  \vert <1$, then there exists $1 \le j \le T$ such that \(| x_{j-1} | \ge 1\) , \(|x_j| < 1\), and \( \| \nabla \widehat f^{\rm snc}(\vx) \| \ge | \partial \widehat f^{\rm snc}(\vx) / \partial x_j | \ge 1 \).
\end{enumerate}
\end{lem}
These properties are sufficient to show the optimal $\Omega(\epsilon^{-4})$ lower bound for nonconvex stochastic optimization \citep{arjevani2023lower}. In the following, we show that by leveraging the bilevel structure and setting the lower-level function with a simple quadratic function to rescale the upper-level variable from $\vx$ to $\kappa_\vy \vx$, it is equivalent to rescaling the stationarity threshold from $\eps$ to $\eps/\kappa_\vy$ and resulting in an $\Omega(\kappa_\vy^4 \eps^{-4})$ lower bound for stochastic NC--SC bilevel problems. Formally, we set up the upper- and lower-level functions $f: \sR^T \rightarrow \sR$ and $g: \sR^T \times \sR^T \rightarrow \sR$ as:
\begin{align} \label{eq:our-fg-stoc}
\begin{split}
      f(\vy)&: = \frac{L_1 \beta^2}{152} \widehat f^{\text{snc}}( \vy / \beta ) \quad {\rm and} \quad 
    g(\vx,\vy): = \frac{\mu_\vy}{2} \Vert \vy \Vert^2 - L_1 \langle \vx, \vy \rangle.
\end{split}
\end{align}
Since the lower-level objective is quadratic and $\mu_\vy$-strongly convex in $\vy$, its minimizer is \(\vy^*(\vx) = \kappa_\vy \vx \) for $\kappa_\vy = L_1/ \mu_\vy $. Substituting $\vy^*(\vx)$ into $f$ gives the hyper-objective as
\begin{align} \label{eq:scaled-phi-stoc}
   F(\vx) =  f(\vy^*(\vx)) = \frac{L_1 \beta^2}{152} \widehat f^{\rm snc}(\kappa_\vy \vx / \beta).
\end{align}
This rescaling directly leads to the following result.

\begin{thm}[SHVP NC--SC lower bound ] \label{thm:NC-SC-stoc}
There exists a constant $b_0 \in \sR_+$, such that, for any $L_1 \ge \mu_\vy > 0$, $\sigma, \Delta >0$, and $L_0 \ge b_0\sqrt{L_1 \Delta}$, there exists a smooth nonconvex--quadratic bilevel problem $(f, g) \in \gF^{\text{nc--sc}}(L_0, L_1,\mu_\vy, \Delta)$ such that, the complexity of any SHVP algorithm $\texttt{A} \in \gA^{\rm shvp}$ to find an $\epsilon$-stationary point of the hyper-objective $F(\vx)$ with probability $0.5$ is at least
\begin{align*}
    \Omega(\kappa_\vy^2 L_1 \Delta \epsilon^{-2} + \kappa_\vy^4\sigma^2  L_1  \Delta \epsilon^{-4}),
\end{align*}
where $\kappa_\vy = L_1/\mu_\vy$ is the condition number of the lower-level problem.
\end{thm}

\begin{proof}
It is not hard to show that Lemma \ref{lem:our-alg-class-reduction} also holds for stochastic oracles $\widehat \nabla f$ when $\gX = \sR^d$. Therefore, it is sufficient to show lower bounds for the zero-respecting SFO algorithms defined as Eq. (\ref{eq:zero-respect-sfo}). For these algorithms, by Lemma \ref{eq:zr-sfo-tail-lb}, for all $t \le (T - \log 2 ) /(2p)$, we have $\vx^{(t)}_T =0$ with probability $0.5$, and thus $\Vert \nabla_\vx f(\vy^*(\vx)) \Vert \ge L_1 \beta \kappa_\vy/152$. Therefore, we can set
\begin{align} \label{eq:para-beta}
    \beta = \frac{152 \epsilon}{L_1 \kappa_\vy}.
\end{align}
to obtain an $\Omega( T/p )$ lower bound for finding an $\epsilon$-stationary point of $f(\vy^*(\vx))$. Next, by Lemma~\ref{lem:Arj-random} (item 1), we have $f(\vy^*(\vx)) - \inf_{\vx \in \sR^T} f(\vy^*(\vx)) \le 12 L_1 \beta^2 T / 152 $. Hence, we can fulfill the constraint $F(\vzero) - \inf_{\vx \in \sR^T} F(\vx) \le \Delta$ by setting
\begin{align} \label{eq:para-T-stoc}
      T = \left \lfloor 
     \frac{152 \Delta}{12 L_1 \beta^2} 
    \right \rfloor = \Theta \left( \frac{\Delta \kappa_\vy^2 L_1}{\epsilon^2 }   \right).
\end{align}
We return the exact gradients of $g$ and perturbed stochastic gradients of $f$ constructed by Lemma \ref{le:pchain}. Then, by combining Lemma \ref{lem:Arj-random} (item 3), we know that the variance of SFO $\widehat \nabla f(\vy)$ is upper bounded by $ (23 L_1 \beta / 152)^2 /p $. Therefore, we can let the variance be bounded by $\sigma^2$ by setting
\begin{align} \label{eq:para-p-stoc}
    p = \min \left\{  1, \left( \frac{23 L_1 \beta}{\sigma 152} \right)^2 \right\}.
\end{align}
Combining all the parameters, we can conclude the lower bound of 
\begin{align*}
    \Omega\left( \frac{T}{p} \right) = \Omega \left( \frac{\Delta}{L_1 \beta^2} \left(1+ \frac{\sigma^2}{L_1^2 \beta^2} \right)  \right) = \Omega \left( \frac{\kappa_\vy^2 L_1 \Delta }{\epsilon^2} \left( 1+ \frac{\kappa_\vy^2 \sigma^2}{\epsilon^2} \right)\right).
\end{align*}
Finally, since $\Vert \widehat \nabla f^{\rm snc}(\vx) \Vert_\infty \le 23$ from Lemma \ref{lem:Arj-random} (item 3), we have $\Vert \widehat f^{\rm snc}(\vx) \Vert_\infty \le 23 \sqrt{T}$ and thus $\Vert \nabla f(\vy) \Vert \le 23 L_1 \beta \sqrt{T} / 152$. Recalling that $T = \Theta(\Delta / (L_1 \beta^2)$, we know that the upper-level function is $L_0$-Lipschitz continuous if $L_0 = \Omega(\sqrt{L_1 \Delta})$. Hence, the bilevel problem \((f,g)\) we constructed lies in the class $\gF^\text{nc--sc}(L_0,L_1,\mu_\vy,\Delta)$.
\end{proof}

The $\Omega(\kappa_\vy^4 \eps^{-4})$ lower bound for stochastic NC--SC bilevel problems shown above has much higher condition number dependency compared with both the $\Omega(\kappa_\vy^{1/3} \eps^{-4})$ lower bound \citep{li2021complexity} or the $\tilde \gO(\kappa_\vy \eps^{-4})$ upper bound \citep{zhang2022sapd+} for NC--SC minimax optimization problems. This result confirms our previous findings in the deterministic setting that bilevel optimization is more challenging than minimax optimization. In particular, for stochastic problems, a large problem complexity comes from the larger variance for hyper-gradient estimation. However, compared to the currently best  $\tilde \gO(\bar \kappa_\vy^9 \eps^{-4})$ upper bound \citep{ji2021bilevel}, a large gap exists in the condition number dependency.

{ \subsection{Lower Bounds for SFO Methods}}

The $\Omega(\kappa_\vy^4\eps^{-4})$ lower bound shown in the previous section holds for both SHVP and SFO algorithms. However, the previous construction does not leverage the nested chain technique introduced in Section \ref{subsec:nested-chain}. In this section, we further use this technique to show an improved lower bound of $\Omega(\kappa_\vy^{9/2} \eps^{-4})$ for SFO methods. At a high level, we set up a nested chain like Construction (\ref{eq:our-class-det}) in Section \ref{subsec:lb-det}, and perturb the lower-level function in the stochastic gradients. To satisfy the bounded variance constraint, we follow \citet{li2021complexity} to restrict the domains of $\vx$ and $\widetilde \vy$ to the following bounded hypercubes:
\begin{align}
\begin{split}
\gX =& \gC_{\beta R_x/ \alpha}^{T} := \{\vx\in \mathbb{R}^{T} : \|\vx\|_\infty \le \beta R_x / \alpha \},\\
\gY =& \gC_{\beta R_y}^{nT} := \{\widetilde \vy\in \mathbb{R}^{nT} : \|\widetilde \vy\|_\infty \le \beta R_y \}.   
\end{split}
\end{align}
Therefore, our hard instance to prove lower bounds for SFO methods is a constrained bilevel problem. In light of Lemma \ref{lem:our-alg-class-reduction}, we restrict the algorithm class to the following zero-respecting algorithm class on the domain $\gX \times \gY$.
\begin{dfn}\label{dfn:zr-sfo-constraint}
Let $f,g: \gX\times\gY\rightarrow\mathbb{R}$ be the upper- and lower-level functions. Given SFO for $f$ and $g$, a zero-respecting SFO algorithm $\texttt{A}$ is initialized at \(\vz^{(0)} = \vzero\), and generates
the iterate $\vz^{(t)} =  (\vx^{(t)},\vy^{(t)})$ such that for $t \ge 1$ and $
\gZ = \gX \times \gY$ satisfy 
\[
\begin{aligned}
\vz^{(t)}
&\in
\Bigl\{
\gP_{\gZ}(\vu):
{\rm supp}(\vu)
\subseteq
\bigcup_{s <t}
{\rm supp}(\vz^{(s)})\cup 
{\rm supp}(\widehat \nabla f(\vz^{(s)})) \cup {\rm supp}(\widehat \nabla g(\vz^{(s)}))
\Bigr\}.
\end{aligned}
\]
We denote $\gA^\text{zr--sfo}$ as the set of all zero-respecting SFO algorithms that follow this protocol.
\end{dfn}
It is easy to see that the above definition is the generalization of the single-level zero-respecting algorithms in Eq. (\ref{eq:zero-respect-sfo}) to bilevel problems. Since we rely on a constrained problem to show the lower bounds, we also need to generalize the previous notion of hyper-stationarity to constrained settings using the projected gradient as follows.

\begin{dfn}[{\citet[Definition 2.2.3]{nesterov2018lectures}}] \label{dfn:projected-grad}
Let $\gX \subseteq \sR^T$ be a convex compact set and $f: \gX \rightarrow  \sR $ be a continuously differentiable function whose gradients are $L$-Lipschitz continuous. We define the projected gradient as
\[
\gG_{f,\gX}(\vx) = L (\vx - \gP_\gX ( \vx - (1/L) \nabla f(\vx) )).
\]
\end{dfn}
It can be used as a natural stationarity measure for constrained problems \citep{nesterov2018lectures}, and we can also extend the $\eps$-stationarity points in Definition \ref{dfn:fo-sta} to $\Vert \gG_{F,\gX} (\vx) \Vert \le \epsilon$, where $F = f(\vx,\vy^*(\vx))$ is the hyper-objective function for bilevel optimization. 

Now, we are ready to give our construction for stochastic NC--SC bilevel problems and show our lower bound. Same as the nested chain (\ref{eq:our-class-det}) for deterministic problems, we let $\vx=[x_1,...,x_{T}]\in\mathbb{R}^{T}$ be the upper-level variable, and the block vector $\widetilde \vy = [\vy^{(1)},\vy^{(2)},....,\vy^{(T)}]$ with each $\vy^{(i)}\in\mathbb{R}^n$ be the lower-level variable. In addition, we set the lower-level function the same as Eq. (\ref{eq:our-class-det}) while constrained on the hypercube $\gY = \gC_{\beta R_y}^{nT}$. For the upper-level function, we replace the quadratic links in deterministic construction (\ref{eq:our-class-det}) with the link function $\Lambda(y_+,y_-) = -\Psi_\beta(y_+) \Phi_\beta(2y_-) + \Psi_\beta(-y_+) \Phi(-2y_-)$, where $\Psi_\beta$ and $\Phi_\beta$ are the shorthands for the scaled functions $\beta \Psi(x / \beta)$ and $\beta \Phi(x / \beta)$, respectively. Formally, we construct
\begin{align}  \label{eq:our-fg-sfo}
\begin{split}
    &f(\tilde \vy) = -\delta \Psi_\beta(\alpha) \Phi_\beta(2y_1^{(1)}) + \delta \sum_{i=1}^{T-1}[- \Psi_\beta(y_n^{(i)} ) \Phi_\beta(2y_1^{(i+1)}) + \Psi_\beta(-y_n^{(i)}) \Phi_\beta(-2y_1^{(i+1)}) ], \\
    &g(\vx,\widetilde \vy) =  \delta \sum_{i=1}^T \left\{ g_n(x_i,\vy^{(i)}) \triangleq \frac{1}{2} (\vy^{(i)})^\top \mA_n \vy^{(i)} -  x y^{(i)}_n \right\}, 
\end{split}
\end{align}
where $\mA_n \in \sR^{n \times n}$ is the inverse-endpoint-amplified tridigagonal matrix defined in Eq. (\ref{eq:matrix-Bk}) and $\Psi,\Phi: \sR \rightarrow \sR$ are the one-dimensional component functions in Definition \ref{dfn:Arj-NC-chain-stoc}. Before moving to subsequent analyses, we first verify that the unconstrained minimizer $\widetilde \vy^*(\vx)$ lies in the domain $\gY = \gC_{\beta R_y}^{nT}$. Then we can apply Lemma \ref{lem:matrix-A}(a) to obtain the explicit form of the hyper-objective function $F(\vx) = f(\widetilde \vy^*(\vx))$.

\begin{lem} \label{lem:yindomain}
Let $\alpha = 2 (n-1)^2/3$. If $R_y \ge 3 R_x$, then for any $\vx$ lying in the hypercube $\gC_{\beta R_x/ \alpha}^T$, the lower-level minimizer $\widetilde \vy^*(x) = \arg \min_{\widetilde \vy \in (\sR^{n})^T} g(\vx,\vy)$ lies in the hypercube $\gC_{\beta R_y}^{nT}$.
\end{lem}

\begin{proof} Since the function $g(\vx,\widetilde \vy) = \delta \sum_{i=1}^T g_n(x_i,\vy^{(i)})$ are separable in each $i \in [T]$. We only need to show that for any $x \in \gC^1_{\beta R_x/ \alpha}$, the minimizer $\vy^*(\vx) = \arg \min_{\vy \in \sR^n} g_n(x,\vy) $ lies in $\gC^n_{\beta R_y}$. We know from Lemma \ref{lem:matrix-A}(b) that
\[
\Vert \mA_n \Vert^{-1} = \frac{1}{\lambda_{\min}(\mA_n)}\le 2(n-1)^2.
\]
Then, for each coordinate of the lower-level solution $\vy^*(x) = x (\mA_n^{-1} \ve_n)$, we have
\[
y_i^*(x) = x (\mA_n^{-1})_{i,n} = x \langle \ve_i , \mA_n^{-1} \ve_n \rangle \le | x| \cdot \Vert \mA_n^{-1} \Vert \le | x | \cdot 2(n-1)^2.
\]
Hence, for $\alpha =  2(n-1)^2/3$ we have $x \in \gC_{\beta R_x/ \alpha}^1$, and for $R_y \ge 3R_x$ we have $\vy^*(x) \in \gC_{\beta R_y}^n$.
\end{proof}

In addition, Lemma \ref{lem:matrix-A}(a) gives the explicit endpoint values of $\mA_n^{-1}$, which indicates that $\vy_n^{(i)} = y_1^{(i)}= \alpha x_i$. Then, by substituting $\vy^*(\vx)$ into $f$ gives the hyper-objective as
\begin{align} \label{eq:hyper-stoc-form}
    F(\vx) = f(\widetilde \vy^*(\vx)) = \delta \beta^2 \widehat f^{\rm snc}(\alpha \vx / \beta),
\end{align}
where \(\widehat f^{\rm snc} : \sR^T \rightarrow \sR \) is the stochastic nonconvex zero-chain defined in Definition \ref{dfn:Arj-NC-chain-stoc}. Next, we show that the projected gradient in Definition \ref{dfn:projected-grad} for $F(\vx)$ is well-defined since the hyper-objective function has Lipschitz continuous gradients. Its projected gradient norm is lower-bounded when the algorithm has not yet reached the end of the chain.

\begin{lem} \label{lem:hyper-smooth}
Let $\alpha = 2 (n-1)^2/3$. If $R_y \ge 3 R_x$ and $ R_x \ge 2$, then our construction satisfies
\begin{enumerate}
    \item  The hyper-objective $F(\vx) = f(\widetilde \vy^*(\vx))$ has $152 \delta \alpha^2$-Lipschitz continuous gradients.
\item For any $\vx \in \sR^T$ such that $x_{T}= 0$, we have $\Vert \gG_{F,\gX}(\vx) \Vert \ge  \delta \alpha \beta$.
\end{enumerate}
\end{lem}

\begin{proof}
Item 1 directly follows from Eq. (\ref{eq:hyper-stoc-form}) that $F(\vx) = \delta \beta^2 \widehat f^{\rm snc}(\alpha \vx / \beta)$ and Lemma~\ref{lem:Arj-random} (item 2) that the gradients of the function $\widehat f^{\rm snc}$ are $152$-Lipschitz continuous. Then, the projected gradient in Definition \ref{dfn:projected-grad} is well-defined, and we provide the lower bound of its norm when  $x_{T}= 0$ as follows.
Let $L = 152 \delta \alpha^2$ be the gradient Lipschitz constant in item~1. By Lemma~\ref{lem:Arj-random} (item 4), there exists $1 \le j\le T$ such that \(| x_{j-1}  | \ge \beta /\alpha\) , \(|x_j| < \beta /\alpha\), and \( | \partial F(\vx) / \partial x_j | \ge \delta \alpha \beta \). Recalling that $\gX = \gC_{\beta R_x/ \alpha}^{ T}$, if \( | x_j-(1/L) F(\vx) / \partial x_j |< R_x\), we have 
\[
\| \gG_{F,\gX}(\vx) \| \ge L \left | x_j -  \gP_{\gC_{\beta R_x/ \alpha}^1}[x_j-(1/ L){\partial F(\vx)}/{\partial x_j}]  \right | =
\left | \partial  F(\vx) / \partial x_j \right | \ge \delta \alpha \beta.
\]
Otherwise, when $R_x >1$, if  \( |x_j-(1/L) \partial F(\vx) / \partial x_j | \ge R_x\), we have
\[
\|  \gG_{F,\gX}(\vx) \| \ge 
L \left | x_j -  \gP_{\gC_{\beta R_x/ \alpha}^1}[x_j-(1/ L){\partial F(\vx)}/{\partial x_j}]  \right | \ge L \left( \frac{\beta R_x}{\alpha} - |x_j | \right) \ge \frac{L \beta (R_x-1)}{\alpha},
\]
which is no smaller than $\delta \alpha \beta$ since $L = 152 \delta \alpha^2$ and $R_x \ge 2$.
\end{proof}

Based on all the above auxiliary lemmas, we prove our main theorem as follows.
\begin{thm}[SFO NC--SC lower bound] \label{thm:stoc-NC-Q}
There exist two numerical constants $a_0 \in (0,1)$ and $b_0 \in \sR_+$ such that, for any $L_1, \sigma, \Delta >0$,  $\mu_\vy \in (0, a_0 L_1]$, and $L_0 \ge b_0 \sqrt{L_1 \Delta}$, there exists a smooth nonconvex--quadratic bilevel problem $(f, g) \in \gF^{\text{nc--sc}}(L_0,L_1,\mu_\vy, \Delta)$ on the domain $\gX \times \gY$ such that, the complexity of any zero-respecting SFO algorithm $\texttt{A} \in \gA^\text{zr-sfo}$ to find an $\epsilon$-stationary point of the  hyper-objective $F(\vx)$ such hat $\Vert \gG_{F,\gX} (\vx) \Vert \le \eps $ is at least
\begin{align*}
    \Omega(\kappa_\vy^{5/2} L_1 \Delta \epsilon^{-2} + \kappa_\vy^{9/2} \sigma^2  L_1  \Delta \epsilon^{-4}),
\end{align*}
where $\kappa_\vy = L_1/\mu_\vy$ is the condition number of the lower-level problem.
\end{thm}

\begin{proof} We construct the hard instance as Eq. (\ref{eq:our-fg-sfo}). As for the SFO, we only perturb $\widehat \nabla_{\widetilde \vy} g(\vx,\widetilde \vy)$, the partial gradient of the lower-level function with respect to $\vy$, while returning the exact gradients for the others, including $\widehat \nabla f(\widetilde \vy)$ and $\widehat \nabla_\vx g(\vx,\widetilde \vy)$. Let $R_x = 2$ and $R_y = 6$ such that the conditions in Lemma \ref{lem:hyper-smooth} are satisfied. Then, we know that in the domain \(x \in \gC_{\beta R_x/\alpha}^1 \) and \(\vy \in \gC^n_{\beta R_y}\), the partial gradient norm of each $g_n (x,\vy)$ in $\vy$ is bounded by:
\[
\|\nabla_\vy g_n(x,\vy) \| \le \beta \Vert \mA_n \Vert  R_y + \frac{\beta R_x \Vert \ve_n \Vert}{\alpha} \le 5.5 \beta R_y,
\]
where the last step uses Lemma \ref{lem:matrix-A}(b) that $\Vert \mA_n \Vert \le 5$. Then, constructing the SFO according to Lemma \ref{le:pchain}, the variance of $\widehat \nabla_{\widetilde \vy} g(\vx,\widetilde \vy)$ is bounded by $5.5^2 \beta^2 \delta^2 R_y^2 / p$.
Consequently, to ensure that the SFO variance is bounded by $\sigma^2$, it suffices to choose
\begin{align}
\label{eq:valuep}
p
=
\min\left\{
1,\,
\frac{5.5 \beta^2 \delta^2 R_y^2}{\sigma^2}
\right\}.
\end{align}
Note that by simply perturbing the lower-level variable $\vy$, we can still turn the zero-chain with length $nT$ to a probability-$p$ zero-chain with length $(n-1)T$, which corresponds to the number of edges in Figure \ref{fig:nested-chain} that are connected by $\nabla_{\widetilde \vy} g(\vx, \widetilde \vy)$. Then, by Lemma \ref{le:pchain}, with probability at least $0.5$, any zero-respecting algorithm satisfy $\vx_T^{(t)} = 0$ whenever
\begin{align}
\label{eq:tzero}
t
\le
\frac{(n-1)(T-1)-\log 2}{2p} = \Omega\left( \frac{nT}{p} \right).
\end{align}
Next, we properly set the values of $n$ and $T$ and obtain the lower bound as the right-hand side above. Recall that Lemma \ref{lem:matrix-A}(b) indicates that the lower-level function $g(\vx,\widetilde \vy)$ is $\mu_\vy$-strongly convex under the setting $n = \Theta(\sqrt{\kappa_\vy})$.
It leads to the scaling $\alpha = \Theta(\kappa_\vy)$ since $\alpha = 2 (n-1)^2/ 3$.
For the hyper-objective $f(\widetilde \vy^*(\vx))$, Lemma \ref{lem:hyper-smooth} suggests that 
\[
\vx_T = 0 \quad \Longrightarrow \quad
\Vert \gG_{F,\gX}(\vx) \Vert \geq \delta \alpha \beta.
\]
Therefore, whenever Eq. (\ref{eq:tzero}) holds, we have \(\vx_T^{(t)} = 0\) for $t = o(nT/p)$, and thus $\Vert \gG_{F,\gX}(\vx^{(t)}) \Vert \ge \eps$ under the parameter setting
\[
\beta = \frac{\eps}{\alpha \delta} = \Theta \left(\frac{\eps}{\kappa_\vy \delta} \right).
\]
From Lemma \ref{le:phipsi} and Lemma \ref{lem:matrix-A}, we know that both the nonconvex upper-level function $f$ and the strongly-convex lower-level function $g$ have $L_1$-Lipschitz continuous gradients under the parameter setting $\delta = \Theta(L_1)$. 
It remains to select parameters $T$ as large as possible, such that the function pair $(f,g) \in {\mathcal{F}}^\text{nc--sc}(L_0,L_1,\mu_\vy,\Delta)$. By Lemma~\ref{lem:Arj-random} (item 1), we have \(f(\widetilde \vy^*(\vzero)) - \inf_{\vx \in \sR^T}f(\widetilde \vy^*(\vx)) \le 12 \delta \beta^2 T \), which is bounded by $\Delta$ under the setting 
\[
T = \left \lfloor \frac{\Delta}{12 \delta \beta^2}\right \rfloor = \Theta \left( \frac{\Delta L_1 \kappa_\vy^2}{\eps^{2}} \right).
\] 
In addition, Lemma \ref{le:phipsi} (item 4) indicates that the upper-level function $f(\widetilde \vy)$ satisfies that $\vert \partial f(\widetilde \vy) / \partial \vy_1^{(i)} \vert, \vert \partial f(\widetilde \vy) / \partial \vy_n^{(i)} \vert = \gO(\delta \beta)$ for all $i \in [T]$.
Consequently, we have 
\[
\Vert \nabla f(\widetilde \vy) \Vert = \gO\left(\delta \beta \sqrt{T} \right) = \Theta \left( \sqrt{L_1 \Delta} \right),
\]
which is bounded by $L_0$ if $L_0 \ge b_0 \sqrt{L_1 \Delta}$ for a sufficiently large constant $b_0 \in \sR_+$. It completes the proof of  $(f,g) \in {\mathcal{F}}^\text{nc--q}(L_0,L_1,\mu_\vy,\Delta)$. 
Putting all the pieces together, it leads to the lower bound of 
\[ 
\Omega \left( \frac{nT}{p} \right) 
= \Omega \left( \frac{\sqrt{\kappa_\vy}\Delta}{L_1 \beta^2} \left(1 + \frac{\sigma^2}{L_1^2\beta^2} \right)  \right) = 
\Omega \left( \frac{\kappa_\vy^{5/2} L_1 \Delta}{\epsilon^2} \left( 1 + \frac{\kappa_\vy^2 \sigma^2}{\eps^2} \right)  \right).
\]
\end{proof}
The $\Omega(\kappa_\vy^{9/2} \epsilon^{-4})$ shown in the theorem improved the $\Omega(\kappa_\vy^4 \eps^{-4})$ of Theorem \ref{thm:NC-SC-stoc} by a factor of $\sqrt{\kappa_\vy}$ by leveraging the nested chain technique and the nature of SFO algorithms. However, compared with the $\tilde \gO(\bar \kappa_\vy^{11} \eps^{-6})$ upper bound of SFO methods \citep{chen2023near,kwon2024complexity}, significant gaps in both $\kappa_\vy$ and $\eps$ persist. 
Although \citet{kwon2024complexity} showed that the $\eps$ dependency in second term of the lower bound can be improved to $\Omega(\eps^{-6})$ under $\vy^*$-aware oracle model, this oracle is only unbiased when $\vy$ lies in an $\gO(\eps)$ neighborhood of $\vy^*(\vx)$ and are adversarial outside. Therefore, a lower bound of $\Omega(\eps^{-6})$ under standard SFO remains open, which is also explicitly mentioned by \citet[Conjecture~1]{kwon2024complexity}.

\section{Conclusion and Future Work}

In this work, we introduce a novel nested chain framework to amplify the effective length for zero-chains and apply it to prove new lower bounds for bilevel optimization in condition number dependency. For deterministic NC--SC problems, we show lower bounds of $\Omega(\kappa_\vy^{5/2} \epsilon^{-2})$, $\Omega(\kappa_\vy^{9/4} \epsilon^{-7/4})$, and $\Omega(\kappa_\vy^{13/6} \epsilon^{-5/3})$ for first-order, second-order, and arbitrarily smooth problems, respectively. We also extended our lower bound to $\Omega(\kappa_\vy^{3/2} /\sqrt{\epsilon})$ for deterministic C--SC problems, and to $\Omega(\kappa_\vy^4 \epsilon^{-4})$ and $\Omega(\kappa_\vy^{9/2} \epsilon^{-4})$ for stochastic NC--SC problems for SHVP and SFO methods, respectively. We also show that our lower bounds are tight (up to logarithmic factors) for deterministic problems with quadratic lower-level functions. However, a significant gap still exists for general strongly convex lower-level problems.

Several important questions remain open. The most immediate challenge is to obtain a similar characterization for the general setting mentioned above by leveraging non-quadratic lower-level functions to show higher lower bounds. It would also be interesting to study other setups, such as 
mean-squared smooth problems \citep{khanduri2021near,yang2021provably, dagreou2022framework,yang2023achieving,chu2024spaba}, 
structured NC--NC problems \citep{kwon2024penalty,chen2024finding,jiang2025beyond}, and zeroth-order bilevel optimization \citep{aghasi2025fully,aghasi2025optimal} in the future.


\section*{Acknowledgments}

The authors greatly thank Jeongyeol Kwon for his helpful comments on the preliminary version of this manuscript. The technical Lemmas \ref{lem:argmin-close-form}, \ref{lem:rela-stationarity},  and \ref{lem:matrix-A}) are discovered under assistance from GPT 5.5 Pro and then verified by the authors. 

\bibliographystyle{plainnat}
\bibliography{sample}

\appendix

\newpage
\section{Examples with Lower-Level Quadratics} \label{apx:example-Q}

This section provides a collection of (S)C--SC and NC--SC bilevel optimization applications in which the lower-level function is quadratic. 

\subsection{Convex--Strongly-Convex Examples} \label{apx:exmple-C-SC}
The (S)C--SC setting is often assumed without any examples in the literature \citep{ghadimi2018approximation,hong2023two}. To the best of our knowledge, the only example is a fair resource allocation problem~\citep{srikant2014communication} given by \citet{ji2023lower}. In the following, we provide several machine learning problems that provably satisfy the assumptions, which verify the rationality of assuming the (strong) convexity of the hyper-objective in previous work \citep{ghadimi2018approximation,hong2023two}.

Let $(\mA,\vb)$ denote a dataset with $n$ samples, where each row of $\mA \in \sR^{n \times d}$ represents a $d$-dimensional feature of a sample, and the element in $\vb \in \sR^n$ represents the corresponding label. We consider the following examples.

\begin{exmp}
Given $T$ different tasks, with each task $i$ corresponding to a training set $(\mA_i^{\rm tr}, \vb_i^{\rm tr})$ and  a testing set $(\mA_i^{\rm test}, \vb_i^{\rm  test})$. Define the covariance matrices of each task as $\mG_i^{\rm test}  =\left(\mA_i^{\rm test}  \right)^\top \mA_i^{\rm test} $ and $\mG_i^{\rm tr}  =\left(\mA_i^{\rm tr}  \right)^\top \mA_i^{\rm tr}$. 
We consider the meta-learning for linear regression~\citep{rajeswaran2019meta}:
    \begin{align*}
        &\quad \min_{\vx \in \sR^d} \frac{1}{T} \sum_{i=1}^T \frac{1}{2} \Vert \mA_i^{\rm tr} \vy_i^*(\vx) - \vb_i^{\rm tr} \Vert^2, \\
        &{\rm s.t.} \quad \vy_i^*(\vx) = \arg \min_{\vy' \in \sR^d} \frac{1}{2} \Vert \mA_i^{\rm test} \vy' - \vb_i^{\rm test} \Vert^2 + \frac{\lambda}{2} \Vert \vy' - \vx \Vert^2.
    \end{align*}
In this example, the hyper-objective $\varphi(\vx) = f(\vx,\vy^*(\vx))$ is (strongly) convex.
\end{exmp}

\begin{proof}
The lower-level solution has a closed form
\begin{align*}
    &\quad \vy_i^*(\vx) = (\mG_i^{\rm test} + \lambda \mI_d)^{-1} \left(\lambda \vx + \left(\mA_i^{\rm test} \right)^\top \vb_i^{\rm test} \right).
\end{align*}
Plugging the above equation into the upper-level function gives the closed form of $\varphi(\vx)$, which is a convex quadratic function because its Hessian matrix is positive (semi)-definite:
\begin{align*}
    &\quad \mH = \frac{1}{T} \sum_{i=1}^T \lambda^2 (\mG_i^{\rm test} + \lambda \mI_d)^{-1} \mG_i^{\rm tr} (\mG_i^{\rm test} + \lambda \mI)^{-1} \succeq \mO_d.
\end{align*}.
\end{proof}

\begin{exmp}
Given a dataset $(\mA,\vb)$ with $n$ samples, we  
consider the Stackelberg game \citep{bruckner2011stackelberg} for robust linear regression. In this model, the learner's goal is to learn to good linear predictor $\vx$ with small error, but the label $\vy$ is modified by an adversarial data provider, whose goal is to make the loss (per sample weighted by $w_i^2$) high, with a cost to modify the original label $\vb$. 
Let $\mD = {\rm diag}(w_1,\cdots, w_n)$.
It can be formulated as the following bilevel problem
\begin{align*}
    &\quad \min_{\vx \in \sR^d} \frac{1}{2} \Vert \mA \vx - \vy^*(\vx) \Vert^2 \\
    &{\rm s.t.} \quad \vy^*(\vx) = \arg \min_{\vy' \in \sR^n} - \frac{1}{2} \Vert \mD (\mA \vx - \vy') \Vert^2 + \frac{\lambda}{2} \Vert \mD(\vy'  -\vb) \Vert^2.
\end{align*}
In this example, the hyper-objective $\varphi(\vx) = f(\vx,\vy^*(\vx))$ is (strongly) convex if $\lambda >1$.
\end{exmp}

\begin{proof}
When $\lambda >1$, the lower-level solution has a closed form
\begin{align*}
    \vy^*(\vx) = \frac{\mD  (\lambda \vb - \mA  \vx)}{\lambda -1}.
\end{align*}
Plugging the above equation into the upper-level function gives the closed form of $\varphi(\vx)$, which is a convex quadratic function because since its Hessian matrix is positive (semi)-definite:
\begin{align*}
    \mH = \frac{1}{(\lambda-1)^2}  \mA^\top ( (\lambda -1)^2 \mI_n + \mD_n^2 ) \mA \succeq \mO_d.
\end{align*}
\end{proof}

\begin{exmp}
    \citet{yang2021graph} reformulates graph convolution network (GCN) training~\citep{kipf2016semi} into a bilevel optimization problem, where the lower-level problem is to find the optimal node embeddings that minimize a graph energy function, and the upper-level problem is the loss of the embedding. For the node regression problem, we have the dataset $(\mA,\vb)$ with $n$ nodes. The bilevel problem is:
    \begin{align*}
        &\quad \min_{\vx \in \sR^d} \frac{1}{2} \Vert \vy^*(\vx) - \vb \Vert^2 \\
        &{\rm s.t.} \quad \vy^*(\vx) = \arg \min_{\vy \in \sR^{n}} \frac{1}{2} \Vert \mA \vx - \vy \Vert^2 + \frac{\lambda}{2} \vy^\top \mL \vy,
    \end{align*}
where $ \mL$ is the graph Laplacian matrix. In this example, the hyper-objective $\varphi(\vx) = f(\vx,\vy^*(\vx))$ is (strongly) convex if $\lambda>0$.
\end{exmp}

\begin{proof}
When $\lambda >0$, the lower-level solution has a closed form 
\begin{align*}
    \vy^*(\vx) = (\mI_n + \lambda \mL)^{-1} \mA \vx.
\end{align*}
Plugging the above equation into the upper-level function gives the closed form of $\varphi(\vx)$, which is a convex quadratic function because since its Hessian matrix is positive (semi)definite:
\begin{align*}
    \mH = \mA^\top (\mI_n + \lambda \mL)^{-1} (\mI_n + \lambda \mL)^{-1} \mA \succeq \mO_d.
\end{align*}
\end{proof}

\subsection{Nonconvex--Strongly-Convex Examples} \label{apx:exmple-NC--SC}

The NC--SC bilevel problems have many applications in the actor-critic method \citep{konda1999actor}, which is a class of important algorithms in reinforcement learning.  For this algorithm class, there is an actor in the upper level that optimizes a (typically nonconvex) cost function to select a good policy, and a critic in the lower level that evaluates by estimating the Q-function. The Q-function can be calculated based on the Bellman equation, which corresponds to solving a strongly convex function in the basic linear Markov decision process (MDP) setup.

Formally, we define the MDP as a tuple with five elements $(S,A, \gamma, P, R)$, where $S$ is the state space, $A$ is the action space, $\gamma \in [0,1)$ is the discount factor of the rewards, $P$ defines the transition kernel such that $P(s' \mid s,a)$ denotes the probability of transitioning to state $s'$ after taking action $a$
in state $s$, and $R$ defines the reward function such that $R(s,a) $ denotes the immediate reward received after taking action $a$ in state $s$. 
Let $Q^\pi$ denote the Q-function for policy $\pi$. 
Under the linear setup, $Q^\pi$ is assumed to have a linear representation $Q^\pi(s,a) = \langle \theta^\pi, \phi(s,a) \rangle$, where $ \theta^\pi \in \sR^d$ is the weight vector of parameters and $\phi(s,a) \in \sR^d$ is the 
feature vector for the state-action pair $(s,a)$. Then we can show that the induced bilevel optimization problem is nonconvex-strongly-convex with a quadratic lower-level problem.

\begin{exmp}[Actor-critic for linear MDP]
Let the initial state $s_0$ be drawn from a fixed distribution $\rho_0$, and let $\Pi \subset \sR^{\vert S \vert \times \vert A \vert}$ be the policy space. We consider the following problem:
\begin{align*}
    &\quad \max_{\pi \in \Pi} \sum_{a \in A} Q^{\pi}(s_0,a) \pi(a \mid s_0) \rho_0(s_0). \\
    &{\rm s.t.} \quad Q^\pi(s,a) = R(s,a) + \gamma \sum_{s' \in S}  \sum_{a \in A} Q^\pi(s',a) P(s' \mid s ,a) \pi(a \mid s).
\end{align*}
Using linear function approximation to $Q^\pi$, the lower-level function is a linear system in $\theta^\pi$, which can also be reformulated as an equivalent (strongly) convex quadratic optimization.
\end{exmp}

There are many similar examples in reinforcement learning, 
where the lower-level problem is exactly a linear system using linear function approximation.
Interested readers can refer to references \citep{hong2023two,zeng2024two,zeng2024fast}.

\section{Proof of Lemma \ref{lem:argmin-close-form}} \label{apx:proof-close-argmin}

\begin{proof}
Note that the function $\widehat Q_m(x_-,\,\cdot\,,x_+)$ is a quadratic function with a tridiagonal Hessian matrix, whose minimizer $\vx^*$ can be calculated by solving a simple recursion. Let
\begin{align} \label{eq:close-form-x-star}
    \vx^*_j = \theta x_- + \frac{m+j -1}{3m-1} \theta (x_+ - x_-).
\end{align}
We can verify that the vector $\vx = (x_1,\cdots,x_m)$ satisfies the first-order optimality condition:
\begin{align*}
    \begin{cases}
        \dfrac{1}{m} (x_1 - \theta x_-) = (x_2 - x_1), & i = 1;\\
        x_i - x_{i-1} = x_{i+1}  -x_i, & 1 < i <m; \\
        x_m - x_{m-1} = \dfrac{1}{m} (\theta x_+ - x_m), & i =m.
    \end{cases}
\end{align*}
Substituting the above values into the form of $\widehat Q_m(x_-,\,\cdot\,,x_+)$ leads to 
\begin{align*}
   \widehat Q_m(x_-, \vx^*,x_+) = (3m-1) (x_+-x_-)^2 = \frac{\theta^2 (x_2^* - x_1^*)^2}{3m-1}.
\end{align*}
\end{proof}

\section{Proof of Lemma \ref{lem:rela-stationarity}} \label{apx:proof-rela-sta}

\begin{proof}
Let $\mH_m \in \sR^{m \times m}$ be the Hessian matrix of the function $\widehat Q_m(x_-,\,\cdot\,,x_+)$ defined in Eq. (\ref{eq:our-coupling-chain}):
\begin{align*}
    \mH_m = 
    \begin{bmatrix}
        ~1+1/m & -1 \\
        -1 & 2 & -1 \\
        & -1 & \ddots & \ddots  \\
        & & \ddots & 2 & -1 \\
        & &  & -1 & 1+1/m \\
    \end{bmatrix}.
\end{align*}
Then the partial gradient with respect to $\vx$ is
\begin{align*}
    \nabla_\vx \widehat Q_m(x_-, \vx, x_+) = 2 \mH_m \vx.
\end{align*}
Recall that $\vx = (\bar \vx, \underline{\vx})$ and $F(\bar \vx) = \min_{\underline{\vx} \in (\sR^m)^T} \widehat f^\text{hs-nc}_{\nu,r}(\bar \vx, \underline{\vx})$, where the function $\widehat f^\text{hs-nc}_{\nu,r}(\vx)$ in Eq. (\ref{eq:our-coupling-chain}) takes the form of 
\begin{align*}
    \widehat f^\text{hs-nc}_{\nu,r}(\vx) = \frac{1}{2} \sum_{i=1}^T \widehat Q_m(\bar x_{i}, \underline{\vx}^{(i)} , \bar x_{i+1} ) + \Phi_{\nu,r}(\bar \vx),
\end{align*}
where $\Phi_{\nu,r}(\bar \vx)$ is a function that depends only on $\bar \vx$. Taking the partial gradient of this function yields
\begin{align*}
  \nabla_{\bar x_i} \widehat f^\text{hs-nc}_{\nu,r}(\vx) = \nabla_{\bar x_i} \Phi_{\nu,r}(\bar \vx)+
\frac{\theta}{m} \cdot
\begin{cases}
   2\theta \bar x_i - \underline{x}_1^{(i)} - \underline{x}_m^{(i-1)}, 
     & i <T; \\[1mm]
    \theta \bar x_i -  \underline{x}_1^{(i)}, & i =T.
    \\
\end{cases}
\end{align*}
In addition, using Danskin's theorem, we can 
know that $\nabla_{\bar x_i} F(\bar \vx) $ takes the same form except for $\underline{x}_1^{(i)}$ and $\underline{x}_m^{(i)}$ being replaced by $(\underline{x}_1^{(i)})^*$ and $(\underline{x}_m^{(i)})^*$, where 
where $(\underline{\vx}^{(i)})^*$ is the minimizer to the quadratic $\widehat Q_m(x_i, \,\cdot\, , x_{i+1})$ whose closed form is given by Eq.~(\ref{eq:close-form-x-star}).
In the vector form, we define $\vv^{(i)} = (\underline{\vx}^{(i)})^* - \underline{\vx}^{(i)}$ and then take the difference between the two partial gradients to obtain that
\begin{align*}
    \Vert \nabla_{\bar \vx} \widehat f^\text{hs-nc}_{\nu,r}(\bar \vx, \underline{\vx}) -\nabla F(\bar \vx) \Vert \le 
    \frac{\theta}{m} \sqrt{2\sum_{i=1}^T \left(v_1^{(i)}\right)^2 + \left(v_m^{(i)}\right)^2}.
\end{align*}
Let $\ve_j$ be the vector whose $j$-th element is 1 and the others are zero.
Using $\vert v_j \vert = \vert \langle \ve_j, \vv \rangle \vert \vert = \vert \langle \mH_m^{-1} \ve_j, \mH_m \vv \rangle \vert $ and Cauchy-Schwarz inequality, we can 
further obtain that
\begin{align*}
     \Vert \nabla_{\bar \vx} \widehat f^\text{hs-nc}_{\nu,r}(\bar \vx, \underline{\vx}) -\nabla F(\bar \vx) \Vert \le 
    \frac{2\theta}{m} \sqrt{ \sum_{i=1}^T \Vert \mH_m^{-1} \ve_1 \Vert^2 \Vert \mH_m \vv^{(i)} \Vert^2}.
\end{align*}
Note that each coordinate of the vector $\mH_m^{-1} \ve_1$ can be obtained by solving a tridiagonal linear system given by the matrix $\mH_m$, which has the closed form of 
\begin{align*}
    (\mH_m^{-1} \ve_1)_j = \frac{2m-j}{3-1/m}, \quad j = 1,\cdots,m.
\end{align*}
Hence, we can apply the sum of squares identity $\sum_{i=1}^n i^2 = n(n+1)(2n+1)/6$ and obtain that $\Vert \mH_m^{-1} \ve_1 \Vert \le m^{3/2}$, which then implies that
\begin{align*}
     \Vert \nabla_{\bar \vx} \widehat f^\text{hs-nc}_{\nu,r}(\bar \vx, \underline{\vx}) -\nabla F(\bar \vx) \Vert \le&
    2\theta \sqrt{m} \sqrt{ \sum_{i=1}^T \Vert \mH_m \vv^{(i)} \Vert^2} \\
    \le& \theta \sqrt{m} \sqrt{\sum_{i=1}^T \left\Vert \nabla_{\underline{\vx}^{(i)}} \widehat Q_m(\bar x_{i}, \underline{\vx}^{(i)} , \bar x_{i+1} ) \right\Vert^2} \\
    =& \theta\sqrt{m} \Vert \nabla_{\underline{\vx}} \widehat f^\text{hs-nc}_{\nu,r}(\bar \vx, \underline{\vx}) \Vert.
\end{align*}
Therefore,
\begin{align*}
    \Vert \nabla F(\bar \vx) \Vert \le& \Vert \nabla_{\bar \vx} \widehat f^\text{hs-nc}_{\nu,r}(\bar \vx, \underline{\vx}) \Vert +  \Vert \nabla_{\bar \vx} \widehat f^\text{hs-nc}_{\nu,r}(\bar \vx, \underline{\vx}) -\nabla F(\bar \vx) \Vert \\
    \le& \Vert \nabla_{\bar \vx} \widehat f^\text{hs-nc}_{\nu,r}(\bar \vx, \underline{\vx}) \Vert + \theta \sqrt{m}  \Vert \nabla_{\underline{\vx}} \widehat f^\text{hs-nc}_{\nu,r}(\bar \vx, \underline{\vx}) \Vert \\
    \le& (1 + \theta \sqrt{m}) \Vert \nabla \widehat f^\text{hs-nc}_{\nu,r}(\vx) \Vert,
\end{align*}
which means that $\Vert \nabla \widehat f^\text{hs-nc}_{\nu,r}(\vx) \Vert \ge \Vert \nabla F(\bar \vx) \Vert /2$ if $\theta \sqrt{m} \le 1$.
\end{proof}

\section{Proof of Theorem \ref{thm:lb-highligh-smooth-NC}} \label{apx:lb-highly-smooth-nc}

\begin{proof}
We split the proof of the lower bound into several pieces, following the analysis in \citet{carmon2020lower,carmon2021lower}. We will fix the following parameters
\begin{align} \label{eq:setting-m-theta}
       m =
    \left\lfloor
        \frac{1}{\sqrt{3 \nu}}
    \right\rfloor,
    \quad
    \theta
    =
    \sqrt{\nu (3m-1)}.
\end{align}
For sufficiently small $\nu$, we have $m \ge 2$. Moreover, 
\begin{align} \label{eq:m-theta-bound}
\theta^2 m
=
m(3m-1) \nu
\le
3m^2 \nu 
\le 1
\end{align}
Therefore \(\theta\sqrt m\le1\). Then the conditions of
Lemma~\ref{lem:argmin-close-form} and~\ref{lem:rela-stationarity} are both satisfied.
\paragraph{Zero-chain and large gradient.} 
Order the coordinates as
\[
\bar x_1,\ \underline{x}^{(1)}_1,\ldots,\underline{x}^{(1)}_m,\ \bar x_2,\ 
\underline{x}^{(2)}_1,\ldots,\underline{x}^{(2)}_m,\ \ldots,\ 
\bar x_T,\ \underline{x}^{(T)}_1,\ldots,\underline{x}^{(T)}_m,\ \bar x_{T+1}.
\]
With respect to this order, \(f\) is a first-order zero-chain. Therefore, before a zero-respecting method has revealed the tail of the chain, we have \( \bar x_T=\bar x_{T+1}=0\). 
From Lemma \ref{lem:argmin-close-form}, we know that $F(\bar \vx) = \nu \widehat f^{\rm nc}_{r}(\bar \vx)$, where $\widehat f^{\rm nc}_{r}$ is the nonconvex zero-chain in Definition \ref{dfn:nc-zo-chain}. Then we can further apply Lemma \ref{lem:Car-Upsion} to obtain that $\Vert \nabla F(\bar \vx) \Vert \ge \nu /4$ when an algorithm makes no more than $m (T-1)$ gradient queries, which in conjunction with Lemma \ref{lem:rela-stationarity} indicates that $\Vert \nabla \widehat f^\text{hs-nc}_{\nu,r}(\vx) \Vert \ge \Vert \nabla F(\bar \vx) \Vert/2 \ge \nu/8$.

\paragraph{Rescaling the function.} We then rescale the function by 
\[
f(\bar \vx, \underline{\vx})= \delta \beta^2 \widehat f^\text{hs-nc}_{\nu,r}(\vx / \beta ). 
\]
Under this scaling, we can set $\beta = 8 \eps/ (\delta \nu )$ to ensure that  $\Vert \nabla F(\vx) \Vert \ge \epsilon$.

\paragraph{Bound on suboptimality gap.} Combining the facts $ \min_{\underline{\vx} \in (\sR^m)^T} \widehat f^\text{hs-nc}_{\nu,r}(\bar \vx, \underline{\vx}) = \nu \widehat f^{\rm nc}_{r}(\bar \vx) $ by Lemma \ref{lem:argmin-close-form} and $\min_{\bar \vx \in \sR^{T+1}} \widehat f^{\text{nc}}_{r}(\bar \vx) = 0$ by Lemma \ref{lem:Car-Upsion}, for $d = m (T+1) + 1$, we have
\begin{align} \label{eq:gap-our-instance}
    F(\vzero) - \min_{\vx \in \sR^d} F(\vx) = \delta \nu \beta^2 \widehat f^\text{nc}_r(\vzero) = \delta \nu \beta^2 \left( \frac{1}{2} +  T \Upsilon_r(0) \right) \le \delta \nu \beta^2 \left( \frac{1}{2} + 10 T  \right),
\end{align}
where we use Lemma \ref{lem:Car-Upsion} that $\Upsilon_r(0) \le 10$ in the last step. Therefore, our hard instance can satisfy the constraint $F(\vzero) - \min_{\vx \in \sR^d} F(\vx) \le \Delta$ if 
\begin{align} \label{cond:delta-T}
 \delta \nu \beta^2 \le \Delta \quad {\rm and} \quad T = \left \lfloor \frac{\Delta}{20 \delta \nu \beta^2} \right \rfloor. 
\end{align}
\paragraph{Smoothness Constraints.} We bound the smoothness of the function by its quadratic part and non-quadratic part. Let \(f^{\rm quad}_\nu\) denote the quadratic
part of \(f\), namely
\begin{align} \label{eq:dfn-F-quad}
f^{\rm quad}_\nu(\vx)
=
\frac{\delta \nu}{2}(\beta-\bar x_1)^2
+
\frac{\delta}{2}\sum_{i=1}^T \widehat Q_m(\bar x_i,x^{(i)},\bar x_{i+1}).    
\end{align}
Under the setting $\theta \sqrt{m} \le 1$, it is not  hard to verify that 
\(\Vert \nabla^2f^{\rm quad}_\nu(\vx) \Vert \le 4 \delta \). 
Consequently, choosing $\delta = L_1/8$ ensures that the quadratic part contributes at most \(L_1/2\) to the gradient
Lipschitz constant. The only derivatives of order at least three come from the nonconvex regularizer $\Upsilon_r$. By
Lemma \ref{lem:Car-Upsion}, for every \(q\ge1\), the Lipschitz constant of the $q$-th order derivatives of the component function \(\Upsilon_r\!\left({\,\cdot\,}\right)\) is bounded by \( \ell_q\), where $\ell_q$ is a numerical constant that only depends on $q$. Therefore, the Lipschitz constant of the $q$-th order derivatives of the component function \(\delta \nu \beta^2 \Upsilon_r\!\left({\,\cdot\,}/ \beta\right)\) is bounded by
\begin{align} \label{eq:contri-NC-regularizer}
    \delta \nu \beta^{1-q}r^{3-q}\ell_q
    =
    \frac{(\delta \nu)^qr^{3-q} \ell_q }{(8\epsilon)^{q-1}},
\end{align}
where the equation is due to the parameter setting \(\beta=8 \eps/(\delta \nu)\). 
We now choose the other parameters case by case, depending on the different orders of smoothness, and then conclude the lower bounds for each case.
\begin{enumerate}
    \item \textbf{Case \(p=1\).} Set $r=1$, $\delta =L_1/8$, and $\beta = 8 \eps/ (\delta \nu)$.
Then, the nonconvex regularizer's contribution to the gradient Lipschitz constant in Eq.~(\ref{eq:contri-NC-regularizer}) is
\[
    \delta \nu \ell_1 
    \le 
    \frac{L_1}{2}, \quad \text{if} \quad \nu = \Theta(1).
\]
Combined with \(\Vert \nabla^2f^{\rm quad}_{\nu}(\vx) \Vert \le L_1/2\) for $\delta = L_1/8$, we can conclude that the total gradient Lipschitz constant is at most \(L_1\). For
sufficiently small~\(\epsilon\), we can select $T$ according to Eq. (\ref{cond:delta-T}), which leads to the lower bound of 
\begin{align*}
    m(T-1) = \Omega \left( {\Delta L_1 \sqrt{\nu}}{\eps^{-2}}  \right) = \Omega \left( \Delta L_1 \epsilon^{-2} \right).
\end{align*}
\item \textbf{Case $p=2$.}  Set $r=1$, $\delta =L_1/8$, and $\beta = 8 \eps/ (\delta \nu)$.
According to Eq. (\ref{eq:contri-NC-regularizer}), the Hessian Lipschitz constant contributed by the nonconvex regularizer is
\[
   \frac{(\delta \nu)^2 \ell_2}{8\eps} 
    \le L_2, \quad \text{if} \quad \nu = \Theta \left( \frac{\sqrt{L_2 \eps}}{L_1}
    \right).
\]
The remaining first-order smoothness contribution from the regularizer is
\[
    \delta\nu \ell_1
    =
    \Theta\left( \sqrt{L_2 \eps} \right),
\]
which is at most \(L_1/2\) for sufficiently small \(\epsilon\). Therefore, $f$ belongs to the class of second-order smooth problems with parameters $L_1$ and $L_2$. Recalling the setting of $T$ in Eq. (\ref{cond:delta-T}), we have
\begin{align} \label{eq:lb-Tm}
    m(T-1)
    = \Omega \left( {\Delta L_1 \sqrt{\nu}}{\eps^{-2}}  \right)
     = \Omega\!\left(
        \Delta L_1^{1/2}L_2^{1/4}\epsilon^{-7/4}
    \right). 
\end{align}
\item \textbf{Case $p=3$.} 
Set $r=1$, $\delta =L_1/8$, and $\beta = 8 \eps/ (\delta \nu)$. According to (\ref{eq:contri-NC-regularizer}), the third-derivative Lipschitz constant contributed by the nonconvex regularizer is
\[
    \frac{(\delta \nu)^3 \ell_3}{(8\epsilon)^2}
    \le L_3, \quad \text{if} \quad \nu = \Theta \left(
    \frac{L_3^{1/3} \eps^{2/3}}{L_1}
    \right).
\]
The lower-order contributions of the regularizer are
\[
    \delta \nu \ell_1 
    =
    \Theta(L_3^{1/3}\epsilon^{2/3}), \quad {\rm and} \quad \frac{(\delta \nu)^2 \ell_2}{8\eps}
    =
    \Theta(L_3^{2/3}\epsilon^{1/3}).
\]
Thus, for sufficiently small \(\epsilon\), these contributions are at most
\(L_1/2\) and \(L_2\), respectively.
Therefore, $f$ belongs to the class of third-order smooth problems with parameters $L_1$, $L_2$, and $L_3$. The final lower bound  is
\[
m (T-1) = \Omega(\Delta L_1 \sqrt{\nu} \eps^{-2}) = \Omega\!\left(
    \Delta L_1^{1/2}L_3^{1/6}\epsilon^{-5/3}
    \right).
\]
\item  \textbf{Case $p \ge 4$.} 
Set $r = (8 \eps)^{-2/3}$, $\delta = L_1/ 8$, and $\beta = 8 \eps/ (\delta \nu)$. According to (\ref{eq:contri-NC-regularizer}), the third-derivative Lipschitz constant contributed by the nonconvex regularizer is
\[
\left( \frac{\delta \nu}{8 \eps r} \right)^q
 8 \epsilon r^3 
\ell_q = \left( \frac{\delta \nu}{(8 \eps)^{2/3}} \right)^q \ell_q \le L_q /2, \quad \text{if} \quad \nu = \Theta_p \left( \frac{\eps^{2/3} \Lambda_p}{L_1}  \right),
\]
where \(\Lambda_p = \min_{1 \le q \le p} (L_q / \ell_q)^{1/q} \). Therefore, $f$ belongs to the class of $p$th-order smooth problems with parameters $L_1,\cdots,L_p$. The final lower bound  is
\[
m(T-1) = \Omega \left( {\Delta L_1 \sqrt{\nu}}{\eps^{-2}}  \right) = \Omega_p\!\left(
        \Delta L_1^{1/2}
        \Lambda_p^{1/2}
        \epsilon^{-5/3} \right).
\]
\end{enumerate}
\end{proof}

\section{Proof of Lemma \ref{lem:matrix-A}} \label{apx:proof-property-Bk}
\begin{proof}
Note that the matrix $\mA_n \in \sR^{n\times n}$, regardless of the endpoints, is equivalent to the 
canonical tridiagonal graph Laplacian matrix for the simple path. Denote such a matrix as $\mB_{n-2} \in \sR^{(n-2) \times (n-2)}$, which is
\begin{align} \label{eq:matrix-Bk}
    \mB_{n-2}
    =
    \begin{bmatrix}
        2 & -1 \\
        -1 & 2 & -1 \\
        & -1 & \ddots & \ddots  \\
        & & \ddots & 2 & -1 \\
        & &  & -1 & 2\\
    \end{bmatrix}.
\end{align}
Then, we can equivalently write the matrix $\mA_n$ as
\begin{align} \label{eq:matrix-Ak-dense-form}
    \mA_n = 
    \begin{bmatrix}
        \omega & -q \ve_1^\top & 0 \\
        -q \ve_1 & \mB_{n-2} & -q \ve_{n-2} \\
        0 & -q \ve_{n-2}^\top & \omega
    \end{bmatrix}.
\end{align}
A trick we use to simplify the subsequent analysis is to reorder the coordinates from 
$\vy = (y_1,y_2,\ldots,y_{n})$ to $(\vu,\vv)$, where the new order is defined by
\[
    \vu=(y_1,y_{n})\in\mathbb{R}^2
    \quad {\rm and} \quad
    \vv=(y_2,\ldots,y_{n-1})\in\mathbb{R}^{n-2},
\]
where $\ve_i \in \sR^{n-2}$ denotes the $i^{\text{th}}$ standard basis vector, whose sole nonzero entry equals $1$.
Let
\(\mR=
    \begin{bmatrix}
    \ve_1,
    \ve_{n-2}
    \end{bmatrix}
    \in\mathbb{R}^{(n-2) \times 2}\).
After this permutation from $\vy \rightarrow (\vu,\vv)$, the matrix \(\mA_n\) becomes
\begin{align} \label{eq:matrix-hat-Ak}
    \widehat \mA_n
    =
    \begin{bmatrix}
    \omega \mI_2 & -q \mR^\top\\
    -q\mR & \mB_{n-2}
    \end{bmatrix}. 
\end{align}
Let \(\mP\) be the corresponding permutation matrix, so that
\[
    \widehat \mA_n=\mP \mA_n \mP^\top.
\]
Then, the leading \(2\times 2\) block of \(\widehat \mA_n^{-1}\) is exactly
the endpoint block of \(\mA_n^{-1}\), namely
\begin{align} \label{eq:hat-Bk-matrix}
    \left(\widehat \mA_n^{-1}\right)_{uu}  =
    \left(  \mP \mA_n^{-1} \mP^\top \right)_{uu}
    =
    \begin{bmatrix}
    (\mA_n^{-1})_{1,1} & (\mA_n^{-1})_{1,n}\\
    (\mA_n^{-1})_{n,1} & (\mA_n^{-1})_{n,n}
    \end{bmatrix}. 
\end{align}
Therefore, we can analyze the values of  $(\mA_n^{-1})_{1,1}$ and $(\mA_n^{-1})_{1,n}$ by calculating the submatrix $(\widehat \mA_n^{-1})_{uu}$. Moreover, we note that the permutation does not change the eigenvalues of the matrix.

\paragraph{Proof of Part (a).} Since the submatrix \(\mB_{n-2}\) is invertible, the Schur complement of \(\mB_{n-2}\) in
\(\widehat \mA_n\) is
\begin{align} \label{eq:matrix-S}
    \mS
    =
    \omega \mI_2-q^2 \mR^\top \mB_{n-2}^{-1} \mR. 
\end{align}
In light of Eq. (\ref{eq:hat-Bk-matrix}), we calculate the inverse of $\mS$, which by the block inverse formula is
\begin{align} \label{eq:A-S-inverse}
\begin{bmatrix}
    (\mA_n^{-1})_{1,1} & (\mA_n^{-1})_{1,n}\\
    (\mA_n^{-1})_{n,1} & (\mA_n^{-1})_{n,n}
    \end{bmatrix}=
    \left(\widehat \mA_n^{-1}\right)_{uu}
    =
    \mS^{-1}.    
\end{align}
Now, let us compute \(\mS\) in Eq. (\ref{eq:matrix-S}), where the nontrivial part is $\mR^\top \mB_{n-2}^{-1} \mR^\top$. Since $ \mR=
    \begin{bmatrix}
    \ve_1,
    \ve_{n-2}
    \end{bmatrix}$, we have
\begin{align} \label{eq:RBR-B}
\mR^\top \mB_{n-2}^{-1} \mR^\top = 
\begin{bmatrix}
(\mB_{n-2})^{-1}_{1,1} & (\mB_{n-2})^{-1}_{1,n-2} \\
(\mB_{n-2})^{-1}_{n-2,1}   & (\mB_{n-2}^{-1})_{n-2,n-2}
\end{bmatrix}.    
\end{align}
We now explicitly calculate each element of the above $2 \times 2$ matrix. 
For a fixed $j \in [n-2]$, we define a vector~$\vg^{(j)}$ elementwise as
\begin{align} \label{eq:dfn-gj}
\begin{split}
       g^{(j)}_i = \frac{i \land j (n-1  - i \lor j)}{n-1} =
    \begin{cases}
    \dfrac{i (n-1 - j)}{n-1}, & i\le j ;\\
    \dfrac{j (n-1-i)}{n-1}, & i >j.
\end{cases}. 
\end{split}
\end{align}
It is easy to see that $g_i^{(j)}$ are both linear in $i$ when $i < j$ or $i >j$. Therefore, we have
\begin{align} \label{eq:identity-ij}
    2 g_i^{(j)} - g_{i-1}^{(j)}  - g_{i+1}^{(j)} = 0, \quad \forall i< j ~~{\rm or}~~ i>j.
\end{align}
Moreover, by explicitly substituting the values of $g_{j-1}^{(j)}$, $g_j^{(j)}$, it is easy to verify that
\[
2 g_i^{(j)} - g_{i-1}^{(j)}  - g_{i+1}^{(j)} = 1, \quad i = j.
\]
Writing Eq. (\ref{eq:identity-ij}) in matrix form, we have
\begin{align} \label{eq:Age}
    \mB_{n-2} \vg^{(j)} = \ve_j.
\end{align}
Therefore, $\vg^{(j)} = \mB_{n-2}^{-1} \ve_j$ is exactly the $j$th column of the inverse matrix $ \mB_{n-2}^{-1}$, which means that 
\begin{align*}
    &(\mB_{n-2}^{-1})_{1,1} = (\mB_{n-2}^{-1})_{n-2,n-2} = g_1^{(1)} = \frac{n-2}{n-1} ; \\
    &(\mB_{n-2}^{-1})_{1,n-2}  = (\mB_{n-2}^{-1})_{n-2,1} = g_1^{(n-2)} = \frac{1}{n-1}.
\end{align*}
Substituting into Eq. (\ref{eq:RBR-B}), we can obtain
\[
    \mR^\top \mB_{n-2}^{-1} \mR
    =
    \begin{bmatrix}
    \dfrac{n-2}{n-1} & \dfrac{1}{n-1}\\[6pt]
    \dfrac{1}{n-1} & \dfrac{n-2}{n-1}
    \end{bmatrix}.
\]
For $q = 1/\sqrt{n-1}$ and $\omega = n / (n-1)^2$, we have
\begin{align} \label{eq:form-schur-comple}
    \mS
    &=
    \frac{n}{(n-1)^2} \mI_2
    -
    \frac{1}{n-1}
    \begin{bmatrix}
    \dfrac{n-2}{n-1} & \dfrac{1}{n-1}\\[6pt]
    \dfrac{1}{n-1} & \dfrac{n-2}{n-1}
    \end{bmatrix}  =
    \frac{1}{(n-1)^2}
    \begin{bmatrix}
    2&-1\\
    -1&2
    \end{bmatrix}.
\end{align}
Therefore, its inverse has the simple form of
\[
    \mS^{-1}
    =
    \frac{(n-1)^2}{3}
    \begin{bmatrix}
    2&1\\
    1&2
    \end{bmatrix}.
\]
Hence, combining with Eq. (\ref{eq:A-S-inverse}), we can conclude part (a) that
\[
    (\mA_n^{-1})_{1,1}
    =
    (\mA_n^{-1})_{n,n}
    =
    \left(\mS^{-1} \right)_{1,1} = 
    \frac{2(n-1)^2}{3},
\]
and
\[
    (\mA_n^{-1})_{1,n}
    =
    (\mA_n^{-1})_{n,1}
    =
    \left(\mS^{-1} \right)_{1,2} =
    \frac{(n-1)^2}{3}.
\] 
\paragraph{Proof of Part (b).} 
In the following, we give the spectral upper and lower bounds of the matrix $\mA_n$. At a high level, since $\mA_n$ is only an extension of $\mB_{n-2}$ by 2 dimensions, the spectral bounds of $\mA_n$ is also close to those of $\mB_{n-2}$. Note that the spectrum of the tridiagonal matrix $\mB_{n-2}$ can be calculated as
\[
\lambda_k(\mB_{n-2}) = 2 - 2 \cos \frac{k \pi}{n-1} = 4 \sin^2 \frac{k \pi}{2 (n-1)}, \quad k = 1,\cdots,n-2.
\]
Therefore, using the elementary fact $ 2 x / \pi \le \sin x$ for $ 0 \le x \le \pi/2$, we have 
\begin{align} \label{eq:Bk-succeq}
   \frac{4}{(n-1)^2} \mI_{n-2} \preceq \mB_{n-2} \preceq 4 \mI_{n-2}.
\end{align}

Since the permutation does not change the eigenvalues of the matrix, we can give the upper bound of the maximal eigenvalue of $\mA_n$ by that of $ \widehat \mA_n $ equivalently. Note that 
the quadratic form of the matrix $\widehat \mA_n$ under the permutation from $\vy \rightarrow (\vu,\vv)$ satisfies that
\begin{align} \label{eq:quadratic-Bk-uv}
    \begin{bmatrix}\vu \\ \vv\end{bmatrix}^\top
    \widehat \mA_n
    \begin{bmatrix}\vu\\ \vv\end{bmatrix}
    =&
     \omega \|\vu\|^2-2q\,\vu^\top \mR \vv+\vv^\top \mB_{n-2} \vv.
\end{align}
Hence, it suffices to upper-bound the maximal eigenvalues of the matrices $\mR$ and $\mB_{n-2}$. Since the first matrix is simply $ \mR=
    \begin{bmatrix}
    \ve_1,
    \ve_{n-2}
    \end{bmatrix}$, we can immediately know that $\Vert \mR \Vert \le 1$. For the second matrix~$\mB_{n-2}$, we know from inequality (\ref{eq:Bk-succeq}) that $\Vert \mB_{n-2} \Vert \le 4$. Therefore, for  $q = 1/\sqrt{n-1} \le 1/2$ and $\omega = n / (n-1)^2 \le 1$ for $n \ge 5$, Eq.~(\ref{eq:quadratic-Bk-uv}) can be upper-bounded by
\[
\begin{aligned}
  \begin{bmatrix}\vu \\ \vv\end{bmatrix}^\top
    \widehat \mA_n
    \begin{bmatrix}\vu\\ \vv\end{bmatrix}
    \le&
    \omega \|\vu\|^2-2q \Vert \mR \Vert \Vert \vu \Vert \Vert \vv \Vert +\Vert \mB_{n-2} \Vert \Vert \vv \Vert^2  \\
    \le & \Vert \vu \Vert^2 +   \Vert \vu \Vert \Vert \vv \Vert + 4 \Vert \vv \Vert^2  \\
    \le& 1.5 \Vert \vu \Vert^2 + 4.5 \Vert \vv \Vert^2 \\
    \le& 5 (\Vert \vu \Vert^2 + \Vert \vv \Vert^2).
\end{aligned}
\]
Therefore, the maximal eigenvalue of $\mA_n$ can be bounded by
\begin{align} \label{eq:upper-eigen-bound}
 \lambda_{\max}\left(\mA_n\right) =
    \lambda_{\max} \left(\widehat \mA_n \right) \le 5.   
\end{align}
As for the minimal eigenvalue, we relate $\lambda_{\min}(\mA_n)$ to $\lambda_{\min}(\mB_{n-2})$ using the following quadratic identity derived from Eq.~(\ref{eq:quadratic-Bk-uv}):
\begin{align} \label{eq:quadratic-W}
\begin{split}
\begin{bmatrix}\vu \\ \vv\end{bmatrix}^\top
    \widehat \mA_n
    \begin{bmatrix}\vu\\ \vv\end{bmatrix}
    =&
     \omega \|\vu\|^2-2q\,\vv^\top \mR \vu+\vv^\top \mB_{n-2} \vv \\
     =& (\vv - q \mB_{n-2}^{-1} \mR \vu)^\top \mB_{n-2} ( \vv - q \mB_{n-2}^{-1} \mR \vu) + \vu^\top \mS \vu.
\end{split}
\end{align}
Via a technical calculation (which is deferred to the end), we can show that the matrix $\mW= q \mB_{n-2}^{-1} \mR \in \sR^{(n-2) \times 2}$ satisfies $\Vert \mW \Vert^2 \le 1/2$. Recalling the explicit form of \(\mS\) in Eq. (\ref{eq:form-schur-comple}) and the matrix inequality for \(\mB_{n-2}\) in Eq.~(\ref{eq:Bk-succeq}), we have
\[
\mS \succeq \frac{1}{(n-1)^2} \mI_2 \quad {\rm and} \quad \mB_{n-2} \succeq \frac{4}{(n-1)^2} \mI_{n-2}. 
\]
Therefore, the quadratic form in Eq. (\ref{eq:quadratic-W}) can be lower-bounded by
\begin{align*}
    \begin{bmatrix}\vu \\ \vv\end{bmatrix}^\top
    \widehat \mA_n
    \begin{bmatrix}\vu\\ \vv\end{bmatrix}
    \ge& \frac{1}{(n-1)^2} \left( \Vert \vv - \mW \vu \Vert^2 + \Vert \vu \Vert^2 \right).
\end{align*}
To give a lower bound of the minimal eigenvalue of $\widehat \mA_n$, it remains to compare the right-hand side with $(\Vert \vu \Vert^2 + \Vert \vv \Vert^2)$ using the fact that $\Vert \mW \Vert^2 \le 1/2$. Continuing from the above inequality, we obtain
\begin{align*}
    \begin{bmatrix}\vu \\ \vv\end{bmatrix}^\top
    \widehat \mA_n
    \begin{bmatrix}\vu\\ \vv\end{bmatrix}
    \ge& \frac{1}{2(n-1)^2} \left( 2\Vert \vv - \mW \vu \Vert^2 + 2\Vert \vu \Vert^2 \right) \\
    \ge& \frac{1}{2(n-1)^2} \left( 2\Vert \vv - \mW \vu \Vert^2 + 2 \Vert \mW \Vert^2 \Vert \vu  \Vert^2 + \Vert \vu \Vert^2 \right) \\
    \ge& \frac{1}{2(n-1)^2} \left( 2\Vert \vv - \mW \vu \Vert^2 + 2 \Vert \mW  \vu   \Vert^2 + \Vert \vu \Vert^2 \right) \\
    \ge& \frac{1}{2(n-1)^2} \left( \Vert \vv \Vert^2 + \Vert \vu \Vert^2 \right),
\end{align*}
which indicates that $\lambda_{\min}(\widehat \mA_n) \ge 1/ (2(n-1)^2)$. Then, the proof of part (b) is completed by combining the bounds for $\lambda_{\min}(\widehat \mA_n)$ in the above and $\lambda_{\max}(\widehat \mA_n)$ in Eq. (\ref{eq:upper-eigen-bound}). Finally, it remains to show that $\Vert \mW \Vert^2 \le 1/2$ in the following.

\paragraph{Upper Bound of $\Vert \mW \Vert^2$.} Recall that $ \mR=
    \begin{bmatrix}
    \ve_1,
    \ve_{n-2}
    \end{bmatrix}$ and $\mW= q \mB_{n-2}^{-1} \mR = 
    \begin{bmatrix}
        \vw_1, \vw_2
    \end{bmatrix}
    $. Using the identity (\ref{eq:Age}) with $\vg^{(j)} = \mB_{n-2}^{-1} \ve_j$ for the vector \(\vg^{(j)}\) defined in Eq. (\ref{eq:dfn-gj}), we know that the vectors $\vw_1$ and $\vw_2$ are 
    \[
    \begin{aligned}
        (\vw_1)_i = & q (\mB_{n-2}^{-1} \ve_1)_i = q \cdot g_i^{(1)} = q \cdot \frac{n-1-i}{n-1}; \\
        (\vw_2)_i = & q (\mB_{n-2}^{-1} \ve_{n-2})_i = q \cdot g_i^{(n-2)} = q\cdot \frac{i}{n-1}.
    \end{aligned}
    \]
Therefore, for $q = 1/\sqrt{n-1}$, we can calculate that
\[
\begin{aligned}
    &\Vert \vw_1 \Vert^2 = \Vert \vw_2 \Vert^2 = q^2 \sum_{i=1}^{n-2} \left( \frac{i}{n-1} \right)^2 = \frac{(n-2)(2n-1)}{6(n-1)^2}; \\
    &\langle \vw_1, \vw_2 \rangle = q^2 \sum_{i=1}^{n-2} \frac{(n-1-i)i}{(n-1)^2} = \frac{(n-2)n}{6(n-1)^2}.
\end{aligned}
\]
It means that
\[
\mW^\top \mW = 
\begin{bmatrix}
    \Vert \vw_1 \Vert^2 & \langle \vw_1,\vw_2 \rangle \\
    \langle \vw_1,\vw_2 \rangle  & \Vert \vw_2 \Vert^2
\end{bmatrix}
= \frac{n-2}{6(n-1)^2} 
\begin{bmatrix}
    2n-1 & n \\
    n & 2n-1
\end{bmatrix},
\]
whose maximal eigenvalue can be bounded by
\[
\Vert \mW \Vert^2 = \frac{n-2}{6(n-1)^2} \times 3 (n-1) = \frac{n-2}{2(n-1)} \le \frac{1}{2}.
\]
\end{proof}

\section{Improved Upper Bounds} \label{apx:ub-refined}

Most existing studies on the upper bounds of bilevel optimization primarily focus on optimizing the $\epsilon$ dependency. In this section, we first review the approaches for establishing these upper bounds, showing their two-source complexity that corresponds to the core idea of our lower bound. Then, we also provide refinements of them by improving the lower-level gradient descent sub-solver with Nesterov acceleration \citep{nesterov1983method}.

\subsection{Overview of Upper Bounds}

The main challenge of designing a fully first-order method for bilevel optimization is that the algorithm only has access to first-order oracles of $f$ and $g$, while the calculation of hyper-gradient $\nabla F$ in Eq. (\ref{eq:hyper-grad}) requires second-order information of $g$. This issue was addressed by \citet{kwon2023fully}, who proposed the following fully first-order estimator:
\begin{align} \label{eq:f2ba-grad}
    \nabla F_{\lambda} (\vx) = \nabla_\vx f(\vx, \vy_{\lambda}^*(\vx)) + \lambda ( \nabla_\vx g(\vx, \vy_{\lambda}^*(\vx)) - \nabla_\vx g(\vx,\vy^*(\vx)),
\end{align}
where $ \vy_{\lambda}^*(\vx) =  \arg \min_{\vy \in \sR^{d_y}} \{f(\vx,\vy) + \lambda g(\vx,\vy)\}$. \citet{kwon2023fully} showed that $\nabla F_{\lambda}$ is a good approximation to $\nabla F$ for a large penalty $\lambda$, as we recall below.
\begin{lem}[{\citet[Lemma 3.1]{kwon2023fully}}] Let $\lambda \ge 2 L_1/ \mu_\vy$. For first-order smooth NC--SC bilevel problem $(f,g) \in$ $\gF^{\text{nc--sc}}(L_0,L_1,L_2,\mu_\vy,\Delta)$, 
$\nabla F_{\lambda}$ is well-defined and satisfies
\begin{align*}
    \Vert \nabla F_{\lambda}(\vx) - \nabla F(\vx) \Vert \le 8 \bar \kappa_\vy^3 L_1 / \lambda.
\end{align*}
\end{lem}
Based on this key lemma, \citet{chen2023near} proposed the Fully First-order Bilevel Approximation (F${}^2$BA) method that achieved a complexity of
$\tilde \gO(\bar \kappa_\vy^4 \epsilon^{-2})$. The F${}^2$BA method consists of two loops:
(1) The outer loop performs gradient descent (GD) on $\vx$ using the hyper-gradient estimator~$\nabla F_\lambda$, under the setting $\lambda = \Omega(\epsilon^{-1})$ that ensures  $\Vert \nabla F_\lambda(\vx) - \nabla F(\vx) \Vert \le \epsilon$. (2) The inner loop computes the lower-level solutions $\vy_{\lambda}^*(\vx)$ and $\vy^*(\vx)$ by GD and uses them to construct the estimator $\nabla F_\lambda(\vx)$. 

Actually, better upper bounds on condition number dependency can be obtained by replacing the lower-level GD with its acceleration AGD \citep{nesterov1983method}, due to the following lemma for Algorithm~\ref{alg:AGD}.



\begin{algorithm*}[htbp]  
\caption{AGD $(h,K,  \vz_0)$} \label{alg:AGD}
\begin{algorithmic}[1] 
\STATE $ \tilde \vz_0 = \vz_0$ \\[1mm]
\STATE \textbf{for} $ k =0,1,\cdots,K-1 $ \\[1mm]
\STATE \quad $\vz_{k+1} = \tilde \vz_k - \frac{1}{L_h} \nabla h(\tilde \vz_k)$ \\[1mm]
\STATE \quad $\tilde \vz_{k+1} = \vz_{k+1} + \frac{\sqrt{\kappa_h}-1}{\sqrt{\kappa_h} +1} (\vz_{k+1} - \vz_k) $ \\[1mm]
\STATE \textbf{end for} \\[1mm]
\STATE \textbf{return} $\vz_K$ 
\end{algorithmic}
\end{algorithm*}

\begin{lem}[{\citet[Theorem 2.2.4]{nesterov2018lectures}}] \label{lem:AGD}
Running Algorithm \ref{alg:AGD} on an $L_h$-gradient Lipschitz and $\mu_h$-strongly convex objective $h: \sR^d \rightarrow \sR$ outputs $\vz_K$ satisfying
\begin{align*}
    \Vert \vz_K - \vz^* \Vert^2 \le (1 + \kappa_h) \left( 1 - \frac{1}{\sqrt{\kappa_h}} \right)^K \Vert \vz_0 - \vz^* \Vert^2,
\end{align*}
where $\vz^* = \arg \min_{\vz \in \sR^d} h(\vz)$ and $\kappa_h = L_h/ \mu_h$ is the condition number of $h$.
\end{lem}

\subsection{Improved Upper Bounds for First-order Smooth Problems} \label{subsec:F2BA}


We present the F${}^2$BA${}^+$ method in Algorithm \ref{alg:F2BA} for NC--SC bilevel problems, which uses AGD to solve the lower-level problem in the F${}^2$BA method \citep{chen2023near}. The superscript~``$^+$'' is used to distinguish from the original algorithm in \citet{chen2023near}.
The analysis mainly follows that of \citep{chen2023near} and relies on the following two lemmas. The first lemma gives the smoothness of $F(\vx)$.

\begin{lem} \label{lem:smooth-phi}
For a smooth NC--SC bilevel problem $(f,g) \in$ $\gF^{\text{nc--sc}}(L_0,L_1,L_2,\mu_\vy,\Delta)$, the hyper-objective $F(\vx) = f(\vx,\vy^*(\vx))$ satisfies
\begin{align} \label{eq:hyper-smooth}
    \Vert \nabla F(\vx) - \nabla F(\vx') \Vert \le L_F \Vert \vx-  \vx' \Vert, \quad \forall \vx,\vx' \in \sR^{d_x},
\end{align}
where $L_F = \gO( \bar \kappa_\vy^3 L_1 )$ and $\bar \kappa_\vy = \max\{L_0,L_1,L_2 \} / \mu_\vy$. In particular, if the lower-level function is quadratic and $(f,g) \in \gF^{\text{nc--q}}(L_1,\mu_\vy, \Delta)$, we have $L_F = \gO(\kappa_\vy^2 L_1)$, where $\kappa_\vy = L_1/ \mu_\vy$.
\end{lem}

\begin{proof}
In the general case, $L_F = \gO( \bar \kappa_\vy^3 L_1)$ can be verified by Eq. (\ref{eq:hyper-grad}) or see {\citep[Lemma 2.2]{ghadimi2018approximation}}. For quadratic lower-level problems defined in Eq.~(\ref{eq:quadratic-g}), we have $\nabla F(\vx) = \nabla_\vx f(\vx,\vy^*(\vx)) - \mJ \mH^{-1} \nabla_\vy f(\vx,\vy^*(\vx))$. Then, from the facts $\Vert \mJ \Vert \le L_1$, $\Vert \mH^{-1} \Vert \le 1/\mu_\vy$, and $\Vert \nabla \vy^*(\vx) \Vert \le \kappa_\vy$, we can easily prove that $L_F = \gO(\kappa_\vy^2 L_1)$.
\end{proof}

The second lemma from \citep{chen2023near}, states that if the lower-level solver achieves sufficient accuracy, then the outer loop of Algorithm \ref{alg:F2BA} requires $\gO(L_F \Delta \epsilon^{-2})$ iterations to find an $\epsilon$-stationary point.

\begin{algorithm*}[t]  
\caption{F${}^2$BA${}^+$ $(f,g, \eta_x, T, K,\vx_0,\vy_0)$} \label{alg:F2BA}
\begin{algorithmic}[1] 
\STATE $ \vz_0 = \vy_0$ \\[1mm]
\STATE \textbf{for} $ t =0,1,\cdots,T-1 $ \\[1mm]
\STATE \quad $\vz_{t+1} = \text{AGD}( g(\vx_t,\,\cdot\,), K, \vz_t)$ \\[1mm]
\STATE \quad $\vy_{t+1} = \text{AGD} (f(\vx_t,\,\cdot\,) + \lambda g(\vx_t, \,\cdot\,) , K, \vy_t    )$ \\[1mm]
\STATE \quad $ \hat \nabla F_\lambda(\vx_t)= \nabla_\vx f(\vx_t,\vy_{t+1}) + \lambda ( \nabla_\vx g(\vx_t,\vy_{t+1}) - \nabla_\vx g(\vx_t,\vz_{t+1}) )$ \\[1mm]
\STATE \quad  $\vx_{t+1} = \vx_t -  \eta_x  \hat \nabla F_{\lambda}(\vx_t)$ \\[1mm]
\STATE \textbf{end for} \\[1mm]
\STATE \textbf{return} $\vx_{\rm out} = \frac{1}{T} \sum_{t=0}^{T-1} \vx_{t}$ \\[1mm]
\end{algorithmic}
\end{algorithm*}

\begin{lem}[{\citet[Theorem 4.1]{chen2023near}}]\label{lem:F2BA}
For a first-order smooth NC--SC bilevel problem $(f,g) \in$ $\gF^{\text{nc--sc}}(L_0,L_1,L_2,\mu_\vy,\Delta)$, 
if in each iteration $t$ of Algorithm \ref{alg:F2BA} the lower-level solver outputs 
$\vz_t^K$ and $\vy_t^K$ such that 
\begin{align*}
    \Vert \vz_{t+1} - \vy^*(\vx_t) \Vert^2+ \Vert \vy_{t+1} - \vy_{\lambda}^*(\vx_t) \Vert^2 \le \gamma \Vert \vz_t - \vy^*(\vx_t) \Vert^2+ \Vert \vy_t - \vy_{\lambda}^*(\vx_t) \Vert^2
\end{align*}
with $\gamma \lesssim \min\{ 1/2, \epsilon^2/ (\bar \kappa_\vy^3 L_1^3) \} / \lambda^2$, then setting $\lambda \asymp L_1 \bar \kappa_\vy^3 \epsilon^{-1}$ and $\eta_x \asymp 1/L_F$ can ensure the algorithm to find an $\epsilon$-stationary point of $F(\vx)$ in $T = \gO( L_{F} \Delta \epsilon^{-2})$ outer-loop iterations.
\end{lem}



Combining Lemma \ref{lem:AGD}, \ref{lem:smooth-phi}, and \ref{lem:F2BA}, we can easily prove that the total complexity of F${}^2$BA${}^+$ is $\tilde \gO(\sqrt{\kappa_\vy} L_F \Delta \epsilon^{-2})$,
where $L_F = \gO(\bar \kappa_\vy^3 L_1)$ for general lower-level functions and $L_F = \gO(\kappa_\vy^2 L_1) $ for quadratic lower-level functions. We summarize the complexity of F${}^2$BA${}^+$ in the following, which improves the $\tilde \gO(\kappa_\vy L_F \Delta \epsilon^{-2})$
complexity of F${}^2$BA \citep{chen2023near}
 by a $\sqrt{\kappa_\vy}$ factor due to the use of lower-level AGD. 



\begin{restatable}{thm}{thmubp}  \label{thm:ub-p1}
Under the setting of Lemma \ref{lem:F2BA}, Algorithm \ref{alg:F2BA} can find an
$\epsilon$-stationary point of $F(\vx)$ in 
$\tilde \gO ( \bar \kappa_\vy^{7/2} L_1 \Delta \epsilon^{-2}  )$
first-order oracle complexity. In particular, if the lower-level function  $g(\vx,\vy)$ is quadratic as in Eq.  (\ref{eq:quadratic-g}), the complexity can be refined to $\tilde \gO ( \kappa_\vy^{5/2} L_1 \Delta \epsilon^{-2}  )$.
\end{restatable}

\begin{proof}
Recall that $\vy_\lambda^*(\vx) = \arg \min_{\vy \in \sR^{d_y}} f(\vx,\vy) + \lambda g(\vx,\vy)$. For any $\lambda \ge 2 L_1/ \mu_\vy$, we know that sub-problem for variable $\vy$ is $2\lambda L_1$-smooth and $\lambda \mu_\vy/2$-strongly convex. Therefore, we can apply Lemma \ref{lem:AGD} with condition number $4 L_1 / \mu_\vy$, in conjunction with \ref{lem:F2BA}, to prove the upper complexity bound of $\tilde \gO(\sqrt{\kappa_\vy} L_F \Delta \epsilon^{-2})$, where $L_F$ is given by Lemma \ref{lem:smooth-phi}.
\end{proof}

The above theorem refines that in \citet{chen2023near} by a factor of $\sqrt{\kappa_\vy}$ in the general setting. {In addition, for NC--Q problems, we demonstrate a further reduction in the condition number dependency and show that it can be refined from $\bar \kappa_\vy$ to $\kappa_\vy$  due to a smaller smoothness constant when $\nabla^2 g(\vx,\vy)$ remains constant.}
Due to similar reasons, we can also improve the guarantee for second-order smooth problems in \citet{chen2023near} as we show.


\subsection{Upper Bounds for Second-order Smooth Problems}

Now, we consider second-order smooth NC--SC bilevel problems, where the hyper-objective $F(\vx)$ is proved to have Lipschitz continuous Hessians.


\begin{lem}[{\citet[Lemma 2.4]{huang2025efficiently}}] \label{lem:smooth-phi-second}
For a second-order smooth NC--SC bilevel problem $(f,g) \in \gF^{\text{nc--sc}}(L_0,\cdots,L_3,\mu_\vy,\Delta)$, the hyper-objective $F(\vx) = f(\vx,\vy^*(\vx))$ satisfies
\begin{align} \label{eq:hyper-smooth-second}
    \Vert \nabla^2 F(\vx) - \nabla^2 F(\vx') \Vert \le \rho_F \Vert \vx-  \vx' \Vert, \quad \forall \vx,\vx' \in \sR^{d_x},
\end{align}
where $\rho_F = \gO( \bar \kappa_\vy^5 L_2 )$ and $\bar \kappa_\vy = \max\{L_0,L_1,L_2,L_3 \} / \mu_\vy$. In particular, if the lower-level function is quadratic and $(f,g) \in \gF^{\text{nc--q}}(L_1,L_2,\mu_\vy, \Delta)$, we have $\rho_F = \gO(\bar \kappa_\vy^3 L_2)$.
\end{lem}

Following \citep{nesterov2006cubic}, we target to find stronger solutions known as second-order stationary points for Hessian smooth nonconvex functions. Compared with first-order stationary points in Definition \ref{dfn:fo-sta}, the second-order stationary points in the following Definition \ref{dfn:so-sta} exclude saddle points and can be equivalent to global minimizers in certain applications \citep{ge2015escaping}.

\begin{dfn} \label{dfn:so-sta}
We call 
a point $\vx \in \sR^d$ an $\epsilon$-second-order stationary point of $F$ if 
\begin{align*}
    \Vert \nabla F(\vx) \Vert \le \epsilon \quad \text{and} \quad \nabla^2 F(\vx) \succeq -\sqrt{\rho_{F} \epsilon} \mI_{d}.
\end{align*}
\end{dfn}

For second-order smooth problems, \citet{chen2023near} proposed AccF${}^2$BA to achieve a $\tilde \gO( \kappa_\vy \sqrt{L_F} \rho_F^{1/4} \epsilon^{-7/4} )$ upper bound by applying nonconvex AGD \citep{carmon2017convex,jin2018accelerated,li2023restarted} in the upper level.
But again, their lower-level GD sub-solver is suboptimal. Based on this observation, we propose 
AccF${}^2$BA${}^+$ in Algorithm \ref{alg:AccF2BA}, which replaces the GD sub-solver in AccF${}^2$BA with AGD.
By keeping the outer loop unchanged, we can directly apply the results from \citep{chen2023near} to obtain the following guarantee of the outer loop, provided that the lower-level sub-solver achieves sufficient accuracy.


\begin{algorithm*}[t]  
\caption{AccF${}^2$BA${}^+$} \label{alg:AccF2BA}
\begin{algorithmic}[1] 
\STATE $ \vz_0 = \vy_0$, $\vx_{-1} = \vx_0$ \\[1mm]
\STATE \textbf{while} $ t < T $ \\[1mm]
\STATE \quad $\tilde \vx_{t} = \vx_t + (1-\theta) (\vx_t -\vx_{t-1})$ \\[1mm]
\STATE \quad $\vz_{t+1} = \text{AGD}( g(\tilde \vx_t,\,\cdot\,), K_t, \vz_t)$ \\[1mm]
\STATE \quad $\vy_{t+1} = \text{AGD} (f(\tilde \vx_t,\,\cdot\,) + \lambda g(\tilde \vx_t, \,\cdot\,) , K_t, \vy_t    )$ \\[1mm]
\STATE \quad $ \hat \nabla F_\lambda(\tilde \vx_{t})= \nabla_\vx f(\tilde \vx_{t},\vy_{t+1}) + \lambda ( \nabla_\vx g(\tilde \vx_{t},\vy_{t+1}) - \nabla_\vx g(\tilde \vx_{t},\vz_{t+1}) )$ \\[1mm]
\STATE\quad  $\vx_{t+1} = \tilde \vx_{t} - \eta_x \hat \nabla F_\lambda(\tilde \vx_{t}) $ \\[1mm]
\STATE \quad $ t = t+1 $ \\ [1mm]
\STATE \quad \textbf{if} $ t \sum_{j=0}^{t-1} \Vert \vx_{j+1} - \vx_j \Vert^2 \ge B^2$ \\[1mm]
\STATE \quad \quad $t=0$, $\vx_{-1} = \vx_0 = \vx_t + \xi_t \vone_{ \Vert \hat \nabla F_\lambda(\tilde \vx_{t}) \Vert \le \frac{B}{2 \eta_x} }$, {\rm where} $\xi_t \sim \sB(r)$  \\[1mm] 
\STATE \quad \textbf{end if} \\[1mm]
\STATE \textbf{end while} \\[1mm]
\STATE $T_0 = \arg \min_{\lfloor \frac{T}{2} \rfloor \le t \le T-1} \Vert \vx_{t+1} - \vx_t \Vert$ \\[1mm]
\STATE \textbf{return} $\vx_{\rm out} = \frac{1}{T_0+1} \sum_{t=0}^{T_0} \tilde \vx_{t}$ \\[1mm]
\end{algorithmic}
\end{algorithm*}

\begin{lem}{\citep[Thm. 5.2]{chen2023near}, \citep[Thm. 4.1]{yang2023accelerating}} \label{lem:AccF2BA}
Let $\sqrt{\epsilon \rho_F} \le L_F$.
For a second-order smooth NC--SC bilevel problem $(f,g) \in \gF^{\text{nc--sc}}(L_0,L_1,L_2,L_3,\mu_\vy,\Delta)$, 
if in each iteration $t$ of Algorithm~\ref{alg:AccF2BA} the lower-level solver outputs 
$\vz_{t+1}$ and $\vy_{t+1}$ such that 
\begin{align*}
    \Vert \vz_{t+1} - \vy^*(\vx_t) \Vert^2+ \Vert \vy_{t+1} - \vy_{\lambda}^*(\vx_t) \Vert^2 \le \gamma \Vert \vz_t - \vy^*(\vx_t) \Vert^2+ \Vert \vy_t - \vy_{\lambda}^*(\vx_t) \Vert^2
\end{align*}
with $\gamma \lesssim {\rm poly}(\epsilon) / {\rm poly}(\log d_x, \bar \kappa_\vy, L_1) $, then setting the hyper-parameters as
\begin{align*}
    \lambda \asymp \max\{ L_1 \bar \kappa_\vy^3 /\epsilon, ~\bar \kappa_\vy^{7/2} \sqrt{L_2/\epsilon} \}, ~~ \eta_x \asymp 1/L_F, ~~ \chi \asymp \log(\nicefrac{d_x}{\delta \epsilon}),  \\
    \theta \asymp \sqrt{\eta_x} ({\rho_F} \epsilon)^{1/4}, ~~ T \asymp \chi / \theta,  ~~ B \asymp \sqrt{\epsilon/\rho_F} / \chi^2, ~~ {\text{and}} ~~ r\asymp \epsilon 
\end{align*}
ensures the algorithm to output an $\epsilon$-second-order stationary point with probability $1-\delta$
in \[
\gO( \sqrt{L_{F}} \rho_{F}^{1/4} \Delta \epsilon^{-7/4} \log^6(\nicefrac{d_x}{\delta \epsilon}))\]
steps of update in the upper-level variable $\vx$.
\end{lem}

Now, by combining Lemma \ref{lem:AGD}, \ref{lem:smooth-phi-second}, and \ref{lem:AccF2BA}, we naturally arrive at the following theorem for the total complexity of AccF${}^2$BA${}^+$. 

\begin{restatable}{thm}{thmubpacc}  \label{thm:ub-p2}
Under the same setting of Lemma \ref{lem:AccF2BA}, Algorithm \ref{alg:AccF2BA} can find an
$\epsilon$-second-order stationary point of $F(\vx)$ in 
$\tilde \gO ( \bar \kappa_\vy^{13/4} \sqrt{L_1} L_2^{1/4} \Delta \epsilon^{-7/4}  )$
first-order oracle complexity. In particular, if the lower-level function $g(\vx,\vy)$ is quadratic as in Eq.  (\ref{eq:quadratic-g}), the complexity can be refined to $\tilde \gO ( \bar \kappa_\vy^{9/4} \sqrt{L_1} L_2^{1/4} \Delta \epsilon^{-7/4} )$.
\end{restatable}

\begin{proof}
By combining Lemma \ref{lem:AGD} and \ref{lem:AccF2BA}, we know that the total complexity of Algorithm~\ref{alg:AccF2BA} is \[
\gO( \sqrt{\kappa_\vy L_F} \rho_F^{1/4} \Delta \epsilon^{-7/4} \log^6 (\nicefrac{d_x}{\delta \epsilon})  ).
\]
Then we can obtain the claimed complexity by substituting the value of $L_F = \gO(\bar \kappa_\vy^3 L_1)$, $\rho_F = \gO(\bar \kappa_\vy^5 L_2)$ for general lower-level problems and $L_F = \gO( \kappa_\vy^2 L_1 )$, $\rho_F = \gO(\bar \kappa_\vy^3 L_2)$ according to Lemma \ref{lem:smooth-phi} and \ref{lem:smooth-phi-second}.
\end{proof}


In the following section, we also provide numerical experiments in Appendix \ref{apx:exp} to verify the effectiveness of replacing GD with AGD in the inner loop. However, we remark that the use of lower-level AGD for reducing condition number dependency is limited to the deterministic setting, perhaps because the bottleneck of NC--SC stochastic bilevel optimization lies in the variance instead of the bias.

\subsection{Numerical Experiments} \label{apx:exp}

In this section, we conduct experiments to verify the benefit of using AGD to solve the lower-level problems. We first conduct a synthetic experiment where the condition number~$\kappa_\vy$ can be explicitly controlled to validate the improved $\kappa_\vy$ dependence, and we then provide numerical experiments on a practical scenario.

{
\paragraph{Hyperparameter Setting.} For each run, we tune the penalty and stepsizes hyperparameters via a base-10 log-scale grid search (\textit{e.g.}, in the real-world experiment, we tune $\lambda \in \{10^1, 10^2, 10^3\}$, $\eta_x,\eta_y \in \{10^{-3}, 10^{-2}, 10^{-1}, 10^0, 10^1, 10^2, 10^3\}$) and report the best performance, and directly choose $\theta = 0.1$ and $B=0$ in AccF${}^2$BA and AccF${}^2$BA${}^+$.}

{\subsubsection{Synthetic Experiments}

We first consider a synthetic bilevel optimization problem of the form:
\begin{align*}
    \min_{\vx \in \sR^m} \frac{1}{2} \Vert \vy \Vert^2 \quad {\rm s.t.} \quad \vy \in \arg \min_{\vy \in \sR^n} \frac{1}{2} \vy^\top \mH \vy+ \vx^\top \mJ \vy,
\end{align*}
where $\mathbf{H} \in \mathbb{R}^{n \times n}$ is a diagonal matrix with entries uniformly distributed in $[1/\kappa_\vy, 1]$, $\mathbf{J}$ is a random matrix with a spectral norm of $1$, and both $\mathbf{x}_0$ and $\mathbf{y}_0$ are initialized as standard Gaussian vectors. We compare the gradient-based methods F${}^2$BA and AccF${}^2$BA \citep{chen2023near}, and our improved versions with lower-level AGD denoted by F${}^2$BA${}^+$ and AccF${}^2$BA${}^+$ in terms of the number of total iterations to reach an $\epsilon$-stationary point in two experiment sets in Figure \ref{fig:synthetic}: the first one with $\epsilon = 0.01$ and $\kappa_\vy$ varying from $[1,10^2]$, the second one with $\epsilon =0.1$ and $\kappa_\vy$ varying from $[1,10^3]$.}

\begin{figure}[htbp]
    \centering
    \includegraphics[scale=0.26]{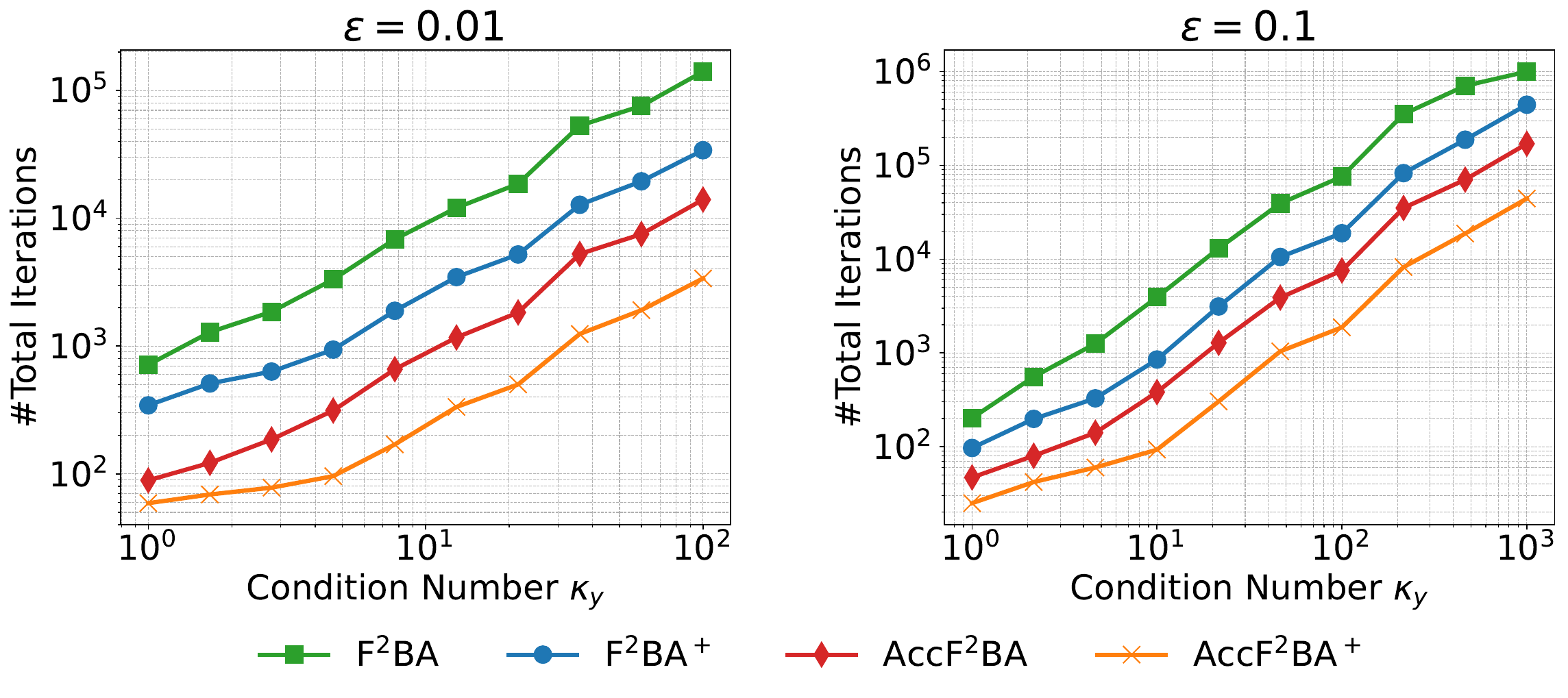} 
\caption{Performances of different algorithms on the synthetic experiments.}
\label{fig:synthetic}
\end{figure}

\subsubsection{Real-World Experiments}
To further verify the effectiveness of using lower-level AGD, we then consider the ``learn-to-regularize'' problem 
\citep{grazzi2020iteration,ji2021bilevel,liu2022bome,chen2023near}, on the ``20 Newsgroup'' dataset, which solves the following bilevel optimization problem:
\begin{align*}
    \min_{\vx \in \sR^{p}} \frac{1}{m} \sum_{i=1}^m \ell_i^{\rm val}(\mY), \quad {\rm s.t.} \quad \mY \in \arg \min_{\mY \in \sR^{p \times q}} \frac{1}{n} \sum_{i=1}^n \ell_i^{\rm tr}(\mY) + \Vert \mW_x \mY  \Vert^2_F ,
\end{align*}
where $\mW_\vx = {\rm diag} \left(\sqrt{\exp(\vx)}\right)$, 
$\vy \in \sR^{p \times q}$ is the parameter matrix of a multi-class logistic model for $p$ features and $q$ classes, $\{ \ell_i^{\rm val}\}_{i=1}^m$ and $\{ \ell_i^{\rm tr}\}_{i=1}^n$ are the validation and training loss on every sample each sample with $m,n$ being the size of the datasets.
In our experiments, we have $p= 130, 107$, $q = 20$, and $m = n = 5,657$. We compare six algorithms, including a baseline without learning the regularization ``w/o Reg'', the HVP-based method AID \citep{ji2021bilevel}, two gradient-based methods F${}^2$BA and AccF${}^2$BA \citep{chen2023near}, and our improved versions with lower-level AGD denoted by F${}^2$BA${}^+$ and AccF${}^2$BA${}^+$.
We report the optimal performance of each algorithm after being well-tuned in Figure \ref{fig:l2reg}.

\begin{figure}[htbp]
    \centering
    \includegraphics[scale=0.26]{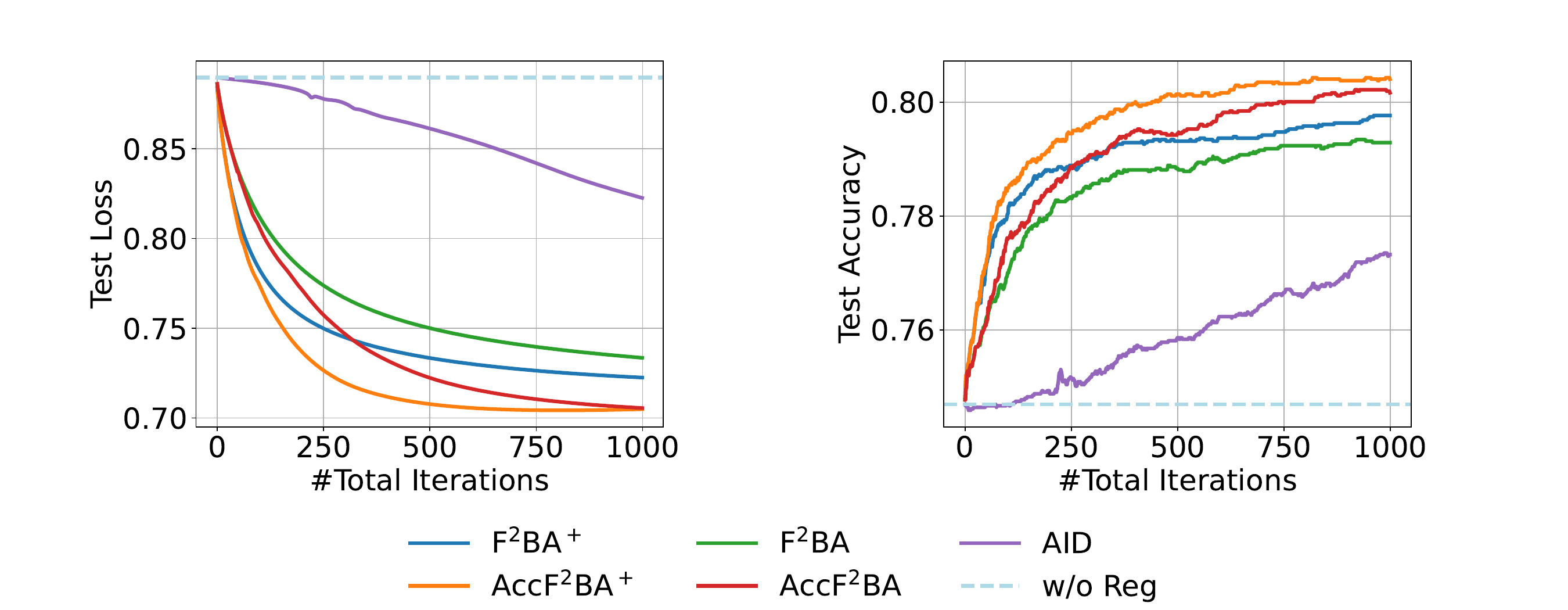} 
\caption{Performances of different algorithms on learning the optimal regularization of the ``20 Newsgroup'' dataset.}
\label{fig:l2reg}
\end{figure}

\newpage

\end{document}